\documentclass{conm-p-l}
\usepackage{amssymb, amsmath, verbatim}

\newtheorem{theorem}{Theorem}[subsection]

\newtheorem{prop}[theorem]{Proposition}
\newtheorem{cor}[theorem]{Corollary}

\theoremstyle{definition}
\newtheorem{definition}[theorem]{Definition}

\theoremstyle{remark}
\newtheorem{remark}[theorem]{Remark}

\numberwithin{equation}{section}

\copyrightinfo{2000}{American Mathematical Society}
\def\<{\left<}
\def\>{\right>}
\def\cstar{$C^*$-algebra}
\def\esg{$E_0$-semigroup}
\def\ind{{\rm{index}}}
\def\tr{{\rm{trace}}}

\begin{document}

\title{Four Lectures on Noncommutative Dynamics}

\author{William Arveson}
\address{Department of Mathematics, University of California, 
Berkeley, CA, 94720}
\email{arveson@mail.math.berkeley.edu}
\thanks{Supported by NSF grant DMS-0100487.}

\date{January 1, 1994 and, in revised form, June 22, 1994.}



\begin{abstract}
These lectures concern basic aspects of the theory 
of semigroups of endomorphisms of type $I$ factors 
that relate to causal dynamics, 
dilation theory, and the problem of classifying 
$E_0$-semigroups 
up to cocycle conjugacy.  
We give only a few proofs here; full details 
can be found in the author's upcoming monograph 
{\em Noncommutative Dynamics and $E$-semigroups\/}, 
to be published in the Springer series 
Monographs in Mathematics.  
\end{abstract}

\maketitle
\tableofcontents

\specialsection{Dynamical Origins: Histories and Interactions}

The dynamics of one-parameter groups of 
automorphisms of $\mathcal B(H)$ 
has been completely understood since the 
1930s, following work of Wigner, Stone, and 
the multiplicity theory of Hahn and Hellinger.  After 
reviewing these basic issues, 
we show that if one takes into account a natural notion 
of causality for such dynamical groups,
one encounters entirely new phenomena.  We describe 
the basic properties of these ``causal" dynamical systems and 
their connections with $E_0$-semigroups, 
we discuss positive results concerning their existence 
and nontriviality, and we point out some of 
the basic mathematical problems that remain unsolved.  

\subsection{The Dynamics of Quantum Systems}

Let us recall the basic mathematical setting of quantum physics.  
The observables of quantum theory are self-adjoint operators 
acting on a separable Hilbert space $H$.  Observables such 
as linear or angular momentum 
arise as generators of one-parameter unitary groups and are typically 
unbounded and only densely defined.  However, there is no essential loss 
in restricting attention to bounded functions of these unbounded 
operators.  It goes without saying 
that one has to make use of the spectral theorem in order 
to define bounded functions of unbounded self-adjoint 
operators, and we do so freely.  

The quantum replacement for the distribution of a random variable requires 
specifying not only an observable $X$ but also a unit vector $\xi\in H$.  
Once $X$ and $\xi$ are fixed, there is a unique probability measure 
$\mu_{X,\xi}$ defined on the real line by specifying its integral 
with respect to bounded Borel functions $f: \mathbb R\to\mathbb R$ 
as follows
$$
\int_{-\infty}^\infty f(t)\,d\mu_{X,\xi}(t) = \<f(X)\xi,\xi\>.  
$$  
For a Borel set $S\subseteq\mathbb R$, one interprets 
$\mu_{X,\xi}(S)$ as representing the probability of finding an observed 
value of $X$ in the set $S$, given that the system is in the 
pure state associated with $\xi$. 

Given two observables $X,Y$ and a unit vector $\xi$, there is 
no ``joint distribution" $\mu_{X,Y,\xi}$ defined on $\mathbb R^2$.  
On the level of physics, this phenomenon is associated with the theory 
of measurement and is a consequence of the uncertainty principle.  
From the point of view of operator theory, 
since the operators $X$ and $Y$ normally fail to commute there is no 
way of using them to define a spectral measure on $\mathbb R^2$.  This 
nonexistence of joint distributions is one of the fundamental 
differences between quantum theory and probability theory.  

Turning now to dynamics, consider the way the flow of time 
acts on the algebra of observables.  Every symmetry of quantum theory 
corresponds to either a $*$-automorphism or a $*$-anti-automorphism 
of the algebra $\mathcal B(H)$ of all bounded operators on $H$.  If 
we are given a one-parameter group of such symmetries, then since each 
one of them is the square of another it follows that 
all of the symmetries must be $*$-automorphisms.  
Thus, the flow of time on a quantum sytem is given by a one-parameter 
family $\alpha=\{\alpha_t: t\in\mathbb R\}$ of automorphisms of 
$\mathcal B(H)$ such that $\alpha_s\circ\alpha_t=\alpha_{s+t}$, and 
which satisfies the natural continuity condition: 
for every $A\in\mathcal B(H)$ and every pair of vectors $\xi,\eta\in H$ 
the function $t\in\mathbb R\mapsto \<\alpha_t(A)\xi,\eta\>$ is 
continuous.

Let us consider the possibilities: How does one  
classify one-parameter groups of automorphisms of $\mathcal B(H)$?  
In the late 1930s, Eugene Wigner proved that
every such  one-parameter group is implemented by a strongly continuous 
one-parameter unitary group $U=\{U_t: t\in\mathbb R\}$ in the 
sense that 
$$
\alpha_t(A)=U_tAU_t^*, \qquad A\in\mathcal B(H), \quad t\in R.  
$$

Earlier, Marshall Stone had shown that a strongly continuous 
one-parameter unitary group $U$ is the Fourier transform of a 
unique spectral measure $E$ defined on the Borel subsets of 
the real line as follows
$$
U_t=\int_{-\infty}^\infty e^{i\lambda t}\,dE(\lambda).  
$$
Equivalently, Stone's theorem implies that for the unbounded 
self-adjoint operator $X=\int_{\mathbb R}\lambda\,dE(\lambda)$, 
we have $U_t=e^{itX}$, $t\in\mathbb R$.  

Together, these two results imply that 
every one-parameter group $\alpha$ 
of automorphisms of $\mathcal B(H)$ corresponds to an observable 
$X$ as follows
$$
\alpha_t(A)=e^{itX}Ae^{-itX}, \qquad A\in\mathcal B(H),\quad t\in \mathbb R.  
$$
$X$ is not uniquely determined by $\alpha$, since replacing 
$X$ with a scalar translate of the form $X+\lambda\bf 1$ 
with $\lambda\in\mathbb R$ does not 
change $\alpha$.  However, $X$ is uniquely determined 
by $\alpha$ up to such scalar perturbations.

Two one-parameter groups $\alpha$ 
and $\beta$  of $*$-automorphisms
(acting on $\mathcal B(H)$ and $\mathcal B(K)$ 
respectively) are 
said to be {\it conjugate} if there is a $*$-isomorphism 
$\theta$ of $\mathcal B(H)\to\mathcal B(K)$ such that 
$$
\theta(\alpha_t(A)) = \beta_t(\theta(A)), \qquad 
A\in\mathcal B(H),  \quad t\in \mathbb R.  
$$
Recalling that such a $*$-isomorphism $\theta$ must 
be implemented by a unitary operator $W: H\to K$ by way 
of $\theta(A)=WAW^*$, we see that 
Wigner's theorem completely settles the classification issue 
for one-parameter groups of automorphisms of $\mathcal B(H)$.  
Indeed, using that result we may find unbounded self-adjoint 
operators $X$, $Y$ on the respective Hilbert spaces such that 
$\alpha$ and $\beta$ are given by 
$\alpha_t(A)=e^{itX}Ae^{-itX}$ and $\beta_t(B)=e^{itY}Be^{-itY}$, 
$A\in\mathcal B(H)$, $B\in\mathcal B(K)$, $t\in\mathbb R$.   
It is an elementary exercise to 
show that $\theta(A)=WAW^*$ implements 
a conjugacy of $\alpha$ and $\beta$ if and only if 
there is a real scalar $\lambda$ 
such that $WXW^*=Y+\lambda\mathbf 1$.  Thus, the 
classification of one-parameter groups of automorphisms 
is reduced to the classification 
of unbounded self-adjoint operators up to unitary equivalence.  
By the spectral theorem, this is equivalent to the classification 
up to unitary equivalence 
of spectral measures on the real line; and the latter 
problem is completely understood in terms of the 
multiplicity theory of Hahn and Hellinger.  

These remarks show that the most basic aspect of quantum 
dynamics, namely the structure and classification of one-paramter 
groups of $*$-automorphisms of $\mathcal B(H)$, is completely understood.  
We have seen all of the possibilities, 
and they are described by self-adjoint 
operators (or spectral measures on the line) and 
their multiplicity theory in a completely explicit way.

\subsection{Causality: Histories and Interactions}

We now show that
by  introducing a natural notion of causality into such dynamical 
systems, one encounters entirely new phenomena.  These ``causal" 
dynamical systems acting on $\mathcal B(H)$ are only partially
understood.  We have surely not seen all of them, and 
we have only partial information about how to classify  
the ones we have seen.

We are concerned with one-parameter 
groups of $*$-automorphisms, of the algebra $\mathcal B(H)$ 
of all bounded operators on a 
Hilbert space $H$, which carry a particular kind 
of causal structure.  More precisely, 
A {\it history} is a pair $(U, M)$ consisting of 
a one-parameter group $U=\{U_t: t\in\mathbb R\}$ of unitary
operators acting on a separable infinite-dimensional
Hilbert space $H$, together with a type $I$ subfactor
$M\subseteq \mathcal B(H)$ that is invariant under 
the automorphisms 
$\gamma_t(X)=U_tXU_t^*$ for negative $t$, 
and which has the following two properties
\begin{enumerate}
\item[(i)]\label{irred} (irreducibility)
$$
(\bigcup_{t\in \mathbb R}\gamma_t(M))^{\prime\prime} = \mathcal B(H),
$$
\item[(ii)]\label{remPast} (trivial infinitely remote past)
$$
\bigcap_{t\in\mathbb R}\gamma_t(M) = \mathbb C\cdot\mathbf 1.  
$$
\end{enumerate}
We find it useful to think of the group 
$\{\gamma_t: t\in\mathbb R\}$ as representing the flow
of time in the Heisenberg picture, 
and the von Neumann algebra $M$ as representing
bounded observables that are associated with the ``past".  
However, we will be concerned with 
purely mathematical issues  
concerning the dynamical 
properties of histories, with problems concerning 
their existence and construction, and especially 
with the issue of nontriviality (to be defined momentarily).  
Two histories $(U,M)$ (acting on $H$) and 
$(\tilde U,\tilde M)$ (acting on $\tilde H$) 
are said to be {\it isomorphic} 
if there is a 
$*$-isomorphism $\theta: \mathcal B(H)\to\mathcal B(K)$
such that $\theta(M)=\tilde M$ and 
$\theta\circ\gamma_t=\tilde\gamma_t\circ\theta$ 
for every $t\in\mathbb R$, 
$\gamma$, $\tilde\gamma$ denoting the automorphism 
groups associated with $U$, $\tilde U$.  The basic 
problems addressed in these lectures all bear some relation 
to the problem of classifying histories.  We have 
already alluded to the fact that the results are 
far from complete.  

An {\it $E_0$-semigroup} is a one-parameter
semigroup $\alpha=\{\alpha_t: t\geq 0\}$ of unit-preserving 
$*$-endomorphisms of a type $I_\infty$ factor $M$, which 
is continuous in the natural sense.  
The subfactors $\alpha_t(M)$ decrease as 
$t$ increases, and $\alpha$ is called 
{\it pure} if $\cap_t\alpha_t(M)=\mathbb C\mathbf 1$.  
There are two $E_0$-semigroups $\alpha^-$, 
$\alpha^+$ associated with any history,
$\alpha^-$ being the one associated with 
the ``past" by restricting $\gamma_{-t}$
to $M$ for $t\geq 0$ and $\alpha^+$ being
the one associated
with the ``future" by restricting $\gamma_t$ 
to the commutant $M^\prime$ for $t\geq 0$.  

By an {\it interaction} we mean a history with 
the additional property that there are normal
states $\omega_-$, $\omega_+$ of $M$, $M^\prime$
respectively such that $\omega_-$ is invariant 
under the action of $\alpha^-$ 
and $\omega_+$ is invariant under the action 
of $\alpha^+$.  Because of 
(i) and (ii), both $\alpha^-$
and $\alpha^+$ are pure $E_0$-semigroups, and 
it is not hard to show that 
if a pure $E_0$-semigroup $\alpha$ has a normal invariant 
state $\omega$ then that state is an {\it absorbing} state 
in the sense that for every other normal state 
$\rho$ on the domain of $\alpha$ one 
has 
$$
\lim_{t\to\infty}\|\rho\circ\alpha_t-\omega\|=0.
$$  
In particular, $\omega_-$ (resp. $\omega_+$) is the 
unique normal state of $M$ (resp. $M^\prime$) 
that is invariant under the action of $\alpha^-$
(resp. $\alpha^+$).  

In particular, it 
follows from this uniqueness that if one 
is given two interactions $(U,M)$ and 
$(\tilde U, \tilde M)$ with respective 
pairs of normal states $\omega_+$, 
$\omega_-$ and $\tilde\omega_+$, 
$\tilde\omega_-$, then an isomorphism 
of histories $\theta: (U,M)\to (\tilde U,\tilde M)$ 
must associate $\tilde\omega_+$, $\tilde\omega_-$ 
with $\omega_+$, $\omega_-$ in the sense that 
if $\theta_+$ (resp. $\theta_-$) denotes 
the restriction of $\theta$ to $M$ 
(resp. $M^\prime$), then one has 
$\tilde\omega_\pm\circ\theta_\pm=\omega_\pm$.  

\begin{remark}[Normal Invariant States]\label{invariantStates}
Since the state space of a unital \cstar\ is
weak$^*$-compact, 
the Markov-Kakutani fixed 
point theorem implies that every $E_0$-semigroup
has invariant states.  But there is no reason to 
expect that there is a {\it normal} invariant 
state.  Indeed, there are examples 
of pure $E_0$-semigroups
which have no normal invariant states.  
Notice too that $\omega_-$, for example, is defined {\it only}
on the algebra $M$ of the past.  Of course, $\omega^-$ has 
many extensions to normal states of $\mathcal B(H)$, but 
none of these normal extensions 
need be invariant under the action of 
the group $\gamma$.  In fact, we will see that if there 
is a {\it normal} $\gamma$-invariant state defined 
on all of $\mathcal B(H)$ then the interaction must 
be trivial.  
\end{remark}

In order to discuss the dynamics of interactions we must 
introduce a \cstar\ of ``local observables".  For 
every 
compact interval $[s,t]\subseteq \mathbb R$ 
there is an associated von Neumann algebra
\begin{equation}\label{introIntEq0.3}
\mathcal A_{[s,t]} = \gamma_t(M)\cap \gamma_s(M)^\prime.  
\end{equation}
Notice that since $\gamma_s(M)\subseteq\gamma_t(M)$ are 
both type $I$ factors, so is the relative commutant
$\mathcal A_{[s,t]}$.
Clearly $\mathcal A_I\subseteq \mathcal A_J$ if 
$I\subseteq J$, and for adjacent intervals 
$[r,s]$, $[s,t]$, $r\leq s\leq t$ we have 
\begin{equation}\label{introIntEq0.4}
\mathcal A_{[r,t]} = \mathcal A_{[r,s]}\otimes\mathcal A_{[s,t]},
\end{equation}
in the sense that the two factors $\mathcal A_{[r,s]}$ and 
$\mathcal A_{[s,t]}$ mutually commute and generate 
$\mathcal A_{[r,t]}$ as a von Neumann algebra.  The 
automorphism group $\gamma$ permutes the algebras
$\mathcal A_I$ covariantly,
\begin{equation}\label{introIntEq0.5}
\gamma_t(\mathcal A_I) = \mathcal A_{I+t}, \qquad t\in\mathbb R. 
\end{equation}
Finally, we define the local \cstar\ $\mathcal A$
to be the {\it norm} closure of the union of all 
the $\mathcal A_I$, $I\subseteq \mathbb R$.  
$\mathcal A$ is a $C^*$-subalgebra
of $\mathcal B(H)$ that is strongly dense and 
invariant under the 
action of the automorphism group $\gamma$.

\begin{remark}
It may be of interest to compare the local structure
of the \cstar\ $\mathcal A$ to its commutative 
counterpart, namely the local algebras associated 
with a stationary random distribution with independent
values at every point \cite{GelVil}.  More precisely, suppose
that we are given a random distribution $\phi$; i.e., 
a linear map from the space of real-valued test functions on 
$\mathbb R$ to the space of real-valued random variables on 
some probability space $(\Omega, P)$.  With every 
compact interval $I=[s,t]$ with $s<t$ one may consider 
the weak$^*$-closed subalgebra $\mathcal A_I$ of 
$L^\infty(\Omega,P)$ generated by random variables of the form
$e^{i\phi(f)}$, $f$ ranging over all test functions supported 
in $I$.  When the random distribution $\phi$ is stationary 
and has independent values at every point, this family 
of subalgebras of $L^\infty(\Omega,P)$ 
has properties analogous to (\ref{introIntEq0.4}) and 
(\ref{introIntEq0.5}), in that 
there is a one-parameter group of measure preserving 
automorphisms $\gamma=\{\gamma_t: t\in \mathbb R\}$ 
of $L^\infty(\Omega,P)$ which satisfies (\ref{introIntEq0.5}), and 
instead of (\ref{introIntEq0.4}) we have the assertion that the algebras 
$\mathcal A_{[r,s]}$ and $\mathcal A_{[s,t]}$ are 
{\it probabilistically independent} and generate 
$\mathcal A_{[r,t]}$ as a weak$^*$-closed algebra.  

One should keep in mind, however, that this commutative 
analogy has serious limitations.  For example, we have 
already pointed out that 
in the case of interactions there is 
typically no normal 
$\gamma$-invariant state on $\mathcal B(H)$,
and there is no reason to expect 
any normal state of $\mathcal B(H)$ to 
decompose as a product state relative to the 
decompositions of (\ref{introIntEq0.4}).  

There is also some common 
ground with the Boolean algebras of type $I$ factors
of Araki and Woods \cite{arakiWoods}, but here too there are 
significant differences.  For example, the local algebras
of (\ref{introIntEq0.3}) and (\ref{introIntEq0.4}) 
are associated with intervals (and 
more generally with finite unions of intervals), but 
not with more general Borel sets as in \cite{arakiWoods}. 
Moreover, here the translation group acts as automorphisms
of the given structure whereas in \cite{arakiWoods} there is 
no hypothesis of symmetry with respect
to translations.  
\end{remark}

\subsection{Dynamics of Interactions}
The \cstar\ $\mathcal A$ 
of local observables is important because it 
provides a way of comparing $\omega_-$ 
and $\omega_+$.  
Indeed, both states $\omega_-$ and $\omega_+$ extend 
{\it uniquely} to $\gamma$-invariant states 
$\bar\omega_-$ and $\bar\omega_+$ of $\mathcal A$.  
We sketch the proof for $\omega_-$.  

\begin{prop}\label{introIntInvState}
There is a unique $\gamma$-invariant state 
$\bar\omega_-$ of $\mathcal A$ such that 
$$
\bar\omega_-\restriction_{\mathcal A_I}=
\omega_-\restriction_{\mathcal A_I}
$$
for every compact interval $I\subseteq (-\infty,0]$.  
\end{prop}
\begin{proof}
For existence of the extension, choose any compact interval
$I=[a,b]$ and any operator $X\in\mathcal A_I$.  Then for sufficiently
large $s>0$ we have $I-s\subseteq (-\infty,0]$ and for these 
values of $s$ $\omega_-(\gamma_{-s}(X))$ does not depend on 
$s$ because $\omega_-$ is invariant under the 
action of $\{\gamma_t: t\leq 0\}$.  
Thus we can define $\bar\omega_-(X)$ unambiguously by
$$
\bar\omega_-(X)=\lim_{t\to-\infty}\omega_-(\gamma_t(X)).  
$$
This defines a positive linear functional $\bar\omega_-$ on 
the unital $*$-algebra
$\cup_I\mathcal A_I$, and now we extend $\bar\omega_-$ to 
all of $\mathcal A$ be norm-continuity.  The extended state is 
clearly invariant under the action of $\gamma_t$, $t\in\mathbb R$.  

The uniqueness of the state $\bar\omega_-$ is apparent.  
\end{proof}

It is clear from the proof of Proposition \ref{introIntInvState} that 
these extensions of $\omega_-$ and 
$\omega_+$ are {\it locally normal} in the sense that 
their restrictions to any localized subalgebra
$\mathcal A_I$ define normal states on that type $I$ factor.  
Thus, the local \cstar\  $\mathcal A$ has a definite
``state of the past" and a definite ``state of the future" 
in the following sense:  

\begin{prop}
For every $X\in \mathcal A$ and every normal state $\rho$ of 
$\mathcal B(H)$ we have 
$$
\lim_{t\to-\infty}\rho(\gamma_t(X))=\bar\omega_-(X), 
\qquad
\lim_{t\to+\infty}\rho(\gamma_t(X))=\bar\omega_+(X)
$$
\end{prop}
\begin{proof}
Consider the first limit formula.  
The set of all $X\in\mathcal A$ for which this formula holds
is clearly closed in the operator norm, hence it suffices 
to show that it contains $\mathcal A_I$ for every compact 
interval $I\subseteq\mathbb R$.  

We will make use of the fact (discussed more fully 
at the beginning of section 5) that if $\rho$ is 
any normal state of $M$ and $A$ is an operator in 
$M$ then 
$$
\lim_{t\to-\infty}\rho(\gamma_t(A))=\omega_-(A),
$$
see formula (4.1).  
Choosing a real number $T$ sufficiently negative that 
$I+T\subseteq (-\infty,0]$, the preceding remark shows
that for the operator $A=\gamma_T(X)\in M$ we have 
$\lim_{t\to-\infty}\rho(\gamma_t(A))=\omega_-(A)$, 
and hence
$$
\lim_{t\to-\infty}\rho(\gamma_t(X))=
\lim_{t\to-\infty}\rho(\gamma_{t-T}(\gamma_T(X)))=
\omega_-(\gamma_T(X))
=\bar\omega_-(X). 
$$
The proof of the second limit formula is similar.  
\end{proof}

\begin{definition}
The interaction $(U,M)$, with past and future states 
$\omega_-$ and $\omega_+$, is said to be trivial 
if $\bar\omega_-=\bar\omega_+$.  
\end{definition}

More generally, the norm $\|\bar\omega_- -\bar\omega_+\|$
gives some measure of the ``strength" of the 
interaction, and of course we have
$0\leq \|\bar\omega_- -\bar\omega_+\| \leq 2$.  

{\it If} there 
is a normal state $\rho$ of $\mathcal B(H)$ that is invariant under 
the action of $\gamma$, then since $\omega_-$ (resp. 
$\omega_+$) is the unique normal invariant state of 
$\alpha_-$ (resp. $\alpha_+$) we must have 
$\rho\restriction_M=\omega_-$, 
$\rho\restriction_{M^\prime}=\omega_+$, and hence 
$\bar\omega_-=\bar\omega_+=\rho\restriction_\mathcal A$ 
by the uniqueness part of Proposition \ref{introIntInvState}.     
In particular, 
{\it if the interaction is nontrivial then 
neither $\bar\omega_-$ nor $\bar\omega_+$ can be extended
from $\mathcal A$ to a normal state of its strong closure $\mathcal B(H)$}.   

Thus, whatever (normal) state $\rho$ one chooses to watch evolve
over time on operators in $\mathcal A$, it settles down to become 
$\bar\omega_+$ in the distant future, it must have 
come from $\bar\omega_-$ in the remote past, and the 
limit states do not 
depend on the choice of $\rho$.  For a trivial 
interaction, nothing happens over the long term: for 
fixed $X$ and $\rho$ 
the function $t\in\mathbb R\mapsto \rho(\gamma_t(X))$
starts out very near some value 
(namely $\bar\omega_-(X)$), exhibits transient fluctuations over 
some period of time, and then 
settles down near the same value again.  
For a nontrivial interaction, there will be a definite change 
from the limit at $-\infty$ to the limit at $+\infty$ - 
at least for some choices of $X\in\mathcal A$.

\begin{remark}[Existence of Nontrivial Interactions]
We have seen that every interaction gives rise to a 
pair of pure \esg s $\alpha^-$, $\alpha^+$, representing 
its ``past" and its ``future".  However, we have seen 
no examples of interactions, and in particular, we have 
said nothing to indicate that nontrivial 
interactions exist.  In the remainder of this lecture, we 
describe a body of results that address this key issue 
of existence in the simplest cases, namely where both 
past and future \esg s are cocycle perturbations of 
the CAR/CCR flows.  

In order to exhibit examples of interactions with 
such properties, one has 
to address three questions.  First, how does one 
construct examples of cocycle perturbations of 
CAR/CCR flows which are a) pure, and b) have a normal 
invariant state with specified properties?  
Second, given such a pair of 
\esg s, how does one determine when they can be 
``assembled" into an interaction so that 
one represents its future and the other represents 
its past?  Third, given that there is an interaction 
assembled from a pair of \esg s $\alpha^-$ and $\alpha^+$ 
in this way, 
how does one determine if 
that interaction is nontrivial?  
In sections 1.4, 1.5, 1.6 we discuss these three questions 
in turn.  
\end{remark}

\subsection{Cocycle Perturbations of CAR/CCR flows}\label{SS:pertCCR}
By a {\it cocycle} for an \esg\ $\alpha$ acting on 
a type $I$ factor $M$ 
we mean a {\it strongly continuous} 
family of unitary
operators $U = \{U_t: t\geq 0\ \}$ in $M$ satisfying 
\begin{equation}\label{basicsCoc2.1}
U_{s+t} = U_s\alpha_s(U_t),\qquad s,t \geq 0. 
\end{equation}
Notice that (\ref{basicsCoc2.1}) 
implies that $U_0 = \mathbf 1$.  The
condition  also implies that the family of endomorphisms 
$\beta = \{\beta_t: t\geq 0\ \}$ defined by
\begin{equation}\label{basicsCoc2.2}
\beta_t(A) = U_t \alpha_t(A) U_t^*,\qquad t\geq 0 
\end{equation}
satisfies the semigroup property 
$\beta_{s+t}=\beta_s\circ\beta_t$.  The notion of cocycle perturbation 
is useful in the more general context of semigroups acting on 
arbitrary factors (see \cite{arvMono}), but here we confine 
the discussion to the case where $M\sim\mathcal B(H)$.  
\begin{definition}[Cocycle Perturbations]
\label{basicsCocConjDef}
An $E_0$-semigroup $\beta$ of the form (\ref{basicsCoc2.2})
called a {\it cocycle perturbation} of $\alpha$.  Two 
\esg s are said to be {\it cocycle conjugate} if one of 
them is conjugate to a cocycle perturbation of the other. 
\end{definition}
\index{cocycle perturbation}%
\index{cocycle conjugacy}%
\index{e@\esg!cocycle perturbation of}%

The fundamental problem in the theory of \esg s is 
to find an effective classification up to cocycle 
conjugacy for \esg s acting on $\mathcal B(H)$.  
We will have more to say about cocycle 
perturbations, cocycle conjugacy, and 
the classification problem in Section \ref{S:prodSys}.  

\subsubsection{Numerical Index}
The numerical index is the simplest example 
of a cocycle conjugacy invariant for 
\esg s acting on $\mathcal B(H)$, and it is defined 
as follows.  By a {\it unit} for an \esg\ 
$\alpha=\{\alpha_t: t\geq0\}$ acting on 
$\mathcal B(H)$ we mean a strongly continuous 
semigroup $T=\{t_t:t\geq0\}$ of bounded operators 
on $H$ such that $T_0=\mathbf 1$ and which 
satisfies the following commutation relation 
\begin{equation}\label{unitEq}
\alpha_t(A)T_t=T_tA,\qquad t\geq0,\quad A\in\mathcal B(H).   
\end{equation}
The set $\mathcal U_\alpha$ of all units is empty if 
$\alpha$ is of type $III$, but it is nonempty otherwise 
(see Lecture 3).  
In the latter cases, a simple argument based on 
the commutation formula leads to the conclusion 
that for every pair $S,T\in\mathcal U_\alpha$ 
there is a unique complex number $c(S,T)$ with the 
property 
$$
T_t^*S_t=e^{tc(S,T)}\mathbf 1,\qquad t\geq 0.  
$$
The function 
$c: \mathcal U_\alpha\times\mathcal U_\alpha\to\mathbb C$ 
is called the {covariance function} of $\alpha$.  The 
covariance function is conditionally positive definite 
in the sense that for every finite set $T_1,\dots,T_n$ 
of units and every set $\lambda_1,\dots, \lambda_n$ of 
complex numbers satisfying $\lambda_1+\dots+\lambda_n=0$, 
one has 
$$
\sum_{j,k=1}^n\lambda_j\bar\lambda_k c(T_j,T_k)\geq 0.  
$$ 
A familiar construction based on these 
inequalities produces a complex 
Hilbert space $H(\mathcal U_\alpha,c)$, and it is possible 
to show that this Hilbert space is separable whenever 
$\mathcal U_\alpha\neq\emptyset$.  The 
index of $\alpha$ is defined as follows
$$
\ind\,\alpha=
\begin{cases}
\dim H(\mathcal U_\alpha,c),\quad &if\  \mathcal U_\alpha\neq
\emptyset\\
2^{\aleph_0},\quad &if\  \mathcal U_\alpha=\emptyset.
\end{cases}
$$
The possible values of the index are 
$\{0,1,2,\dots,\infty=\aleph_0,c=2^{\aleph_0}\}$, 
and the index is denumerable iff $\alpha$ is not of 
type $III$.  It is a nontrivial fact that we have 
unrestricted validity of the logarithmic addition 
formula
$$
\ind\, \alpha\otimes\beta=\ind\ \alpha + \ind\ \beta.  
$$
Notice that, in order to calculate the index of an 
\esg, one has to calculate the entire set $\mathcal U_\alpha$ 
of its units, as well as the covariance function 
$$
c: \mathcal U_\alpha\times\mathcal U_\alpha\to\mathbb C.
$$
It is significant that these calculations can be 
carried out for specific examples, and in particular, 
we will see in Lecture 3 that the index of the CAR/CCR flow 
of rank $r=1,2,\dots,\infty$ is its rank.  From that 
calculation it follows, for example, that 
the CAR/CCR flow of rank 2 is not conjugate to any 
cocycle perturbation of the CAR/CCR flow of rank 1.  

\subsubsection{Eigenvalue Lists of Normal States}
We have already pointed out that one fundamental 
way that quantum probability differs from 
classical probability is that in quantum theory, there is no 
sensible notion of joint distribution.  Another 
fundamental difference is that while in 
probability theory all nonatomic probability measures 
``look the same", that is not so in quantum theory.  

More precisely, if $(X,\mathcal A,P)$ and 
$(Y,\mathcal B,Q)$ are two probability spaces based 
on standard Borel spaces $(X,\mathcal A)$, 
$(Y,\mathcal B)$ and if in both cases the probability 
of finding any singleton $\{p\}$ is zero, then the 
two probability spaces 
are isomorphic in the sense that there is a 
Borel isomorphism $\phi: X\to Y$ that transforms 
one measure into the other: $P(\phi^{-1}(F))=Q(F)$ 
for every $F\in\mathcal B$.  The quantum analogue 
of a nonatomic probability measure is a normal 
state of $\mathcal B(H)$, and  if 
one is given two normal states $\rho$, $\rho^\prime$ 
defined on $\mathcal B(H)$, $\mathcal B(H^\prime)$ 
respectively, then there may or may not exist 
a $*$-isomorphism $\theta: \mathcal B(H)\to\mathcal B(H^\prime)$ 
which carries one to the other in the sense that 
\begin{equation}\label{conjState}
\rho^\prime(\theta(A))=\rho(A),\qquad A\in\mathcal B(H).  
\end{equation}

Indeed, for $\rho$ and $\rho^\prime$ to be so related it 
is necessary for them to have the same eigenvalue list; 
and we now elaborate on this important invariant.  
By an {\it eigenvalue list}
we mean a decreasing sequence of nonnegative
real numbers $\lambda_1\geq \lambda_2\geq\dots$ with finite 
sum.  Every normal state $\omega$ of a 
type $I$ factor is associated with a positive operator of 
trace $1$, whose eigenvalues counting multiplicity can be 
arranged into an eigenvalue list which will be 
denoted $\Lambda(\omega)$.  
If the factor is finite dimensional,
we still consider $\Lambda(\omega)$ to be an infinite list by 
adjoining zeros in the obvious way.  Given two eigenvalue 
lists $\Lambda=\{\lambda_1\geq\lambda_2\geq\dots\}$ and 
$\Lambda^\prime=\{\lambda_1^\prime\geq\lambda_2^\prime\geq\dots\}$, 
we will write 
$$
\|\Lambda - \Lambda^\prime\| = 
\sum_{k=1}^\infty |\lambda_k-\lambda_k^\prime|
$$
for the $\ell^1$-distance from one list to the other.  
A classical result of Hermann Weyl 
implies that if $\rho$ and $\sigma$ are 
normal states of a type $I$ factor $M$, then we have 
$$
\|\Lambda(\rho)-\Lambda(\sigma)\| \leq \|\rho -\sigma\|.   
$$
Now suppose that $\rho$ and $\rho^\prime$ are normal 
states of $\mathcal B(H)$ and $\mathcal B(H^\prime)$ 
respectively.  Since a $*$-isomorphism  
$\theta:\mathcal B(H)\to\mathcal B(H^\prime)$ 
must be implemented by a unitary 
operator from $H$ to $H^\prime$, it follows that 
$\rho$ and $\rho^\prime$ are conjugate as in 
(\ref{conjState}) only when $\Lambda(\rho)=\Lambda(\rho^\prime)$.  
More generally, assuming that both $H$ and $H^\prime$ 
are of dimension $\aleph_0$, there is a sequence of $*$-isomorphisms 
$\theta_n: \mathcal B(H)\to \mathcal B(H^\prime)$ such that 
$$
\lim_{n\to\infty}\|\rho^\prime\circ\theta_n-\rho\|=0
$$
iff $\rho$ and $\rho^\prime$ have the same eigenvalue list.  

Finally, note that the eigenvalue list of a vector state 
of $\mathcal B(H)$ has the form $\{1,0,0,\dots\}$; and 
more generally, 
the eigenvalue list of a state of $\mathcal B(H)$ has only 
a finite number of nonzero terms iff the state is continuous 
in the weak operator topology of $\mathcal B(H)$.

\begin{theorem}[Existence of Cocycle Perturbations]\label{CocPertThm}
Let $N=1,2,\dots,\infty$ and let $\Lambda=\{\lambda_1\geq\lambda_2\geq\dots\}$ 
be an eigenvalue list with only a finite number of nonzero terms, such 
that $\lambda_1+\lambda_2+\dots=1$.  There 
is a cocycle perturbation $\alpha$ of the CAR/CCR flow of index $N$ 
which is pure, and which has an absorbing state with eigenvalue list $\Lambda$.  
\end{theorem}

The proof of Theorem \ref{CocPertThm} is very indirect, and we 
merely outline the four key ideas behind the argument.  
Starting with the finite 
sequence $\lambda_1\geq\dots\geq \lambda_n>0$ of positive terms 
of the given eigenvalue list, one first constructs a semigroup of unital 
completely positive maps 
$P=\{P_t:t\geq0\}$ acting on 
the matrix algebra $M_n(\mathbb C)$ 
which satisfies $P_t(\mathbf 1)=\mathbf 1$ for every $t\geq0$, 
which is {\it pure} in a sense that is appropriate for CP semigroups, 
and which leaves invariant a state of $M_n(\mathbb C)$ with 
precisely the eigenvalue list $\{\lambda_1,\dots,\lambda_n\}$.  
The second step appeals to the Dilation Theory of Lecture 2 in order 
to find a minimal dilation of $P$ to an $E_0$-semigroup $\alpha$ that 
acts on a Hilbert space $H$.  Making use of minimality, it is possible 
to show that $\alpha$ is a pure \esg\ which acts on a type $I$ 
factor and has a normal invariant 
state with an eigenvalue list whose nonzero terms are precisely 
$\lambda_1,\dots,\lambda_n$ - and by a further 
adjustment one can arrange that the index of $\alpha$ 
is the given integer $N$.  Third, one appeals to results of 
\cite{arvIndDil} which imply that $\alpha$ is completely spatial, and 
finally, one may appeal to the 
classification results of \cite{arvI} for 
completely spatial $E_0$-semigroups to infer that $\alpha$ 
is conjugate to a cocycle perturbation of the CAR/CCR flow 
of index $N$.  

The finiteness hypothesis on the eigenvalue list allowed us 
to work with matrix algebras in constructing 
the initial CP semigroup $P$.  We conjecture that similar 
constructions can be carried out with CP semigroups 
having bounded generators that act on an infinite 
dimensional type $I$ factor, but 
this has not been proved.  As a 
test problem for such techniques, we propose 

\vskip0.1in
{\bf Problem:}
Can the hypothesis that $\Lambda$ is finitely nonzero 
be dropped from Theorem \ref{CocPertThm}?  
\vskip0.1in

Theorem \ref{CocPertThm} asserts that for every type $I$ \esg\ 
$\alpha$ and every finitely nonzero eigenvalue list $\Lambda$, 
there is a cocycle perturbation of $\alpha$ that is pure and which 
has an absorbing state with a specified eigenvalue 
list.  It seems 
reasonable to ask if the same is true if one drops the 
hypothesis that $\alpha$ should be of type $I$.  As a 
limited step in that direction, we conjecture that the 
following question has an affirmative answer.  

\vskip0.1in
{\bf Problem:}  
Can an arbitrary $E_0$-semigroup be perturbed by a cocycle into 
a pure $E_0$-semigroup?

\subsection{Existence of Interactions}

Given a pair of pure \esg s $\alpha^-$, $\alpha^+$ acting, 
respectively on $\mathcal B(H^-)$, $\mathcal B(H^+)$, each of which is 
a cocycle perturbation of a CAR/CCR flow, it is natural 
to ask when they can be assembled into a history.  Thus, we seek 
a simple test for determining when there is a one-parameter 
unitary group $U=\{U_t:t\geq0\}$ acting on $H^-\otimes H^+$ 
with the property that its associated automorphsim group 
$\gamma_t=U_t\cdot U_t^*$ should satisfy 
$$
\gamma_{-t}(A\otimes\mathbf 1_{H^+})=\alpha^-_t(A)\otimes\mathbf 1_{H_-}, 
\qquad 
\gamma_t(\mathbf 1_{H_-}\otimes B)=\mathbf 1_{H_-}\otimes\alpha^+_t(B), 
$$
for all $t\geq 0$, $A\in\mathcal B(H^-)$, $B\in\mathcal B(H^+)$.  

More generally, let 
$\mathcal M$ be a type $I$ subfactor of $\mathcal B(H)$, and 
let $\alpha$, $\beta$ be two $E_0$-semigroups acting, 
respectively, on $\mathcal M$ and its commutant $\mathcal M^\prime$.  
In Lecture 3 we will discuss a general result 
- Theorem \ref{ctpSigPairingProp} - which  provides a 
necessary and sufficient condition for the existence of a 
one-parameter unitary group $\{U_t: t\in\mathbb R\}$ 
acting on $H$ whose associated automorphism group
$\gamma_t(A)=U_tAU_t^*$ has $\alpha$ as its past
and $\beta$ as its future in the sense 
that 
\begin{equation}\label{intIndEq1.1}
\gamma_{-t}\restriction_{\mathcal M}=\alpha_t,\qquad
\gamma_t\restriction_{\mathcal M^\prime}=\beta_t,
\qquad t\geq 0,
\end{equation}
the condition being that the product systems 
of $\alpha$ and $\beta$ are {\it anti-isomorphic}.  

What we require for the current discussion is the 
following corollary 
of Theorem \ref{ctpSigPairingProp}, which can be viewed as 
a counterpart for noncommutative 
dynamics of von Neumann's theorem on the existence of self-adjoint
extensions of symmetric operators in terms of 
their deficiency indices.   

\begin{theorem}\label{intIndThm1}
Let $\alpha$, $\beta$ be two $E_0$-semigroups acting, respectively,  
on $\mathcal B(H)$, $\mathcal B(K)$, 
each of which is a cocycle perturbation of 
a CCR/CAR flow.  The following are equivalent.  
\begin{enumerate}
\item[(i)]
There is a one-parameter group 
$\gamma=\{\gamma_t:t\in\mathbb R\}$ of automorphisms
of $\mathcal B(H\otimes K)$ satisfying 
$$
\gamma_{-t}(A\otimes \mathbf 1)=\alpha_t(A)\otimes\mathbf 1,
\qquad 
\gamma_t(\mathbf 1\otimes B)=\mathbf 1\otimes \beta_t(B),
$$
for all $t\geq 0$, $A\in\mathcal B(H)$, $B\in\mathcal B(K)$.
\item[(ii)] 
$\ind\, \alpha=\ind \,\beta$.  
\end{enumerate}
\end{theorem}

In Section 3 we will show how 
Theorem \ref{intIndThm1} follows from the results 
of that section concerning the relation 
between the structure of product 
systems and the dynamics of histories.  We emphasize 
that when the group $\gamma$ of Theorem \ref{intIndThm1} exists, 
it is not uniquely determined by 
$\alpha$ and $\beta$, and that issue will also be discussed 
in Lecture 3.  Together, 
Theorems \ref{CocPertThm} and \ref{intIndThm1} lead 
one to the following conclusion:

\begin{theorem}[Existence of Interactions]\label{intSemOldThmA}
Let $n=1,2,\dots,\infty$ and let $\Lambda_-$ and 
$\Lambda_+$ be two eigenvalue lists having only 
a finite number of nonzero terms.  There is an 
interaction $(U,\mathcal M)$ whose past and future 
$E_0$-semigroups 
are cocycle perturbations of the $CAR/CCR$ flow 
of index $n$, and whose absorbing 
states have eigenvalue 
lists $\Lambda_-$, $\Lambda_+$ respectively.  
\end{theorem}

\begin{proof}[Proof of Theorem \ref{intSemOldThmA}]
Fix $n=1,2,\dots,\infty$ and let 
$\Lambda_+$, $\Lambda_-$ be eigenvalue 
lists with finitely many nonzero terms.  
By Theorem \ref{CocPertThm}, there is a pure 
cocycle perturbation $\alpha^+$ of 
the CAR/CCR flow of index $n$ which has 
an absorbing state with eigenvalue list $\Lambda_+$.  
Similarly, we find a pure cocycle perturbation 
$\alpha_-$ of the same CAR/CCR flow which has an 
absorbing state with list $\Lambda_-$.   
Since both $\alpha^-$ and 
$\alpha^+$ are pure, we conclude from  Theorem 
\ref{intIndThm1} that there is a history 
with past and future \esg s given, respectively, 
by $\alpha^-$ and $\alpha^+$.  This history must 
be an interaction since both $\alpha^-$ and $\alpha^+$ 
have normal invariant states.  
\end{proof}

\begin{remark}[On the Existence and 
Nonexistence of Dynamics]\label{intIndRem1}
Given an arbitrary $E_0$-semigroup $\alpha$, 
it is natural to ask if cocycle perturbations  
of $\alpha$ 
can represent both the past and future of some history, 
as can cocycle perturbations of the CAR/CCR flows.  
More precisely, is there a history whose past and future 
$E_0$-semigroups are {\it both} conjugate to cocycle perturbations 
of $\alpha$?  

Significantly, the answer can be no.  
Boris Tsirelson has given examples 
of product systems 
which are {\it not} anti-isomorphic 
to themselves \cite{tsirelFirst}.  
By Theorem \ref{specExThm1} below, 
such a product system is 
isomorphic to the product system of some $E_0$-semigroup $\alpha$, 
and it follows that
there are $E_0$-semigroups 
$\alpha$ whose product  systems are not anti-isomorphic
themselves.  Since cocycle  
perturbations of $E_0$-semigroups must have isomorphic 
product systems, the general criteria 
of Theorem \ref{ctpSigPairingProp} 
implies that 
such $E_0$-semigroups cannot serve as both the 
past and future of any history.   

On the other hand, every $E_0$-semigroup $\alpha$ acting
on $\mathcal B(H)$ can serve as the ``past" of some 
automorphism group.  
To see why, let 
$E$ be the product system of $\alpha$ and let $E^{\rm{op}}$ be 
the product system opposite to $E$. 
Theorem \ref{specExThm1} implies that there is an 
$E_0$-semigroup $\beta$, acting on $\mathcal B(K)$, 
whose product system is isomorphic to $E^{\rm{op}}$ and therefore 
anti-isomorphic to $E$.  
We may conclude from the Theorem \ref{ctpSigPairingProp} that 
there is a one parameter group of automorphisms 
$\gamma$ acting on $\mathcal B(H\otimes K)$ which 
satisfies (\ref{intIndEq1.1}) by having $\alpha\otimes\mathbf 1$ 
(acting on $\mathcal B(H)\otimes\mathbf 1$) 
as its past and $\mathbf 1\otimes\beta$ 
(acting on $\mathbf 1\otimes\mathcal B(K)$) as its future.  
We point out that 
a more elementary proof of this extension result, based 
on the theory of spectral subspaces, can 
be found in \cite{arvKish}.  
\end{remark}

\subsection{The Interaction Inequality}

Given a pair of eigenvalue lists $\Lambda_-$, $\Lambda_+$, 
each of which has only a finite number of nonzero terms, and 
given a positive integer $n=1,2,\dots,\infty$, 
Theorem \ref{intSemOldThmA} implies that there is an interaction 
whose past and future absorbing states 
$\omega_-$, $\omega_+$, have eigenvalue 
lists $\Lambda_-$ and $\Lambda_+$ respectively.  One might 
expect that such an interaction should be nontrivial when 
$\Lambda_-\neq \Lambda_+$, but this is is far from self-evident.  

More precisely, we seek a method for computing, or at least 
estimating, the quantity $\|\bar\omega_- - \bar\omega_+\|$ 
in terms of the eigenvalue lists of  
$\omega_-$ and $\omega_+$.  We now describe a general 
solution of this problem involving an inequality that 
appears to be of some interest in its own right.

\begin{theorem}[Interaction Inequality]\label{IntIneq}
Let $(U,M)$ be an interaction with past and future states
$\omega_-$, $\omega_+$ on $M$, $M^\prime$ respectively, 
and let $\bar\omega_-$ and 
$\bar\omega_+$ denote their extensions
to $\gamma$-invariant states of $\mathcal A$.  Then 
$$
\|\bar\omega_- - \bar\omega_+\|\geq
\|\Lambda(\omega_-\otimes\omega_-) - 
\Lambda(\omega_+\otimes\omega_+)\|.  
$$
\end{theorem}

\begin{remark}\label{squareRootRemark}
Notice the tensor product of states on the right.  For 
example, $\Lambda(\omega_-\otimes\omega_-)$ is obtained 
from the eigenvalue list 
$\Lambda(\omega_-)=\{\lambda_1\geq\lambda_2\geq\dots\}$
of $\omega_-$ by rearranging the 
doubly infinite sequence 
of all products $\lambda_i\lambda_j$, $i,j=1,2,\dots$
into decreasing order.  It can be an 
unpleasant combinatorial chore to calculate 
$\Lambda(\omega_-\otimes\omega_-)$ even when 
$\Lambda(\omega_-)$ is relatively simple and 
finitely nonzero; but it is also possible 
to show that if $A$ and $B$ are two positive 
trace class operators such that 
$\Lambda(A\otimes A)=\Lambda(B\otimes B)$, then 
$\Lambda(A) = \Lambda(B)$.  

Note too that we do not assume that the eigenvalue 
lists of $\omega_{\pm}$ are finitely nonzero 
in Theorem \ref{IntIneq}.  
The reader is referred to 
\cite{arvMono} for these details, as well as for 
the proof of Theorem \ref{IntIneq}.  
\end{remark}

The results of the preceding three sections are summarized 
as follows: 

\begin{cor}[Nontriviality of Interactions]\label{bigIntCor}
Let $n=1,2,\dots,\infty$ be a positive integer and let 
$\Lambda_-$, $\Lambda_+$ be two eigenvalue lists with 
finitely many nonzero terms.  Then there is an interaction 
whose past and future \esg s are cocycle perturbations 
of the CAR/CCR flow of index\/ $n$, and whose past and future 
absorbing states $\omega_-$, $\omega_+$ 
have eigenvalue lists $\Lambda_-$, $\Lambda_+$ 
respectively.  

Such an interaction in nontrivial whenever $\Lambda_-\neq \Lambda_+$, 
and in general we have 
$$
\|\bar\omega_- - \bar\omega_+\|\geq \|\Lambda_-\otimes \Lambda_- 
- \Lambda_+\otimes \Lambda_+\|.  
$$
\end{cor}

\begin{remark}[Existence of Strong Interactions]
To illustrate Corollary \ref{bigIntCor}, we use it to show 
that strong interactions exist in the sense that for 
every  $\epsilon>0$ there is an interaction with the property 
\begin{equation}\label{applicationIneq}
\|\bar\omega_- - \bar\omega_+\| > 2-\epsilon.  
\end{equation}
To see that, 
choose positive integers $p<q$ and consider the eigenvalue 
lists 
\begin{align*}
\Lambda_-&=\{1/p,1/p,\dots,1/p,0,0,\dots\}\\
\Lambda_+&=\{1/q,1/q,\dots,1/q,0,0,\dots\},  
\end{align*}
where $1/p$ is repeated $p$ times and $1/q$ is 
repeated $q$ times.  
Corollary \ref{bigIntCor} implies 
that there is an interaction $(U,M)$ whose
past and future $E_0$-semigroups are cocycle perturbations 
of the $CAR/CCR$ flow of index $1$, whose past and 
future absorbing states satisfy 
$\Lambda(\omega_-)=\Lambda_-$ and 
$\Lambda(\omega_+)=\Lambda_+$, and we have 
$$
\|\bar\omega_--\bar\omega_+\|\geq 
\|\Lambda_-\otimes\Lambda_-  -\Lambda_+\otimes\Lambda_+\|.  
$$
If we neglect zeros, 
the eigenvalue list of $\Lambda_-\otimes \Lambda_-$ 
consists of the single eigenvalue $1/p^2$, repeated 
$p^2$ times, and that of $\Lambda_+\otimes\Lambda_+$ consists
of $1/q^2$ repeated $q^2$ times.  Thus 
$$
\|\Lambda_-\otimes\Lambda_- -
\Lambda_+\otimes\Lambda_+\|=
p^2(1/p^2-1/q^2)+(q^2-p^2)/q^2=2-2p^2/q^2,
$$
and the inequality (\ref{applicationIneq}) follows whenever  
$q$ and $p$ satisfy 
$q>p\sqrt{2/\epsilon}$.  
\end{remark}

\specialsection{Generators of Dynamics and Dilation Theory}

In this lecture we describe a new approach to the dilation 
theory of quantum dynamical semigroups 
that is based on the notion of an 
$A$-dynamical system.  These objects provide the 
$C^*$-algebraic structure that underlies much of 
noncommutative dynamics, whether it takes place 
in \cstar\ or a von Neumann algebra, independently 
of issues relating to dilation theory.  
After describing the general properties of $A$-dynamical systems, 
we introduce $\alpha$-expectations and 
noncommutative moment polynomials, and show how these objects 
enter into the construction of $C^*$-dilations.  

Once one is in possession of this $C^*$-algebraic infrastructure, 
one can establish the existence and uniqueness of dilations for 
quantum dynamical semigroups acting on von Neumann algebras in 
a natural way.

\subsection{Generators of Noncommutative Dynamics}

The flow of time in quantum theory is represented by
a one-parameter group of 
$*$-automorphisms $\alpha = \{\alpha_t: t\in\mathbb R\}$
of a \cstar\ $B$.  There is often a $C^*$-subalgebra $A\subseteq B$ 
that can be singled out from physical considerations 
which, together with its time 
translates, generates $B$.  For example, in a nonrelativistic
quantum mechanical system with $n$ degrees 
of freedom the flow of time is represented by a 
one-parameter group of automorphisms of ${\mathcal B}(L^2({\mathbb R}^n))$   
of the form $\alpha_t(T)=e^{itH}Te^{-itH}$, $t\in\mathbb R$, 
where $H$ is a self-adjoint Schr\"odinger 
operator  of the form
$$
H=-\sum_{k=1}^n \frac{\partial^2}{\partial x_k^2}
+V(X_1,\dots,X_n).   
$$
Here, $X_1,\dots,X_n$ denote the configuration observables 
at time $0$ 
$$
X_k:\xi(x_1,\dots,x_n)\mapsto 
x_k\xi(x_1,\dots,x_n),
$$ 
defined appropriately on a common dense domain in 
$L^2({\mathbb R}^n)$, and $V$ denotes the potential 
associated with the interaction forces.  
The functional calculus provides a faithful representation 
of the commutative \cstar\ $C_0({\mathbb R}^n)$ on $L^2({\mathbb R}^n)$ 
by way of 
$f\mapsto f(X_1,\dots,X_n)$, and these functions 
of the configuration operators form 
a commutative $C^*$-subalgebra 
$A\subseteq{\mathcal B}(H)$.  
It is not hard to see that the 
set of all time translates 
$\{\alpha_t(A):t\geq 0\}$ of $A$ 
generates an irreducible $C^*$-subalgebra 
$B$ of ${\mathcal B}(H)$.  In particular, for 
different times $t_1\neq t_2$, the \cstar s $\alpha_{t_1}(A)$ 
and $\alpha_{t_2}(A)$ fail to commute with each other. Indeed, 
{\it no nontrivial relations appear to exist between 
$\alpha_{t_1}(A)$ and $\alpha_{t_2}(A)$
when $t_1\neq t_2$}.  

We now look closely at this phenomenon in general.  
Throughout this lecture $A$ be denote an 
arbitrary but fixed \cstar. 
\begin{definition}\label{dynSystemDef}
An $A$-dynamical system is a triple $(\iota,B,\alpha)$ 
consisting of a semigroup $\alpha=\{\alpha_t: t\geq 0\}$ 
of $*$-endomorphisms acting 
on a \cstar\ $B$ and an injective $*$-homomorphism 
$\iota: A\to B$, such that $B$ is generated 
by $\cup_{t\geq0}\alpha_t(\iota(A))$.  
\end{definition}

Notice that we impose no continuity requirement 
on the semigroup $\alpha_t$ in its time parameter $t$.  
We lighten notation by identifying $A$ with its image 
$\iota(A)$ in $B$, thereby replacing $\iota$ with the inclusion 
map $\iota: A\subseteq B$.  Thus, an $A$-dynamical system 
is a dynamical system $(B,\alpha)$ that contains $A$ as a 
$C^*$-subalgebra in a specified way, with the property that 
$B$ is the norm-closed linear span of finite 
products 
\begin{equation}\label{spanEq}
B=\overline{\rm span}\{\alpha_{t_1}(a_1)\alpha_{t_2}(a_2)\cdots\alpha_{t_k}(a_k)\}
\end{equation}
where $t_1,\dots,t_k\geq 0$, $a_1,\dots,a_k\in A$, $k=1,2,\dots$

Our aim is to say something sensible about the class 
of {\it all} $A$-dynamical systems, and to obtain more 
detailed information about certain of its members.  
The examples described in the preceding 
paragraphs illustrate the fact that in 
even the simplest cases  
where $A$ is $C(X)$, 
the structure of individual 
$A$-dynamical systems can be very complex.   

There is a natural hierarchy in the class of all 
$A$-dynamical systems, defined by 
$(\iota,B,\alpha)\geq (\tilde\iota,\tilde B,\tilde\alpha)$ iff there 
is a $*$-homomorphism $\theta: B\to \tilde B$ satisfying 
$\theta\circ\alpha_t=\tilde\alpha_t\circ\theta$, $t\geq0$, and 
$\theta(a)=a$ for $a\in A$.  Since $\theta$ fixes $A$, 
it follows from (\ref{spanEq}) that  $\theta$ must be 
surjective, $\theta(B)=\tilde B$, hence 
$(\tilde\iota,\tilde B,\tilde\alpha)$ is a {\it quotient} 
of $(\iota,B,\alpha)$.  
Two $A$-dynamical systems are said to be {\it equivalent} if there is 
a map $\theta$ as above that is an isomorphism of \cstar s.  
This will be the case iff each of the $A$-dynamical systems 
dominates the other.  
One may also think of the class 
of all $A$-dynamical systems as a category, whose objects 
are $A$-dynamical systems and whose maps 
$\theta$ are the ones just described.  

\begin{remark}[The Universal $A$-dynamical System]
There is a largest equivalence class in this hierarchy, 
whose representatives are called 
{\it universal} $A$-dynamical systems.  We exhibit one
as follows.  
Consider the free product 
of an infinite family of copies of $A$ indexed 
by the nonnegative real numbers 
$$
{\mathcal P} A=*_{t\geq 0}A_t, \qquad A_t=A.  
$$
Thus, we have a family of $*$-homomorhisms $\theta_t$ 
of $A$ into the \cstar\ ${\mathcal P}A$ such ${\mathcal P}A$ is generated 
by $\cup\{\theta_t(A):t\geq0\}$ and such that the following 
universal property is satisfied: for every  family 
$\bar\pi=\{\pi_t:t\geq 0\}$ of $*$-homomorphisms of $A$ into 
some other \cstar\ $B$, there is a unique  
$*$-homomorphism $\rho: {\mathcal P}A\to B$ such that 
$\pi_t=\rho\circ\theta_t$, $t\geq 0$.  
Nondegenerate representations 
of ${\mathcal P}A$ correspond to families $\bar\pi=\{\pi_t:t\geq 0\}$ 
of representations $\pi_t: A\to {\mathcal B}(H)$ of $A$ on 
a common Hilbert space $H$, subject to no condition other 
than the triviality of their common nullspace 
$$
\xi\in H,\quad \pi_t(A)\xi=\{0\}, \quad t\geq 0 \implies \xi=0.  
$$
A simple argument establishes the existence of ${\mathcal P}A$ by taking 
the direct sum of a sufficiently large set of 
such representation sequences $\bar\pi$.

This definition does not exhibit ${\mathcal P}A$ in concrete 
terms (see \S \ref{S:construction} for that), 
but it does allow us to define a  
universal $A$-dynamical system.  
The universal property of ${\mathcal P}A$ implies that there is 
a semigroup of shift endomorphisms 
$\sigma=\{\sigma_t: t\geq 0\}$ acting on $\mathcal P A$ and 
defined uniquely by $\sigma_t\circ\theta_s=\theta_{t+s}$, 
$s,t\geq 0$.  Using the universal properties of 
${\mathcal P}A$, it is quite easy to verify 
that $\theta_0$ is an injective $*$-homomorphism 
of $A$ in ${\mathcal P}A$,  and we use this map to 
identify $A$ with $\theta_0(A)\subseteq{\mathcal P}A$.  
Thus {\it the triple $(i,{\mathcal P}A,\sigma)$ 
becomes an $A$-dynamical system with the property that every 
other $A$-dynamical system is subordinate to it}.  
\end{remark}

\subsection{$\alpha$-Expectations and $C^*$-Dilations}

Suppose now that we are 
given a semigroup $P=\{P_t: t\geq 0\}$ of completely 
positive contractions acting on a $C^*$-algebra $A$.  
We are interested in singling out certain 
$A$-dynamical systems $(\iota,B,\alpha)$ 
that admit a conditional 
expectation $E: B\to A$ with the property 
$$
E(\alpha_t(a))=P_t(a),\qquad a\in A,\quad t\geq 0.  
$$
Of course, for many $A$-dynamical systems $(\iota,B,\alpha)$ 
there will be no such conditional 
expectation; and even if such an expectation exists, 
there is no reason to expect it to be uniquely 
determined by the preceding formula.  We now 
introduce a class of conditional expectations, 
called $\alpha$-expectations, that occupy a 
central position in the dilation theory of 
CP semigoups.

Let us first 
review some common terminology.  
Let $A\subseteq B$ be an inclusion of \cstar s. 
For any subset $S$ of $B$ we write $[S]$ for 
the norm-closed linear span of $S$.  
The subalgebra 
$A$ is said to be {\it hereditary} 
if for $a\in A$ and $b\in B$, one has 
$$
0\leq b\leq a \implies b\in A.  
$$
The hereditary subalgebra of $B$ generated by 
a subalgebra $A$ is the 
closed linear span $[ABA]$ of all products $axb$, 
$a,b\in A$, $x\in B$, 
and in general 
$A\subseteq [ABA]$.  A {\it corner} of $B$ is a hereditary 
subalgebra of the particular 
form $A=pBp$ where $p$ is a projection 
in the multiplier algebra $M(B)$ of $B$.  

We also make essential 
use of {\it conditional expectations} 
$E: B\to A$.  A conditional expectation 
is an idempotent positive linear map 
with  range $A$, satisfying $E(ax)=aE(x)$
for $a\in A$, $x\in B$.  
Conditional expectations 
are completely positive linear maps of norm 
$1$ whenever $A\neq \{0\}$.  
When $A=pBp$ is a {\it corner} 
of $B$, the map 
$E(x)=pxp$, defines a conditional 
expectation of $B$ onto $A$.  On the other hand, many of
the  conditional expectations encountered 
here do not have this simple 
form, even when $A$ has a unit.  Indeed, 
if $A$ is subalgebra of $B$ that is 
{\it not} hereditary, then there is no natural conditional 
expectation $E: B\to A$.   For example, 
the universal dynamical system $(\iota,{\mathcal P}A,\sigma)$ 
never contains $A$ as a hereditary subalgebra, hence there 
is no ``obvious" conditional expectation $E: {\mathcal P}A\to A$.

\begin{definition}\label{alphaDef} Let 
$(\iota,B,\alpha)$ be an $A$-dynamical system.  An 
$\alpha$-expectation is a conditional 
expectation $E: B\to A$ having the following 
two properties:
\begin{enumerate}
\item[E1.]
Equivariance:  $E\circ\alpha_t=E\circ\alpha_t\circ E$, $t\geq 0$.  
\item[E2.]
The restriction of $E$ to the hereditary subalgebra 
generated by $A$ is multiplicative, 
$E(xy)=E(x)E(y)$, $x,y\in [ABA]$.  
\end{enumerate}
\end{definition}

Note that an {\it arbitrary} 
conditional expectation $E: B\to A$ 
gives rise to a family of linear maps 
$P=\{P_t:t\geq 0\}$ of $A$ to itself by way 
of  $P_t(a)=E(\alpha_t(a))$, 
$a\in A$.  Each $P_t$ 
is a completely positive contraction.   
When $E$ is an $\alpha$-expectation 
property E1 implies that $P_t$ is 
related to $\alpha_t$ by 
\begin{equation}\label{equivariance}
E\circ\alpha_t=P_t\circ E,    \qquad t\geq 0,
\end{equation}
and from the relation 
(\ref{equivariance}) one finds 
that $P$ must satisfy 
the semigroup property $P_s\circ P_t=P_{s+t}$. 

Property E2 is of course automatic 
if $A$ is a hereditary subalgebra of $B$.  
It is a fundamentally {\it noncommutative} 
hypothesis on $B$.  For example, if $Y$ is 
a compact Hausdorff space and $B=C(Y)$, then every 
unital subalgebra $A\subseteq C(Y)$ generates 
$C(Y)$ as a hereditary algebra, and the only linear maps 
$E: C(Y)\to A$ satisfying E2 
are $*$-endomorphisms of $C(Y)$.  

\begin{definition}\label{cstarDilDef}
Let $P=\{P_t: t\geq 0\}$ be a semigroup of completely positive 
contractions acting on a \cstar\ $A$.  By a $C^*$-dilation of 
$(A,P)$ we mean an $A$-dynamical system $(\iota,B,\alpha)$ with 
the additional property that there is an 
$\alpha$-expectation $E: B\to A$ satisfying 
$$
P_t(a)=E(\alpha_t(a)),\qquad a\in A,\quad t\geq 0.  
$$
\end{definition}

We will see in the following section that 
an $\alpha$-expectation $E: B\to A$ is uniquely determined 
by the family of completely positive 
maps $a\in A\mapsto E(\alpha_t(a))$, $t\geq 0$; 
thus, the $\alpha$-expectation associated with a 
$C^*$-dilation of a given CP semigroup $P$ is 
{\it uniquely determined} by $P$.  

The following result implies that $C^*$-dilations 
always exist; and in fact, the universal $A$-dynamical 
system serves as a simultaneous 
$C^*$-dilation of every semigroup 
of completely positive contractions acting on $A$.  

\begin{theorem}\label{universalDilationThm}
For every semigroup of completely positive 
contractions $P=\{P_t: t\geq 0\}$ acting on $A$, there 
is a unique $\sigma$-expectation $E: {\mathcal P}A\to A$ 
satisfying
\begin{equation}\label{compressionEq}
P_t(a)=E(\sigma_t(a)), \qquad a\in A, \quad t\geq 0.  
\end{equation}
\end{theorem}

Both assertions
are nontrivial;  we discuss 
uniqueness in the following section,  
existence is discussed in \S \ref{S:construction}.  

\subsection{Moment Polynomials}\label{S:moments}
The theory of $C^*$-dilations rests on properties 
of certain noncommutative polynomials that are  
defined recursively as follows.

\begin{prop}\label{momPolyProp}
Let $A$ be an algebra over a field $\mathbb F$.  For every 
family of linear maps  $\{P_t:t\geq 0\}$ of $A$ 
to itself satisfying the semigroup property 
$P_{s+t}=P_s\circ P_t$ and $P_0={\rm id}$, 
there is a unique family of multilinear mappings 
from $A$ to itself, indexed by the 
$k$-tuples of nonnegative real numbers, 
$k= 1,2,\dots$, where for a 
fixed $k$-tuple $\bar t=(t_1,\dots,t_k)$
$$
a_1,\dots,a_k\in A\mapsto [\bar t; a_1, \dots,a_k]\in A
$$
is a $k$-linear mapping, all 
of which satisfy 
\begin{enumerate}
\item[MP1.]\label{mp1}
$P_s([\bar t; a_1,\dots,a_k])=[t_1+s,t_2+s,\dots,t_k+s;
a_1,\dots,a_k]$.  
\item[MP2.]\label{mp2}
Given a $k$-tuple for which $t_\ell=0$ for some 
$\ell$ between $1$ and $k$, 
$$
[\bar t; a_1,\dots,a_k]=[t_1,\dots,t_{\ell-1};
a_1,\dots,a_{\ell-1}]a_\ell 
[t_{\ell+1},\dots,t_k; a_{\ell+1},\dots,a_k].
$$
\end{enumerate}
\end{prop}

\begin{remark}
The proofs of both existence and uniqueness are straightforward 
arguments using induction on the number $k$ of variables.  
Note that in the second axiom MP2, 
we make the natural conventions when 
$\ell$  has one of the extreme values $1$, $k$.  For example, if  
$\ell=1$, then MP2 should be interpreted as 
$$
[0,t_2,\dots,t_k; a_1,\dots,a_k]=a_1[t_2,\dots,t_k; a_2,\dots,a_k].  
$$
In particular, in the linear case $k=1$, MP2
makes the assertion
$$
[0;a]=a, \qquad a\in A;
$$
and after applying axiom MP1 one obtains 
$$
[t;a]=P_t(a), \qquad a\in A, \quad t\geq 0.  
$$
One may calculate any particular
moment polynomial explicitly, but 
the computations quickly become a tedious exercise 
in the arrangement of parentheses.  For example, 
\begin{align*}
[2,6,3,4;a,b,c,d]&=P_2(aP_1(P_3(b)cP_1(d))), \\
 [6,4,2,3;a,b,c,d]&=P_2(P_2(P_2(a)b)cP_1(d)).  
\end{align*}

Finally, we remark that when $A$ is a $C^*$-algebra 
and the 
linear maps satisfy $P_t(a)^*=P_t(a^*)$, 
$a\in A$, $t\geq 0$, then the associated moment polynomials 
obey the following symmetry
\begin{equation}\label{symmetryEq}
[t_1,\dots,t_k;a_1,\dots,a_k]^*=[t_k,\dots,t_1;a_k^*,\dots,a_1^*].  
\end{equation}
Indeed, one finds that the 
sequence of polynomials $[[\cdot;\cdot]]$ 
defined by 
$$
[[t_1,\dots,t_k;a_1,\dots,a_k]]=[t_k,\dots,t_1;a_k^*,\dots,a_1^*]^*
$$
also satisfies axioms MP1 and MP2, and hence must coincide 
with the moment polynomials of $\{P_t\}$ by the uniqueness 
assertion of Proposition \ref{momPolyProp}.  
\end{remark}

These polynomials are important because they are the 
expectation values of certain $A$-dynamical systems.  

\begin{theorem}\label{uniquenessThm}
Let $P=\{P_t:t\geq 0\}$ be a semigroup of completely 
positive maps on $A$ satisfying $\|P_t\|\leq 1$, $t\geq 0$, 
with associated 
moment polynomials $[t_1,\dots,t_k;a_1,\dots,a_k]$.  

Let $(i,B,\alpha)$ be an $A$-dynamical system 
and let $E: B\to A$ be an $\alpha$-expectation with the property 
$E(\alpha_t(a))=P_t(a)$, $a\in A$, $t\geq 0$.  Then 
\begin{equation}
E(\alpha_{t_1}(a_1)\alpha_{t_2}(a_2)\cdots\alpha_{t_k}(a_k)) 
=[t_1,\dots,t_k; a_1,\dots,a_k].  \label{nPointEq}
\end{equation}
for every $k=1,2,\dots$, $t_k\geq 0$, $a_k\in A$.  
In particular, there is at most one 
$\alpha$-expectation $E: B\to A$ 
satisfying $E(\alpha_t(a))=P_t(a)$, $a\in A$, $t\geq 0$.  
\end{theorem}

\begin{proof} One applies 
the uniqueness of moment polynomials as follows.  
Properties 
E1 and E2 of Definition \ref{alphaDef} imply that the 
sequence of polynomials $[[\cdot;\cdot]]$ defined by 
$$
[[t_1,\dots,t_k;a_1,\dots,a_k]]=E(\alpha_{t_1}(a_1)\cdots\alpha_{t_k}(a_k))
$$
must satisfy the two axioms MP1 and MP2.  For example, 
notice that E2 implies 
\begin{equation}\label{multipEq}
E(xay)=E(x)aE(y)\qquad x,y\in B, \quad a\in A,
\end{equation} 
since for an approximate 
unit $e_n$ for $A$ we can write $E(xay)$ as follows:
\begin{align*}
\lim_{n\to\infty}e_nE(xae_ny)e_n&=
\lim_{n\to \infty}E(e_nxae_nye_n)
=\lim_{n\to\infty}E(e_nxa)E(e_nye_n)\\
&=\lim_{n\to\infty}e_nE(x)ae_nE(y)e_n=E(x)aE(y).  
\end{align*}
Property (\ref{multipEq}) implies that the 
polynomials $[[\cdot;\cdot]]$ 
satisfy MP2.  Moreover, from 
(\ref{equivariance}) we find that 
\begin{align*}
P_s([[t_1,\dots,t_n;a_1,\dots,a_n]])&=
P_s(E(\alpha_{t_1}(a_1))\cdots\alpha_{t_n}(a_n)))=\\
E(\alpha_s(\alpha_{t_1}(a_1)\cdots\alpha_{t_n}(a_n)))&=
E(\alpha_{t_1+s}(a_1)\cdots\alpha_{t_n+s}(a_n)),
\end{align*}
hence MP1 is satisfied as well.  
Thus 
(\ref{nPointEq}) follows from the uniqueness assertion 
of Proposition \ref{momPolyProp}.  
The uniqueness of the $\alpha$-expectation 
associated with $P$ 
is now apparent from formulas (\ref{nPointEq}) and 
(\ref{spanEq}).  
\end{proof}

\subsection{Construction of $C^*$-dilations}\label{S:construction}

The proof of the existence assertion of 
Theorem \ref{universalDilationThm} 
is based on a construction 
that exhibits ${\mathcal P}A$ as 
the enveloping \cstar\ of a Banach $*$-algebra $\ell^1(\Sigma)$, 
in such a way that the desired $\alpha$-expectation 
appears as a completely positive map on 
$\ell^1(\Sigma)$.  We now describe this construction 
of $\mathcal P A$ in some detail, but we do not
prove that the $\alpha$-expectation 
is completely positive here.  Full details can be found in 
\cite{arvMono}.  

Let $S$ be 
the set of finite sequences $\bar t=(t_1,t_2,\dots,t_k)$
of nonnegative real numbers $t_i$, $k=1,2,\dots$ 
which have distinct neighbors, 
$$
t_1\neq t_2, t_2\neq t_3,\dots, t_{k-1}\neq t_k.  
$$
Multiplication and  
involution are defined in $S$ as follows.  
The product of two elements 
$\bar s=(s_1,\dots,s_k), \bar t = (t_1,\dots,t_\ell)\in S$
is defined by conditional concatenation 
$$
\bar s\cdot \bar t=
\begin{cases}
(s_1,\dots,s_k, t_1,\dots, t_\ell), \qquad\text{if } s_k\neq t_1,\\
(s_1,\dots,s_k, t_2,\dots, t_\ell), \qquad\text{if } s_k= t_1,
\end{cases}
$$
where we make the natural conventions when $\bar t = (t)$ 
is of length 1, namely $\bar s\cdot (t)=(s_1,\dots, s_k,t)$ 
if $s_k\neq t$, and $\bar s\cdot (t)=\bar s$ 
if $s_k=t$.  
The involution in $S$ is defined by reversing the 
order of components  
$$
(s_1,\dots, s_k)^*=(s_k,\dots,s_1).  
$$
One  finds that $S$ is an associative $*$-semigroup.  

Fixing a \cstar\ $A$, we 
attach a Banach space $\Sigma_\nu$ to  
every $k$-tuple $\tau=(t_1,\dots,t_k)\in S$ as follows
$$
\Sigma_\tau = 
\underbrace{A\hat\otimes\cdots\hat\otimes A}_{k\text{\ times}},
$$
the $k$-fold projective tensor product of copies of the 
Banach space $A$.  We assemble the $\Sigma_\tau$ 
into a family of Banach spaces 
over $S$, $p: \Sigma\to S$, by way of 
$\Sigma = \{(\tau, \xi): \tau\in S, \xi\in E_\tau\}$, $p(\tau,\xi)=\tau$.  

We introduce a multiplication 
in $\Sigma$ as follows.  Fix 
$\lambda=(\lambda_1,\dots, \lambda_k)$ and $\mu=(\mu_1,\dots,\mu_\ell)$ in $S$ 
and choose 
$\xi\in \Sigma_\mu$, $\eta\in \Sigma_\nu$.  If $\lambda_k\neq \mu_1$ 
then $\xi\cdot\eta$ is defined as the tensor product 
$\xi\otimes\eta \in \Sigma_{\mu\cdot\nu}$.  If 
$\lambda_k=\mu_1$ then we must tensor over $A$ and make the 
obvious identifications.  More explicitly, in this case 
there is a natural map of 
the tensor product $\Sigma_\mu\otimes_A \Sigma_\nu$ 
onto 
$\Sigma_{\mu\cdot\nu}$ by making identifications 
of elementary tensors as follows:
$$
(a_1\otimes\cdots\otimes a_k)\otimes_A(b_1\otimes\cdots\otimes b_\ell)
\sim
a_1\otimes\cdots\otimes a_{k-1}\otimes 
a_kb_1\otimes b_2\otimes\cdots\otimes b_\ell.  
$$
With this convention  
$\xi\cdot\eta$ is defined by
$$
\xi\cdot\eta=\xi\otimes_A\eta\in \Sigma_{\mu\cdot\nu}.  
$$ 
This defines an associative multiplication in the family 
of Banach spaces $\Sigma$.  There is also a natural 
involution in $\Sigma$, defined on each 
$\Sigma_\mu$, $\mu=(s_1,\dots,s_k)$ as the unique 
antilinear isometry to $\Sigma_{\mu^*}$ satisfying
$$
((s_1,\dots,s_k), a_1\otimes\cdots\otimes a_k)^* =
((s_k,\dots,s_1), a_k^*\otimes\cdots\otimes a_1^*).  
$$
This defines an isometric antilinear mapping of the 
Banach space $\Sigma_\mu$ onto $\Sigma_{\mu^*}$, for 
each $\mu\in S$, and thus the structure $\Sigma$ becomes 
an involutive $*$-semigroup in which each fiber 
$\Sigma_\mu$ is a Banach space.  

Let $\ell^1(\Sigma)$ be the Banach $*$-algebra 
of summable sections.  
The norm and involution are the natural ones
$\|f\|=\sum_{\mu\in \Sigma}\|f(\mu)\|$, 
$f^*(\mu)=f(\mu^*)^*$.  Noting that 
$\Sigma_\lambda\cdot\Sigma_\mu\subseteq\Sigma_{\lambda\cdot\mu}$, 
the multiplication in
$\ell^1(\Sigma)$ is defined by convolution
$$
f*g(\nu)=\sum_{\lambda\cdot\mu=\nu}f(\lambda)\cdot g(\mu), 
$$
and one easily verifies that $\ell^1(\Sigma)$ is a Banach 
$*$-algebra. 

For $\mu=(s_1,\dots,s_k)\in S$ and $a_1,\dots,a_k\in A$ we 
define the function 
$$
\delta_\mu\cdot a_1\otimes\dots\otimes a_k\in \ell^1(\Sigma)
$$
to be zero except at $\mu$, and at $\mu$ it 
has the value $a_1\otimes\dots\otimes a_k\in \Sigma_\mu$.  
These elementary functions have $\ell^1(\Sigma)$ as their 
closed linear span.  Finally, there is a natural family of 
$*$-homomorphisms 
$\theta_t: A\to \ell^1(\Sigma)$, $t\geq 0$, defined by 
$$
\theta_t(a)=\delta_{(t)}\cdot a, \qquad a\in A,\quad t\geq 0, 
$$
and these maps are related to the generating sections by
$$
\delta_{(t_1,\dots,t_k)}\cdot a_1\otimes\dots\otimes a_k=
\theta_{t_1}(a_1)\theta_{t_2}(a_2)\cdots\theta_{t_k}(a_k).  
$$

The algebra $\ell^1(\Sigma)$ fails to have a unit, but it has the 
same representation theory 
as ${\mathcal P}A$ in the following sense.  
Given a family of representations 
$\pi_t:A\to{\mathcal B}(H)$, $t\geq 0$, 
fix $\nu=(t_1,\dots,t_k)\in S$.
There is a unique bounded linear operator   
$L_\nu: \Sigma_\nu\to{\mathcal B}(H)$ of norm 
$1$ that is defined by its action on elementary 
tensors as follows 
$$
L_\nu( a_1\otimes\cdots\otimes a_k)=
\pi_{t_1}(a_1)\cdots\pi_{t_k}(a_k).  
$$
Thus there is a bounded 
linear map $\tilde\pi: \ell^1(\Sigma)\to{\mathcal B}(H)$ 
defined by 
$$
\tilde\pi(f)=\sum_{\mu\in S}L_\mu(f(\mu)), 
\qquad f\in \ell^1(\Sigma).     
$$
One finds that $\tilde\pi$ 
is a $*$-representation of 
$\ell^1(\Sigma)$ with
$\|\tilde\pi\| = 1$.  
This representation satisfies 
$\tilde\pi\circ\theta_t=\pi_t$, 
$t\geq0$.  Conversely, every bounded 
$*$-representation $\tilde\pi$ of $\ell^1(\Sigma)$ 
on a Hilbert space $H$ is associated with 
a family of representations 
$\pi_t$, $t\geq0$, of $A$ on $H$ by way
of $\pi_t=\tilde\pi\circ\theta_t$.

The results of the preceding discussion 
are summarized as follows: 

\begin{prop}
The enveloping \cstar\ $C^*(\ell^1(\Sigma))$, together 
with the family of homomorphisms 
$\tilde\theta_t: A\to C^*(\ell^1(\Sigma))$, $t\geq0$,
defined by promoting the homomorphisms 
$\theta_t: A\to \ell^1(\Sigma)$, has the 
same universal property as the infinite free product 
${\mathcal P}A=*_{t\geq0}a$, and is therefore isomorphic 
to ${\mathcal P}A$.  
\end{prop}

Notice that there is a natural semigroup of $*$-endomorphisms 
of $\ell^1(\Sigma)$ defined by 
$$
\sigma_t:\delta_{(s_1,\dots,s_k)}\cdot\xi\mapsto
\delta_{(s_1+t,\dots,s_k+t)}\cdot\xi, \qquad 
(s_1,\dots,s_k)\in \Sigma,\quad\xi\in \Sigma_\nu
$$ 
and it promotes to the natural 
shift semigroup of ${\mathcal P}A=C^*(\ell^1(\Sigma))$.  
The inclusion of $A$ in $\ell^1(\Sigma)$ is given by 
the map $\theta_0(a)= \delta_{(0)}a\in \ell^1(\Sigma)$, 
and it too promotes 
to the natural inclusion of $A$ in ${\mathcal P}A$.  

Finally, we fix a semigroup of completely positive 
contractions $P_t: A\to A$, $t\geq 0$, and consider the 
associated moment polynomials of Proposition \ref{momPolyProp}.  
Since $\|P_t\|\leq 1$ for every $t\geq0$, a 
straightforward inductive argument based on the 
two properties MP1 and MP2 shows that 
$$
\|[t_1,\dots,t_n; a_1,\dots,a_n]\|\leq \|a_1\|\cdots \|a_n\|,
\qquad t_k\geq0, \quad a_k\in A, 
$$
hence there is a unique 
bounded linear map $E_0: \ell^1(\Sigma)\to A$ 
satisfying 
$$
E_0(\delta_{(t_1,\dots,t_k)}\cdot a_1\otimes\cdots\otimes a_k) 
=[t_1,\dots,t_k; a_1,\dots,a_k], 
$$
for $(t_1,\dots,t_k)\in S$, $a_1,\dots,a_k\in A$, 
$k=1,2,\dots$, and in fact $\|E_0\|\leq 1$.  
Using the axioms MP1 and 
MP2, one finds that the map 
$E_0$ preserves the adjoint (see Equation (\ref{symmetryEq})),  
satisfies the conditional expectation property 
$E_0(af)=aE_0(f)$ for $a\in A$, $f\in\ell^1(\Sigma)$, 
that the restriction of $E_0$ to the 
``hereditary" $*$-subalgebra 
of $\ell^1(\Sigma)$ spanned by 
$\theta_0(A)\ell^1(\Sigma)\theta_0(A)$ is multiplicative, 
and that it is related to $\phi$ 
by $E_0\circ\sigma=\phi\circ E_0$ and 
$E_0(\sigma(a))=\phi(a)$, $a\in A$.  Thus, $E_0$ 
satisfies the axioms of Definition \ref{alphaDef}, suitably 
interpreted for the Banach $*$-algebra $\ell^1(\Sigma)$.  

In view of the basic fact that a bounded completely positive 
linear map of a Banach $*$-algebra to $A$ 
promotes naturally to a completely positive map of 
its enveloping \cstar \ to $A$, the critical 
property of $E_0$ reduces to:

\begin{theorem}\label{cpThm}
For every $n\geq 1$,  
$a_1,\dots,a_n\in A$, and $f_1, \dots,f_n\in\ell^1(\Sigma)$, 
we have 
$$
\sum_{i,j=1}^n a_j^*E_0(f_j^*f_i)a_i\geq 0.  
$$
Consequently, $E_0$ extends uniquely through the 
completion map $\ell^1(\Sigma)\to{\mathcal P}A$ to 
a completely positive map $E_\phi: {\mathcal P}A\to A$ that 
becomes a $\sigma$-expectation 
satisfying Equation (\ref{compressionEq}).  
\end{theorem}

\begin{cor}\label{cpCor}
Every semigroup of completely positive contractions 
acting on a $C^*$-algebra has a $C^*$-dilation.  
\end{cor}

\subsection{Existence of $W^*$-dilations}
We now show how to use the results of 
the preceding section 
to obtain dilations appropriate 
for the category of 
von Neumann algebras.  
The dynamical issues discussed in Lecture 
1 involved pairs of $E_0$-semigroups, 
but now it is appropriate 
to broaden the context.  
While we are primarily interested in the case of 
$E_0$-semigroups acting on type $I_\infty$ factors, 
the basic concepts of dilation theory are best formulated in 
greater generality, and the general formulation has 
certain advantages.  

\begin{definition}[$E$-semigroups]
An $E$-semigroup is 
a semigroup $\{\alpha_t:t\geq 0\}$ of 
normal $*$-endomorphisms of a von Neumann algebra $M$ that 
obeys the natural continuity requirement in its time 
variable; for every normal linear functional 
$\rho$ on $M$, the functions $t\mapsto \rho(\alpha_t(x))$, 
$x\in M$, should be continuous.  
\end{definition}

Let $(M,\alpha)$ be a pair consisting of 
a von Neumann algebra $M$ with 
separable predual and an $E$-semigroup 
$\alpha=\{\alpha_t: t\geq 0\}$ acting 
on it.  
The operators 
$\alpha_t(\mathbf 1)$ form a decreasing 
family of projections in $M$ in general, and 
if one has $\alpha_t(\mathbf 1)=\mathbf 1$ for 
every $t\geq0$, then $\alpha$ is called 
an $E_0$-semigroup, {\it even when 
$M$ is not a type $I$ factor}.  
In order to deal effectively 
with dilation-theoretic issues, we must allow 
for more general von Neumann algebras $M$ and 
the possibility of non-unital semigroups.  
The general issues discussed in 
this section do not depend on spatial aspects of 
$M$, and for the most part  we will not have 
to realize $M$ in any concrete representation 
as a subalgebra of $\mathcal B(H)$.  

A {\it corner} of $M$ is a von Neumann subalgebra 
of the particular form $N=pMp$,  where 
$p$ is a projection in $M$.  The corner is 
said to be {\it full} if the central carrier of 
$p$ is $\mathbf 1$, and in that case $pMp$ is 
a factor iff $M$ is a factor of the same type.  

Given an arbitrary projection $p\in M$, one can 
ask if there is a semigroup of completely 
positive maps $P=\{P_t: t\geq 0\}$ that 
acts on the corner $pMp$ and is related to 
$\alpha$ as follows
\begin{equation}\label{dilCompEq1}
P_t(pxp)=p\alpha_t(x)p,\qquad t\geq 0, \quad x\in M.  
\end{equation}
Such maps $P_t$ 
need not exist in general; for example, 
taking $x=\mathbf 1-p$, one finds that a necessary 
condition for $P_t$ to exist is that 
$p$ should satisfy $p\alpha_t(\mathbf 1-p)p=0$.  
Equivalently, a projection $p\in M$ is said to be 
{\it coinvariant} under $\alpha$ if 
\begin{equation}\label{dilCompEq2}
\alpha_t(\mathbf 1-p) \leq \mathbf 1-p,\qquad t\geq 0.  
\end{equation}
\index{coinvariant projection}%

\begin{remark}[Increasing Projections and \esg s]\label{dilCompRem1}
A projection $p\in M$ is said to be {\it increasing}
if it has the property 
\begin{equation}\label{dilCompEq3}
\alpha_t(p)\geq p,\qquad t\geq 0.  
\end{equation}
\index{increasing projection}%
Notice that in general, an increasing projection must 
be coinvariant.  Indeed, since $\alpha_t(\mathbf 1)\leq \mathbf 1$, 
we will have 
$$
\alpha_t(\mathbf 1-p)=
\alpha_t(\mathbf 1)-\alpha_t(p)\leq \mathbf 1-p
$$
whenever $p$ is an increasing projection.  The converse
is not necessarily true.  But in 
the special case where $\alpha$ is an 
\esg, $\alpha_t(\mathbf 1-p)=\mathbf 1-\alpha_t(p)$; 
we conclude that  {\it a projection is coinvariant under an \esg\ 
iff it is an increasing projection}.  
\end{remark}

Now for any projection $p\in M$, 
one can define a family of linear maps 
$P = \{P_t: t\geq 0\}$ on $N=pMp$ 
by compressing each map $\alpha_t$ as follows
\begin{equation}\label{dilCompEq4}
P_t(a) = p\alpha_t(a)p,\qquad a\in pMp,\quad t\geq 0.  
\end{equation}
Obviously, each $P_t$ is a normal completely positive 
linear map of $pMp$ into itself satisfying 
$\|P_t\|\leq 1$ for every $t\geq 0$.  More significantly, 
one easily establishes: 

\begin{prop}\label{dilCompProp1}
Let $p$ be a coinvariant projection for $\alpha$  
and consider the family of maps $P=\{P_t: t\geq 0\}$ 
of $pMp$ defined by (\ref{dilCompEq4}).  
$P$ is a continuous semigroup of completely 
positive contractions, satisfying (\ref{dilCompEq1}).  
If, in addition, $\alpha$ is an \esg, then 
we have $P_t(p)=p$, $t\geq 0$.  
\end{prop}

Dilation theory in the category 
of von Neumann algebras concerns 
the properties of completely positive semigroups that 
can be obtained from  $E$-semigroups in this particular way.  
By a CP semigroup 
we mean a pair $(N,P)$ where  $P=\{P_t:t\geq0\}$ 
is a semigroup of normal completely positive 
linear maps 
acting on a von Neumann algebra $N$ which satisfies 
$\|P_t\|\leq 1$ for every $t\geq 0$.

\begin{definition}[Dilation and Compression]\label{minDef1.3}
A triple $(M,\alpha,p)$ consisting of 
an $E$-semigroup $\alpha=\{\alpha_t: t\geq 0\}$ acting 
on a von Neumann algebra $M$, together with a 
distinguished coinvariant 
projection $p\in M$, is called a dilation triple.  
Let $N=pMp$ be the corner of $M$ associated with $p$ 
and let $P=\{P_t: t\geq 0\}$ be the semigroup acting 
on $N$ as follows
\begin{equation}\label{dilCompEq5}
P_t(a)=p\alpha_t(a)p,\qquad t\geq 0, \quad a\in N.  
\end{equation}
The CP semigroup $(N,P)$ called a compression of $(M,\alpha,p)$, 
and $(M,\alpha,p)$ is called a dilation of $(N,P)$.  
\end{definition}

We emphasize that the notion of a compression 
to a subalgebra has meaning only when 
(a) the subalgebra 
is a corner $pMp$ of $M$ and (b) the projection $p$ 
satisfies (\ref{dilCompEq2}).  
There are several equivalent notions 
of {\it minimality} that are associated 
with this dilation theory \cite{arvMono}; but in 
this lecture we confine attention 
to the question of existence.  However, 
we point out that an arbitrary dilation 
can always be reduced to a minimal 
one, and that a {\it necessary} condition for 
$(M,\alpha,p)$ to be a minimal dilation 
of $(N,P)$ is that the central carrier of 
$p$ should be $\mathbf 1$.  Thus, in the context 
of minimal dilations, if $N$ is a factor then  
$M$ must be a factor of the same type.  

Starting 
with a CP semigoup $(N,P)$, in order to find a dilation 
$(M,\alpha,p)$ of $(N,P)$ one has 
to find a way of embedding $N$ as a corner $pMp$ of 
a larger von Neumann algebra $M$, on which there is 
a specified action of an $E$-semigroup $\alpha$ 
that is related to $P$ as 
above.  
Notice that Corollary
\ref{cpCor}  provides the following 
infrastructure.  If we view $N$ as a unital 
\cstar\ and $P$ as a semigroup of contractive completely 
positive maps on $N$, forgetting the continuity of $P_t$ in its  
time variable, 
then we can assert that the pair $(N,P)$ has a 
$C^*$-dilation $(\iota,B,\alpha)$.  Certainly, 
$B$ is not a von Neumann algebra and 
$\alpha$ is not an $E$-semigroup; 
thus $(\iota,B,\alpha)$ cannot 
serve as a $W^*$-dilation 
of $(N,P)$.  However, it is possible to make use of 
the $\alpha$-expectation $E: B\to N$ 
to find another dilation of $(N,P)$ that is 
subordinate to $(\iota,B,\alpha)$ and has all the desired properties.  
The results are summarized as follows.  
\begin{theorem}
[Existence of $W^*$-dilations]\label{dilEx2Thm1}
Let $\{P_t:t\geq0\}$ be a contractive $CP$-semigroup acting 
on a von Neumann algebra $N$ with separable predual.  Then $(N,P)$ has 
a dilation $(M,\alpha,p)$.
\end{theorem}

\begin{proof}[Idea of Proof]
Considering $P=\{P_t:t\geq0\}$ as a semigroup of completely 
positive contractions acting on the unital \cstar\ $N$, 
we see from 
Corollary \ref{cpCor} that $P$ has a 
$C^*$-dilation $(\iota,B,\alpha)$.  
We may obviously assume that $N\subseteq\mathcal B(H)$ acts 
concretely and nondegenerately on some 
separable Hilbert space $H$.  We will construct 
a representation $\pi$ of $B$ on a Hilbert space $K\supseteq H$ 
with the property that each $\alpha_t$ can be extended 
to a normal $*$-endomorphism of the weak closure $M$ of $\pi(B)$, 
and this will provide the required dilation of $(N,P)$.  
The representation $\pi$ 
is obtained as follows.  

Let $E: B\to N$ be the $\alpha$-expectation associated with 
$(\iota,B,\alpha)$.  Since we may view $E$ as a completely positive 
map of $B$ to $\mathcal B(H)$, it has a minimal Stinespring 
decomposition $E(x)=V^*\pi(x)V$, $x\in B$, where $\pi$ is a 
representation of $B$ on another Hilbert space $K$ and $V: H\to K$ 
is a bounded linear map such that $VH$ has $K$ as its closed linear 
span.  

Let $M$ be the von Neumann algebra $\pi(B)^{\prime\prime}$.  The 
remainder of the proof amounts to showing first, that there is a 
unique $E$-semigroup $\tilde\alpha$ acting on $M$ that satisfies 
$\tilde\alpha_t(\pi(x))=\pi(\alpha_t(x))$, $t\geq 0$, $x\in B$, second, 
that $p=VV^*$ is a projection in $M$ whose corner 
$pMp$ can be naturally identified with $N$, and that 
after this identification is made, $(M,\tilde\alpha,p)$ becomes 
a dilation triple for $(N,P)$.  
\end{proof}

\subsection*{Historical Remarks}
Several approaches to dilation theory for semigroups 
of completely positive maps have been proposed over the 
years, including work of Evans and Lewis \cite{EvLewis}, 
Accardi et al \cite{AFL}, K\"ummerer \cite{kumMarkovDil}, 
Sauvageot \cite{sauv86}, 
and many others.  Our attention was drawn to these 
developments by work of Bhat \cite{bhatMin}, building 
on work of Bhat and Parthasarathy 
\cite{BhatPar} for noncommutative Markov processes, 
in which the first dilation 
theory for CP semigroups acting on ${\mathcal B}(H)$ 
emerged that was effective for our work on $E_0$-semigroups
\cite{arvPureAbsorbing}, \cite{arvInteractions}.  
SeLegue \cite{selThesis} 
showed how to apply multi-operator dilation theory 
to obtain Bhat's dilation result for CP semigroups 
acting on $\mathcal B(H)$, and 
he calculated the expectation values of the 
$n$-point functions of such dilations.  Recently, 
Bhat and Skeide \cite{BhSkeide} have initiated an 
approach to the subject 
that is based on Hilbert modules over \cstar s 
and von Neumann algebras.  
\vfill

\specialsection{The Role of Product Systems.}\label{S:prodSys}

The fundamental problem in the theory of $E_0$-semigroups 
is to find a complete set of computable 
invariants for cocycle conjugacy.  
There is some optimism that such a classification is 
possible, but we are far from achieving that goal.  
We have not yet seen all possible cocycle 
conjugacy classes of $E_0$-semigroups, and we do not have 
an effective and computable set of invariants for 
the ones we have seen.   

The numerical index is an example of a cocycle conjugacy 
invariant.  However, while it 
serves to classify type $I$ $E_0$-semigroups 
up to cocycle conjugacy, it is far from being a complete 
invariant for examples of type $II$ and it is degenerate 
for examples of type $III$.  The gauge 
group introduced below 
provides a more subtle cocycle conjugacy invariant, 
but it has not been calculated 
except in type $I$ cases.  

In this lecture we describe how the classification problem 
can be reduced to the problem of classifying certain simpler 
objects (product systems) up to natural isomorphism, and 
we discuss the role of product systems in other dynamical 
issues related to the theory of $E_0$-semigroups.  
This reformulaton has already lead to significant 
progress in several directions.  We discuss these aspects 
in this lecture, and complete the discussion in Lecture 4.  

\subsection{Product Systems and Cocycle Conjugacy.}\label{SS:prodSysCC}
It is appropriate to discuss concrete 
product systems within the context of $E$-semigroups 
acting on $\mathcal B(H)$.  
Given an $E$-semigoup $\alpha=\{\alpha_t:t\geq0\}$ and $t>0$, consider
the linear  space of operators 
\begin{equation}\label{basicsPsOpSp}
\mathcal E_\alpha(t)=
\{T\in\mathcal B(H):\alpha_t(A)T=TA,\quad A\in\mathcal B(H)\}.  
\end{equation}
We assemble the various $\mathcal E_\alpha(t)$ 
into a family of vector spaces 
$p: \mathcal E_\alpha\to (0,\infty)$ over the interval 
$(0,\infty)$ in the natural way 
\begin{equation}\label{basicsPsFam}
\mathcal E_\alpha=\{(t,T): t>0,\quad T\in\mathcal E_\alpha(t)\} 
\subseteq (0,\infty)\times\mathcal B(H), 
\end{equation}
where $p$ is the projection $p(t,T)=t$.  

The family $\mathcal E_\alpha$ has three 
important properties.  First, the operator 
norm on each particular space 
$\mathcal E_\alpha(t)$ is actually a Hilbert 
space norm.  
To see how the inner product is defined, 
choose two elements 
$S,T\in\mathcal E_\alpha(t)$ and an 
arbitrary operator $A\in\mathcal B(H)$.   
Writing 
$$
T^*SA=T^*\alpha_t(A)S=(\alpha_t(A^*)T)^*S=
(TA^*)^*S=AT^*S, 
$$
we find that $T^*S$ must be a scalar multiple 
of the identity, and the value of that scalar 
defines an inner product on $\mathcal E_\alpha(t)$: 
\begin{equation}\label{basicsPsInPr}
T^*S=\langle S, T\rangle \cdot\mathbf 1.  
\end{equation}
The operator norm on 
$\mathcal E_\alpha(t)$ coincides with 
the norm defined by this inner product, 
since for $T\in\mathcal E_\alpha(t)$ we 
have 
$$
\|T\|^2 = \|T^*T\|=\|\langle T,T\rangle\cdot\mathbf 1\|
=\langle T,T\rangle.  
$$
In particular, $\langle T,T\rangle =0$ iff $T=0$, 
and we conclude that 
each fiber $\mathcal E_\alpha(t)$ becomes a Hilbert space
relative  to the inner product 
defined by (\ref{basicsPsInPr}).

Second, one may verify directly that 
$\mathcal E_\alpha(s)\mathcal E_\alpha(t)\subseteq\mathcal E_\alpha(s+t)$, 
so that the family 
$\mathcal E_\alpha$ can be made 
into an associative semigroup by 
defining multiplication as follows:
\begin{equation}\label{basicsPsMult}
(s,S)\cdot (t,T) = (s+t,ST),  
\end{equation}
and this multiplication makes $p$ into 
a homomorphism of the multiplicative 
structure of $\mathcal E_\alpha$ onto 
the additive semigroup of positive reals.  

Third, this multiplication {\it acts 
like tensoring} in the sense that for every  
$s,t>0$ there 
is a unique unitary operator 
$W_{s,t}: \mathcal E_\alpha(s)\otimes \mathcal E_\alpha(t)
\to \mathcal E_\alpha(s+t)$ satisfying 
$W_{s,t}(S\otimes T) = ST$ 
for all $S\in\mathcal E_\alpha(s)$, 
$T\in\mathcal E_\alpha(t)$.

The weak 
operator topology on $\mathcal B(H)$ generates 
a $\sigma$-algebra of subsets of $\mathcal B(H)$, 
whose elements we refer to as Borel sets. 
This makes $B(H)$ 
into a standard Borel 
space because $H$ is separable.  
We now describe an appropriate context 
for the structure $\mathcal E_\alpha$.

\begin{definition}[Concrete Product System]
\label{basicsPsCoProdSys}
A concrete product system is a Borel subset 
$\mathcal E$ of the cartesian product of 
Borel spaces $(0,\infty)\times\mathcal B(H)$ 
which has the following properties.  Let 
$p:\mathcal E\to(0,\infty)$ be the natural 
projection $p(t,T)=t$.  We require that 
$p$ should be surjective, and in addition: 
\begin{enumerate}
\item[(i)]
For every $t>0$, the set of operators 
$\mathcal E(t)=p^{-1}(t)$ is a 
norm-closed linear subspace of 
$\mathcal B(H)$ with the property that 
$B^*A$ is a scalar for every 
$A,B\in\mathcal E(t)$.  
\item[(ii)]
For every $s,t>0$, $\mathcal E(s+t)$ is 
the norm-closed linear span of the 
set of products $\mathcal E(s)\mathcal E(t)$.  
\item[(iii)]
As a measurable 
family of Hilbert spaces, 
$\mathcal E$ is isomorphic to the trivial 
family $(0,\infty)\times H_0$, where 
$H_0$ is a separable infinite dimensional 
Hilbert space.  
\end{enumerate}
\end{definition}

Property (iii) requires some elaboration.  
The appropriate notion of isomorphism here is 
the one that belongs with the abstract theory 
of the following section, namely:

\begin{definition}\label{basicsPsConcreteIsoDef}
Two concrete product systems 
$\mathcal E\subseteq(0,\infty)\times\mathcal B(H)$
and $\mathcal F\subseteq (0,\infty)\times\mathcal B(K)$ 
are said to be {\it isomorphic} if there is an isomorphism 
of Borel spaces $\theta:\mathcal E\to \mathcal F$, 
which restricts to a unitary operator from $\mathcal E(t)$ 
to $\mathcal F(t)$ for every $t>0$, and which 
satisfies $\theta (xy)=\theta(x)\theta(y)$, 
$x,y\in\mathcal E$. 
\end{definition}
There are several topologies 
that can used to topologize a given concrete product 
systems $\mathcal E$.  
However, a basic result in the theory of Borel  
structures implies that all of these topologies  
generate the same Borel structure on $\mathcal E$ 
because it is a standard Borel space.  Thus 
any ``topological" isomorphism of product systems must 
be an isomorphism in the sense of 
Definition \ref{basicsPsConcreteIsoDef}.  
This feature of product systems allows for considerable 
flexibility.  

Item (iii) makes the assertion that {\it there 
is an isomorphism of Borel spaces 
$\theta: \mathcal E\to (0,\infty)\times H_0$ 
with the property that the restriction 
of $\theta$ to each fiber $\mathcal E(t)$ 
is a unitary operator with range $H_0$}.  
The fact that the concrete product system 
$\mathcal E_\alpha$ associated with an 
$E$-semigroup $\alpha$ satisfies (iii) is 
nontrivial, and we refer the reader to 
\cite{arvMono} for the proof, as well as 
for other characterizations of this property.

Finally, we point out that every concrete product 
system $\mathcal E\subseteq (0,\infty)\times\mathcal B(H)$ 
arises in the above way from a {\it unique} 
$E$-semigroup $\alpha$ acting on $\mathcal B(H)$, 
in such a way that 
the correspondence $\mathcal E \leftrightarrow \alpha$ is a bijection.  
More precisely: 

\begin{prop}\label{ctpFocRecoverProp}
Let $\mathcal E\subseteq (0,\infty)\times\mathcal B(H)$
be a concrete product system.  
There is a unique $E$-semigroup 
$\alpha=\{\alpha_t:t\geq0\}$ acting on $\mathcal B(H)$ whose 
endomorphisms satisfy the following two conditions 
for every $t>0$
\begin{enumerate}
\item[(i)]
$\alpha_t(A)T=TA$, for every $T\in\mathcal E(t)$, 
$A\in\mathcal B(H)$.  
\item[(ii)]
$\alpha_t(\mathbf 1)$ is the projection on $[\mathcal E(t)H]$.  
\end{enumerate}
Moreover, $\mathcal E=\mathcal E_\alpha$ is 
the concrete product system associated with $\alpha$.  
\end{prop}

In order to define $\alpha_t$ for $t>0$, one chooses 
an orthonormal basis $T_1(t), T_2(t),\dots$ for $\mathcal E(t)$ 
and sets $\alpha_t(A)=T_1(t)AT_1(t)^*+T_2(t)AT_2(t)^*+\cdots$.  
In particular, we conclude that 
{\it $E_0$-semigroups correspond bijectively 
with concrete product systems $p: \mathcal E\to(0,\infty)$ 
with the property that $[\mathcal E(t)H]=H$ for every 
$t>0$}.  

It is significant that the fibers $\mathcal E(t)$ of 
a concrete product system are all infinite dimensional 
except in degenerate cases.  More precisely, let 
$V=\{V_t:t\geq 0\}$ be a strongly continuous semigroup 
of isometries acting on a Hilbert space $H$ and let 
$\alpha$ be the associated $E$-semigroup
\begin{equation}\label{trivEsgEq}
\alpha_t(A)=V_tAV_t^*,\qquad A\in\mathcal B(H),\quad t\geq 0.  
\end{equation}
One may verify by a direct calculation that 
$\mathcal E_\alpha(t)=\mathbb C\cdot V_t$ for every 
$t>0$, and in particular the spaces $\mathcal E_\alpha(t)$ 
are all one-dimensional.  It is less obvious that 
the converse is true.  Indeed, for any concrete 
product system $p: \mathcal E\to (0,\infty)$, the 
fiber spaces $\mathcal E(t)=p^{-1}(t)$ are either 
all infinite dimensional or else they are 
all one-dimensional, 
and in the latter case the underlying 
$E$-semigroup must have the form (\ref{trivEsgEq}).  

Let us now examine the connection between cocycle conjugacy 
and product systems.  Let $\alpha$ be an $E_0$-semigroup 
acting on $\mathcal B(H)$, let $U=\{U_t:t\geq 0\}$ be 
an $\alpha$-cocycle, i.e., a strongly continuous family 
of unitary operators in $H$ which satisfy 
\begin{equation}\label{cocEq}
U_{a+t}=U_s\alpha_s(U_t),\qquad s,t\geq 0,   
\end{equation}
and let $\beta_t(A)=U_t\alpha_t(A)U_t^*$ be the corresponding 
cocycle perturbation of $\alpha$.  We can define a 
map $\theta: \mathcal E_\alpha\to \mathcal E_\beta$ as follows: 
$$
\theta(t,T)=(t,U_tT),\qquad T\in\mathcal E_\alpha(t), \quad t>0.  
$$
Indeed, one sees that 
$U_t\mathcal E_\alpha(t)\subseteq \mathcal E_\beta(t)$ for 
every $t>0$ because 
$$
\beta_t(A)U_tT=U_t\alpha_t(A)T=U_tTA,\qquad A\in\mathcal B(H), 
$$
and in fact $U_t\mathcal E_\alpha(t)=\mathcal E_\beta(t)$.  
Thus, $\theta$ defines an isomorphism of measurable families 
of Hilbert spaces that is unitary on fibers.  $\theta$ also 
preserves multiplication by virtue of the cocycle equation 
(\ref{cocEq}) and the fact that 
$\alpha_s(U_t)S=SU_t$: 
\begin{align*}
\theta((s,S)(t,T))&=\theta(s+t,ST)=(s+t,U_{s+t}ST)=
(s+t,U_s\alpha_s(U_t)ST)\\
&=(s+t,U_sSU_tT)=(s,U_sS)(t,U_tT)=\theta(s,S)\theta(t,T).
\end{align*}
It can be seen by a more involved argument that this 
argument is reversible in the sense that every other 
$E_0$-semigroup 
$\beta$ that acts on $\mathcal B(H)$ 
whose product system 
$\mathcal E_\beta$ is isomorphic to $\mathcal E_\alpha$ 
must actually be a 
cocycle perturbation of $\alpha$ (see \cite{arvMono}).
With these results in hand, one easily deduces:

\begin{theorem}\label{basicsPsCocConj}
Two \esg s $\alpha$ and $\beta$ acting, 
respectively, on $\mathcal B(H)$ and 
$\mathcal B(K)$, 
are cocycle conjugate iff their 
product systems are isomoprhic.  
\end{theorem}

Theorem \ref{basicsPsCocConj} implies that 
the cocycle conjugacy class of an \esg\ is completely 
determined by the {\it structure} of its 
product system.  Thus the classification of 
\esg s up to cocycle conjugacy should begin 
with a systematic development of 
a general theory of continuous tensor products of Hilbert 
spaces that is appropriate for \esg s.  
We turn now to this more general discussion.

\subsection{Abstract Product Systems.}\label{SS:absProdSys}

An effective theory should have at least the following 
three properties.  First, it should be possible to
associate a continuous tensor product 
of Hilbert spaces with every \esg.  Second, 
the structure of the
continuous tensor product associated with 
an \esg\  should 
be an invariant for cocycle conjugacy, 
and optimally a complete invariant.  
Third, every continuous tensor product 
in this category of objects should 
be associated with an \esg.  

Given a theory with these three properties, 
the isomorphism classes of such 
structures will then provide a full description 
of the cocycle conjugacy classes of 
$E_0$-semigroups; hence the problem of classifying 
\esg s to cocycle conjugacy is reduced 
to the problem of classifying these 
considerably simpler objects, their 
symmetries, and their related structures.  

We now describe an 
axiomatic approach to continuous tensor products 
of Hilbert spaces that satisfies these requirements.  
We will 
find that the first two requirements follow 
from the results described in the preceding 
section; the proof 
that the third property is satisfied 
requires additional $C^*$-algebraic 
tools, and will 
be taken up in Lecture 4.

This program has already led to significant 
progress in the classification of 
\esg s up to cocycle conjugacy,  in two directions.  
The classification of type $I$ \esg s
was completely settled through an 
analysis of the structure of their 
product systems \cite{arvI}, and the description 
of decomposable product systems in \cite{arvPaths}.  
In \cite{powNewEx}, Powers 
showed that for every $n=1,2,\dots$ there is a 
type $II$ \esg s whose numerical index 
is $n$; hence there are infinitely many 
type $II$ \esg s that are mutually non 
cocycle conjugate.  In a previous paper 
\cite{powTypeIII}, he showed that type $III$ \esg s exist, 
but the methods of \cite{powTypeIII} do not 
lend themselves to differentiating between 
cocycle conjugacy classes in the 
family of examples exhibited there.  
More recently, Tsirelson 
showed \cite{tsirelFirst} that there 
is a continuum of type $II$ product 
systems that are mutually non-isomorphic, and that 
there is a product system that is not anti-isomorphic 
to itself.  In a subsequent paper 
\cite{tsirelSecond}, he constructed a one 
parameter family of type 
$III$ product systems that are mutually non 
isomorphic; we describe these product systems 
in the following section.    
Given that every product system can be 
associated with an \esg\ (see Theorem 
\ref{specExThm1} below), 
Tsirelson's results on product systems 
have immediate implications for the 
problem of classifying 
\esg s up to cocycle conjugacy.


Heuristically, a {\it product system} is a measurable 
family of Hilbert spaces $E = \{E(t): t > 0\ \}$ which 
behaves as if each $E(t)$ were a continuous tensor product
\begin{equation}\label{ctpProdHeuristicEq}
E(t) = \bigotimes_{0<s<t} H_s,\qquad H_s = H
\end{equation}
of copies of a single separable Hilbert space $H$.  While
this heuristic picture is often useful, one must be careful
not to push it too far.  Indeed, we will see that this picture
is basically correct for the simplest examples of product
systems, but that there are other examples with the 
remarkable property that the ``germ" $H$ fails to exist.  
We first illustrate the essentials of the structure of product
systems in the discrete case, where the positive real line 
is replaced with the discrete set $\mathbb N = \{ 1, 2, \dots\}$ of 
positive integers.  Then we will indicate how 
to change the axioms to pass from $\mathbb N$ to $\mathbb R^+$.  

Let $H$ be a separable Hilbert space.  For every 
$n = 1,2 \dots$ let $E(n)$ be the full tensor product of
$n$ copies of $H$:
$$
E(n) = \underbrace{H\otimes H\otimes\dots\otimes H}_{n \text{ times}}.  
$$
We may organize these spaces into a family of Hilbert spaces 
$p: E\to \mathbb N$ over $\mathbb N$ by setting 
$$
E = \{(t,\xi): t\in\mathbb N, \quad \xi\in E(t)\,\},
$$
with projection $p(t,\xi) = t$.  We introduce an associative 
multiplication on the structure $E$ by making use of the tensor product 
$$
(s,\xi)\cdot(t,\eta) = (s+t,\xi \otimes \eta),
$$
$\xi \in E(s)$, $\eta\in E(t)$.  This multiplication is bilinear on fibers,
and has the two additional properties
\begin{align}
E(s+t)&=\overline{\text{span}}\ E(s)\cdot E(t),\qquad\qquad s,t\in\mathbb N 
\label{ctpProdEq4.1}\\  
<ux,vy>&=<u,v><x,y>,\quad 
 u,v\in E(s),\quad x,y\in E(t).
\label{ctpProdEq4.2}
\end{align}
Notice that the Hilbert space associated with the sections of 
$p:E\to \mathbb N$ is the direct sum
$$
\sum_{t\in\mathbb N} E(t) = \sum_{n=1}^\infty H^{\otimes n},
$$
namely the full Fock space over the one-particle space $H$. 
\index{Fock space!over a one-particle space}%

A {\it unit} is a section $n\in\mathbb N\mapsto u_n\in E(n)$ satisfying 
$$
u_{m+n} = u_m u_n,
$$
and which is not the zero section.  The most general unit has the 
form 
$$
u_n = \underbrace{x\otimes x\otimes\dots\otimes x}_{n\ \text{times}}
$$
$n\geq 1$, where $x$ is a nonzero element of the one-particle space $H$.  
A vector $u\in E(n)$ is called 
{\it decomposable} if for every $k=1,2,\dots,n-1$ there are vectors 
$v_k\in E(k)$, $w_k\in E(n-k)$ such that
$$
u = v_k w_k.
$$
The most general decomposable vector in $E(n)$ is an 
elementary tensor of the form
$$
u = x_1\otimes x_2\otimes\dots\otimes x_n,
$$
where $x_k\in H$ for $k=1,2,\dots,n$.

A product system is a similar structure, except that it is 
associated with the space of positive reals rather than $\mathbb N$.
\begin{definition}\label{ctpProdMainDef}
A  product system is a family of separable
Hilbert spaces $p:E \to (0,\infty)$
over the semi-infinite interval $(0,\infty)$,
with fiber Hilbert spaces $E(t)=p^{-1}(t)$, 
endowed with an associative multiplication that 
restrcts to a bilinear map on fibers
\begin{equation*}
(x,y)\in E(s)\times E(t)\mapsto xy\in E(s+t),\qquad s,t>0, 
\end{equation*}
which acts like tensoring in the sense 
that properties (\ref{ctpProdEq4.1}) and 
(\ref{ctpProdEq4.2}) are satisfied.  
In addition, $E$ should be endowed 
with the structure of a standard Borel space 
that is compatible with 
projection onto $(0,\infty)$, multiplication, the 
vector space operations and the inner 
product, 
and which has the further 
property that there should be a separable Hilbert space $H$ such that
\begin{equation}\label{ctpProd4.3}
E \cong (0,\infty)\times H,
\end{equation}
where $\cong$ denotes an isomorphism of 
measurable families of Hilbert spaces.  
\end{definition}
\index{product system!abstract}%
\index{standard Borel space}%
\begin{remark}\label{ctpProdBasisRem}
For example, measurability of the inner product
means that if one considers the 
subset $\Delta=\{(x,y)\in E\times E: p(x)=p(y)\}$ of 
the standard Borel space $E\times E$, then $\Delta$ 
is a Borel subset because $p:E\to(0,\infty)$ is 
a Borel measurable function, and measurability 
of the inner product means that the complex-valued 
function defined on $\Delta$ by 
$(x,y)\mapsto \langle x,y\rangle$ should be 
Borel-measurable.  

The requirement 
(\ref{ctpProd4.3}) is nontrivial, and is the 
counterpart for this category 
of local triviality of Hermitian vector bundles.  It 
is equivalent to the existence of a sequence of measurable 
sections $t\in(0,\infty)\mapsto e_n(t)\in E(t)$ with 
the property that $\{e_1(t), e_2(t),\dots\}$ is an 
orthonormal basis for $E(t)$, for every $t>0$.    
\end{remark}

\begin{definition}\label{ctpProdAbstractDef}
By an isomorphism of product systems we mean an 
isomorphism of Borel spaces $\theta: E\to F$, such that 
$\theta(xy)=\theta(x)\theta(y)$ for all $x,y\in E$, 
whose restriction to each fiber space is a unitary 
operator $\theta_t: E(t)\to F(t)$, $t>0$.   
\end{definition}

One can show easily that 
a concrete product system in the sense 
of Definition \ref{basicsPsCoProdSys} is a 
product system in the the more abstract 
sense of Definition \ref{ctpProdMainDef}.  
We have also seen that every 
\esg\ $\alpha$ gives rise to a 
concrete product 
system $\mathcal E_\alpha$, and  Theorem 
\ref{basicsPsCocConj} asserts that $\alpha$ 
and $\beta$ are cocycle conjugate iff their product 
systems $\mathcal E_\alpha$ and $\mathcal E_\beta$ 
are isomorphic.  Moreover, we will see in 
Lecture 4 that every abstract product system 
is isomorphic to the product system $\mathcal E_\alpha$  
associated with some \esg\ $\alpha$.  Thus, {\it the 
problem of classifying product systems up to 
isomorphism becomes 
a central problem in noncommutative dynamics}.  

One might expect that it should be possible to write down
a comprehensive list of (continuous) product systems as we have 
done above for their discrete analogues.  
In the discrete case there is, up to isomorphism,
exactly one ``product system" for every integer 
$d = 1, 2, \dots,\aleph_0$, and of course $d$ is the dimension 
of the one-particle space $E(1)$.  In the continuous case, however,
nothing like that is true.  While there is a family of ``natural" 
examples parameterized by the values $d = 1, 2, \dots, \aleph_0$, 
there are many others as well.  

\subsection{Examples of Product Systems}
We begin this section 
by describing the simplest examples of product 
systems and we discuss 
their role in the classification of type $I$ \esg s.  
We then describe 
Tsirelson-Vershik product systems, without technical details.  

\subsubsection{Exponential Product Systems} Recall first 
the basic features of the symmetric Fock space $e^H$ over a 
one-particle space $H$.  It is defined as the direct sum of 
Hilbert spaces 
$$
e^H=\sum_{n=0}^\infty H^n
$$
where $H^n$ denotes the $n$-fold 
symmetric tensor product of copies of $H$ when $n\geq 1$, 
and $H^0=\mathbb C$.  There is a natural 
exponential map $f\in H\mapsto \exp(f)\in e^H$ defined by 
$$
\exp(f)=\sum_{n=0}^\infty \frac{1}{\sqrt n}f^{\otimes n}.   
$$
$e^H$ is the closed linear span of the 
set of exponentials $\exp(f)$, $f\in H$, 
and one has 
$$
\langle \exp(f),\exp(g)\rangle = e^{\langle f,g\rangle}, 
\qquad f,g\in L^2((0,\infty);K).  
$$
This construction is functorial in that for 
every unitary operator $U: H_1\to H_2$ there is a 
natural second quantization  
$\Gamma(U): e^{H_1}\to e^{H_2}$ defined by 
$$
\Gamma(U) = V_0\oplus V_1\oplus V_2\oplus\cdots
$$
where, for $n\geq 1$, $V_n: H_1^n\to H_2^n$ is the $n$-fold 
tensor product of copies of $U$, and where $V_0$ is the identity 
map of $\mathbb C$.  Equivalently, $\Gamma(U)$ is defined 
as the unique unitary operator satisfying 
$\Gamma(U):\exp(f)\to\exp(Uf)$, $f\in H$.  
Another fundamental property of this construction is that 
for any two Hilbert spaces $H_1$, $H_2$, $e^{H_1\oplus H_2}$ 
is naturally identified with the tensor product 
$e^{H_1}\otimes e^{H_2}$; 
indeed, there is a unique unitary 
operator $W: e^{H_1\oplus H_2}\to e^{H_1}\otimes e^{H_2}$ 
that satisfies 
$$
W: \exp(f\oplus g)\mapsto \exp(f)\otimes \exp(g), 
\qquad f\in H_1,\quad 
g\in H_2.  
$$

The simplest product systems are the ones 
associated with the CAR/CCR flows, and 
are described as follows.  
Let $N$ be a positive integer 
or $\infty=\aleph_0$ and let $K$ be a Hilbert space 
of dimension $N$.  We form the 
symmetric Fock space 
$e^{L^2((0,\infty);K))}$ over the one-particle space 
$L^2((0,\infty);K)$ consisting of all square-integrable 
$K$-valued functions $\xi: (0,\infty)\to K$.   
For every $t>0$ let
$E(t)$ be the closed subspace 
$$
E(t)=e^{L^2((0,t);K)}\subseteq e^{L^2((0,\infty);K)}, 
$$
$L^2((0,t);K)$ denoting the subspace of $L^2((0,\infty);K)$ 
consisting of functions vanishing almost everywhere outside 
$(0,t)$, 
and consider the family of Hilbert spaces 
$p:E_N\to(0,\infty)$ defined by 
$$
E_N=\{(t,\xi):t>0, \quad \xi\in E(t)\}, \qquad p(t,\xi)=t.
$$  
The shift semigroup 
$S=\{S_t:t\geq 0\}$ acting on $L^2((0,\infty);K)$ 
gives rise to a semigroup of isometries by way of 
second quantization 
$$
U_t=\Gamma(S_t),\qquad t\geq 0,
$$
and we use $\{U_t:t\geq 0\}$ to introduct 
a multiplication in $E_N$ as follows.  Noting 
that $U_s$ maps $e^{L^2((0,t);K)}$ onto 
$e^{L^2((s,s+t);K)}$, and noting the natural 
identification 
$$
e^{L^2((0,s+t);K)}= e^{L^2((0,s);K)}\otimes e^{L^2((s,s+t);K)},
$$
we can define multiplication 
in $E_N$ as follows: for $f\in e^{L^2((0,s);K)}$ and 
$g\in e^{L^2((0,t);K)}$ we set 
$$
(s,f)\cdot (t,g)=(s+t,f\otimes U_sg).  
$$
$E_N$ is a closed subset of the 
Polish space 
$$
E_N\subseteq (0,\infty)\times e^{L^2((0,\infty);K)},
$$ 
the topology being the obvious one arising from 
the usual metric on $(0,\infty)$ and the 
Hilbert space norm of $e^{L^2((0,\infty);K)}$, 
hence $E_N$ is a standard Borel space.  
It is not hard to verify directly that $E_N$ is 
a product system.  

The CCR flow of rank $N$ is an \esg\ $\alpha$ that acts 
on $\mathcal B(e^{L^2((0,\infty);K)})$.  
It is defined most explicitly in terms of 
the the natural representation of the canonical 
commutation relations on $e^{L^2((0,\infty);K)}$, 
and its concrete product system $\mathcal E_\alpha$ 
is isomorphic to $E_N$: the following result 
exhibits that 
isomorphism $E_N\sim\mathcal E_\alpha$.  Since 
an \esg\ is uniquely defined by its product 
system, it will not be necessary to reiterate 
the definition of $\alpha$ in terms of the 
CCRs here.

\begin{prop}\label{ctpProdCCRidentProp}
For every $t>0$, $f\in L^2((0,t);K)$, 
there is a bounded operator
$T_f$ 
on the symmetric Fock space 
$e^{L^2((0,\infty);K)}$, 
defined uniquely on the spanning 
set of vectors 
$\{\exp(g):g\in L^2((0,\infty);K)\}$ 
by
\begin{equation}\label{basicsCompProdOpDef}
T_f(\exp(g))= f\otimes \exp(S_tg), 
\quad t\geq 0, \quad g\in L^2((0,\infty);K).   
\end{equation} 
 The mapping 
$\theta: E_N\to (0,\infty)\times \mathcal B(e^{L^2((0,\infty);K)})$ 
defined by  
\begin{equation}\label{ctpProdThetaEq}
\theta :(t,f) \mapsto (t,T_f)
\end{equation}
is an isomorphism of $E_N$ onto a concrete 
product system $\mathcal E$ acting 
on the Hilbert space $e^{L^2((0,\infty);K)}$.  
$\mathcal E$ is 
the product system of the CCR flow of rank $N$.   
\end{prop}

\begin{definition}\label{ctpProdExpPSDef}
For every $N=1,2,\dots,\infty$, 
$E_N$ is called the exponential product system 
of rank $N$.  
\end{definition}

\begin{remark}[Type of a Product System]
Let $p: E\to(0,\infty)$ be a product system and choose $t>0$.  
A vector $u\in E(t)$ is said to be {\it decomposable} if for 
every $0<s<t$ there are vectors $v\in E(s)$, $w\in E(t-s)$ such 
that $u=vw$.  The set of all decomposable vectors in $E(t)$ 
spans a closed subspace $D(t)\subseteq E(t)$, which can 
be the trivial subspace $D(t)=\{0\}$, but in general one have 
$$
D(s+t)=\overline{\rm {span}}D(s)D(t), \qquad s,t>0.  
$$  
Moreover, if $D(t_0)\neq\{0\}$ for some particular 
$t_0>0$ then $D(t)\neq\{0\}$ for every $t>0$.  The product 
system $E$ is said to be of type $I$ if $D(t)$ spans $E(t)$ 
for some (and therefore every) positive $t$, of type $II$ 
if it is not of type $I$ but $D(t)\neq \{0\}$ for some 
(and therefore every) postive $t$, and of type $III$ if 
$D(t)=\{0\}$ for some (and therefore every) $t>0$.  An 
\esg\ $\alpha$ is said to be of type $I$, $II$, or $III$ 
according as its product system is of that type.  

In the section \ref{SS:dimFn} below we will 
describe a reformulation of the 
index invariant of $E_0$-semigroups 
into a dimension function defined on the 
category of product systems.  
The basic results on the classification of type $I$ product systems 
and type $I$ \esg s are summarized as follows (see \cite{arvMono}):  
Every type $I$ product system $E$ is isomorphic to an exponential 
product system $E_N$, where $N=\dim E$.  Every type $I$ 
\esg\ $\alpha$ is conjugate to a cocycle perturbation of 
the CAR/CCR flow of rank $N$, where  $N=\ind(\alpha)$.

\end{remark}

\subsubsection{Tsirelson-Vershik Product Systems}

Tsirelson and Vershik \cite{verTs} 
constructed a family of 
continuous tensor product of Hilbert 
spaces that could not be described in terms of the 
classical symmetric Fock space construction.  Responding to 
a question of the author, Tsirelson \cite{tsirelSecond} 
adapted those ideas so as to generate a 
one parameter family of product systems; moreover, by  
very ingenious arguments he was able to show that these 
product systems are not only of type $III$, but they 
are mutually non-isomorphic.  In this section we 
describe these product 
systems in enough detail so that their basic 
features are exposed, avoiding most technicalities.  
The reader is referred to Tsirelson's contribution 
to these proceedings for more detail.  
Our exposition differs somewhat from the original; for 
example, we use complex-valued Gaussian random variables 
rather than real-valued ones, and 
our treatment of quasiorthogonality is formulated
somewhat differently.  But the resulting 
constructions are fundamentally the same as Tsirelson's.  

We first describe the correlation functions 
of Tsirelson and Vershik.  
For every real number $\theta>1$ we fix, once and for 
all, a continous real-valued function 
$C_\theta$, defined on the punctured line $\mathbb R\setminus\{0\}$, 
which vanishes outside some small interval 
$(-\epsilon,+\epsilon)\setminus\{0\}$ 
with $\epsilon<1$ ($\epsilon$ may 
depend on $\theta$), and which 
has the following properties 
\begin{enumerate}
\item[(i)]
The restriction of 
$C_\theta$ to the positive real line 
is nonnegative, continous, decreasing, and convex.  
\item[(ii)]
$C_\theta(-t)=C_\theta(t)$, 
for every $t\in\mathbb R\setminus\{0\}$.  
\item[(iii)]
For some positive number $0<\delta<\epsilon$ we have 
$$
C_\theta(t)=\frac{1}{|t|\cdot |\log |t||^\theta}, \qquad 
0<|t|<\delta.
$$
\end{enumerate}

Obviously, the conditions (i)--(iii) can be achieved with $\epsilon$ 
as small as we please.  
Notice the singularity of $C_\theta$ near the 
origin; and this is an essential feature.  
Thus, such functions $C_\theta$ 
``approximate" the delta function.  
The singularity is mild enough that 
these functions belong to $L^1(\mathbb R)$, 
and therefore they define bounded convolution operators 
$f\mapsto C_\theta*f$ on the Hilbert space 
$L^2(\mathbb R)$.

Heuristically, one thinks of $C_\theta$ as the correlation function 
$$
C_\theta(t-s)=E(X_s X_t), \qquad s,t\in\mathbb R,
$$ 
associated with a stationary Gaussian random process 
$\{X_t: t\in \mathbb R\}$.  However, because of the singularity 
(iii) at $t=0$, this formula cannot be achieved with a 
classical process, but only with a 
stationary random distribution \cite{GelVil}.  
This random distribution shares certain 
properties with white noise, and it is suggestive to think 
of it as ``off-white" noise or, as Tsirelson puts it, 
slightly colored noise.

The second key property of these functions 
is that they are positive definite in the sense that 
the sesquilinear form they define on $L^2(\mathbb R)$ by way 
of 
\begin{equation}\label{innProdEq}
\langle f,g\rangle = \int_\mathbb R (C_\theta* f)(x)\bar g(x)\,dx 
\end{equation}
is a semidefinite inner product on $L^2(\mathbb R)$.  Indeed, any 
real function 
$C_\theta\in L^1(\mathbb R)$ that satisfies properties (i) and (ii) above 
can be seen to be positive definite, and therefore the form 
defined by (\ref{innProdEq}) will 
satisfy $\langle f,f\rangle\geq 0$ for all $f\in L^2(\mathbb R)$ (see 
Chapter 14 of \cite{arvMono} for more detail).  
The completion of $L^2(\mathbb R)$ in the inner product 
$\langle\cdot,\cdot\rangle$ is a Hilbert space that will be denoted 
$H_\theta$.  

Third, notice that the inner product of (\ref{innProdEq}) is invariant 
under the action of the one-parameter 
group of translation operators acting 
on $L^2(\mathbb R)$ in the sense that 
\begin{equation}\label{innProdInvEq}
\langle T_tf, T_tg\rangle = \langle f,g\rangle, \qquad f,g\in L^2(\mathbb R).  
\end{equation}
Thus there is a unique one-parameter unitary group 
$U=\{U_t: t\in\mathbb R\}$ in $\mathcal B(H_\theta)$ 
that extends that action 
of the translation operators to $H_\theta$.  

For a bounded interval $I=(a,b)\subseteq \mathbb R$ we write 
$L^2(I)$ for the subspace of $L^2(\mathbb R)$ consisting of 
functions that vanish almost everywhere on the complement 
of $I$.  
Finally, notice that if $I=(a,b)$ and $J=(c,d)$ are 
two bounded intervals whose separation is greater than $2\epsilon$, 
then $L^2(I)$ and $L^2(J)$ are orthogonal in $H_\theta$:
$$
\langle f,g\rangle=0, \qquad f\in L^2(I), \quad g\in L^2(J).  
$$

\subsubsection{Quasiorthogonal Subspaces}
If the separation between $I$ and $J$ 
is smaller than $2\epsilon$ then this is 
no longer true.  But for arbitrary disjoint intervals 
there is a generalized sense in which it is approximately 
true.  Let $H$ be a Hilbert space, let 
$M_1,\dots, M_n$ be a finite set of closed subspaces of $H$, 
and let $M_1\oplus\cdots\oplus M_n$ be the direct 
sum of Hilbert spaces.  There is a unique bounded linear 
map $L: M_1\oplus\cdots\oplus M_n\to H$ satisfying 
$$
L(\xi_1,\dots,\xi_n)=\xi_1+\dots+\xi_n,\qquad \xi\in M, \quad \eta\in N.  
$$
In general, $L$ is a bounded linear map of $M_1\oplus\cdots\oplus M_n$ 
onto the algebraic sum $M_1+\cdots+M_n\subseteq H$.  

\begin{definition}
A finite set of subspaces $M_1, \dots, M_n$ is said to be 
quasiorthogonal if the linear map $L$ is injective, and 
in addition $\mathbf 1-L^*L$ is a Hilbert Schmidt operator 
on $M_1\oplus\cdots\oplus M_n$.  
\end{definition}

\begin{remark}[Equivalence Operators]\label{tsirGauRem2}
The hypothesis that $\mathbf 1-L^*L$ is compact 
implies that the range of $L$ is closed.  It follows that 
the algebraic linear span $M_1+\dots+M_n$ is a closed 
subspace of $H$ that is linearly isomorphic to the 
orthogonal direct sum of the various $M_k$.  
More generally, let $L: H_1\to H_2$ be a bounded operator 
from one Hilbert space to another.   
We will say 
that $L$ is an {\it equivalence operator} if it is 
one-to-one, has dense range, and 
\index{equivalence operator}%
$\mathbf 1-L^*L$ is a Hilbert Schmidt operator on 
$H_1$.  This forces 
the range of $L$ to be closed, hence $L$ is 
an invertible operator.  One verifies easily that 
the both the inverse and the adjoint 
of an equivalence operator are equivalence 
operators, and that 
the class of equivalence operators is closed 
under composition.  In particular, the 
set of equivalence operators in $\mathcal B(H)$ is a 
subgroup of the general linear group of $H$.  
\end{remark}

The first important property of these subspaces 
of $H_\theta$ is based on ideas originating 
in \cite{verTs}; a detailed proof appears in 
\cite{tsirelSecond}.  

\begin{theorem}\label{tsirFirstThm}
Let $I_1,\dots,I_n$ be a pairwise disjoint sequence of 
bounded intervals in $\mathbb R$ and let $H_\theta(I_k)$ be the 
closure of $L^2(I_k)$ in $H_\theta$.  Then 
$H_\theta(I_1),\dots, H_\theta(I_n)$ is quasiorthogonal set
of subspaces of $H_\theta$.  
\end{theorem}

\subsubsection{Equivalence Operators and Gaussian Measures}
The second key element of this construction 
of product systems 
is the following result that describes Gaussian 
measures that are mutually absolutely continuous 
with respect to a given one, and which is in some sense 
a refinement of Kakutani's characterization of 
mutual absolute continuity for infinite product 
measures \cite{kakProd}.  The result is due 
to Feldman \cite{feldmGau}, H\'ajek \cite{haj}, 
and Segal \cite{segCan}.  

Let $(\Omega,\mathcal B,P)$ be a probability space, 
which we may assume is modelled on a standard Borel 
space $(\Omega,\mathcal B)$.  By a Gaussian random variable 
we mean a complex-valued function in $L^2(\Omega,\mathcal B,P)$ 
having the form $z=x+iy$ where $x$ and $y$ are independent real Gaussian 
random variables with mean zero and equal variances.  A Gaussian 
space is a complex linear subspace $G\subseteq L^2(\Omega,\mathcal B,P)$ 
consisting of Gaussian random variables.  Finally, with any 
subspace $S$ of $L^2(\Omega,\mathcal B,P)$ there is 
an associated sub $\sigma$-algebra $\mathcal B_S$ of 
$\mathcal B$, namely the $\sigma$-algebra generated by 
a sequence of functions in $S$ that has $S$ as its closed 
linear span.  Up to sets of measure zero, $\mathcal B_S$ 
does not depend on the choice of spanning sequence.  

Suppose we are 
given a Gaussian space $G\subseteq L^2(\Omega,\mathcal B,P)$ 
of random variables, along with 
a second probability measure $Q$ 
on $(\Omega,\mathcal B)$ that  is mutually 
absolutely continuous with $P$ and which has 
the further property 
that every random variable in $G$ is also Gaussian 
when viewed as a random variable relative to $Q$.  
We may consider the two 
inner products defined on $G$ by 
$$
\langle z_1,z_2\rangle_{P}=\int_\Omega z_1\bar z_2\,dP,\qquad 
\langle z_1,z_2\rangle_{Q}=\int_\Omega z_1\bar z_2\,dQ.  
$$
Let us write $G_P$ for the Hilbert space structure  of  
$G$ relative to the inner product $\langle\cdot,\cdot\rangle_P$ 
and $G_Q$ for the inner product space associated
with $\langle\cdot,\cdot\rangle_Q$.  
Since $Q\sim P$, it follows that $G$ is also a 
closed subspace of $L^2(\Omega,\mathcal B,Q)$, 
and an application of the closed graph theorem 
implies that 
the identity map of $G$ defines 
an invertible operator from $G_P$ to $G_Q$.  
In particular, the two Hilbert norms on $G$ are equivalent, 
hence there is a unique positive invertible 
operator $B\in\mathcal B(G_P)$ such that 
\begin{equation}\label{tsirGauEq2}
\langle z_1,z_2\rangle_Q=\langle Bz_1,z_2\rangle_P, 
\qquad z_k\in G.  
\end{equation}
The possibilities for $B$ are characterized as follows.  

\begin{theorem}[Feldman, H\'ajek, Segal]\label{tsirGauThm2}
Any operator satisfying (\ref{tsirGauEq2}) has the 
form $B=\mathbf 1+C$  where $C$ is a Hilbert Schmidt
operator in $\mathcal B(G_P)$.  Conversely,  
if $B$ is a positive invertible operator in $\mathcal B(G_P)$ such 
that $\mathbf 1-B$ is Hilbert Schmidt, then there 
is a probability measure $Q$ on 
$(\Omega,\mathcal B)$, 
mutually absolutely continuous with $P$, 
such that $G$ is also 
a Gaussian space relative to $Q$ and for which 
$$
\langle Bz_1,z_2\rangle_P=\int_\Omega z_1\bar z_2\,dQ, 
\qquad z_k\in G.  
$$
This formula determines 
$Q$ uniquely on the $\sigma$-algebra 
$\mathcal B_G$ associated with $G$.  
\end{theorem}

Suppose now that 
$M, N$ are two subspaces 
of a Gaussian space $G\subseteq L^2(\Omega,\mathcal B,P)$, 
with associated triples 
$(\Omega,\mathcal B_M, P_M)$, $(\Omega,\mathcal B_N,P_N)$.   
If $M$ 
and $N$ happen to be orthogonal, then it is a fundamental 
property of Gaussian random variables that 
$$
P(A\cap B)=P(A)P(B),\qquad A\in \mathcal B_M, 
\quad B\in\mathcal B_N 
$$
and this property allows one to 
identify the $L^2$ space 
associated with the sum $M+N$ with the tensor product 
of $L^2$ spaces 
\begin{equation}\label{prodDecompEq}
L^2(\Omega,\mathcal B_{M+N},P_{M+N})\cong 
L^2(\Omega,\mathcal B_M,P_M)\otimes L^2(\Omega,\mathcal B_N,P_N).   
\end{equation}
Indeed, this identification associates a product of functions 
of the form $F_1F_2$ with $F_1$ $\mathcal B_M$-measurable 
and $F_2$ $\mathcal B_N$-measurable, with 
the tensor product $F_1\otimes F_2$.  

If $M$ and $N$ are merely quasiorthogonal, then we have the 
following substitute, which implies that they are independent 
with respect to an equivalent Gaussian measure.

\begin{theorem}\label{tsirGauThm3}
Let $M,N$ be a quasiorthogonal pair of subspaces of 
a Gaussian space $G\subseteq L^2(\Omega,\mathcal B,P)$.  
There is a unique probability measure $Q$ on 
the Borel space $(\Omega,\mathcal B_{M+N})$ 
satisfying the three conditions 
\begin{enumerate}
\item[(i)]
$Q\sim P_{M+N}$, 
\item[(ii)]
$M+N$ is a Gaussian subspace of $L^2(\Omega,\mathcal B_{M+N},Q)$.  
\item[(iii)]
$(\xi,\eta)\in M\oplus N\mapsto \xi+\eta\in M+N
\subseteq L^2(\Omega,\mathcal B_{M+N},Q)$ is an isometry.   
\end{enumerate}
$M$ and $N$ are probabilistically independent Gaussian subspaces 
relative to $Q$, and moreover the restrictions of $Q$ to the 
sub $\sigma$-algebras $\mathcal B_M$ and $\mathcal B_N$ agree 
with $P_M$ and $P_N$ respectively.  
\end{theorem}

\begin{proof}[Sketch of Proof]
Let $L: M\oplus N\to M+N\subseteq L^2(\Omega,\mathcal B,P)$ be 
the natural linear map $L(w_1,w_2)=w_1+w_2$.  Since 
$M$ and $N$ are quasiorthogonal, $L$ 
is an equivalence operator, hence $L^{-1*}L^{-1}=\mathbf 1+C$ where 
$C$ is a Hilbert Schmidt operator on $M+N$.  Theorem 
\ref{tsirGauThm2} implies that there is a unique 
probability measure $Q$ on 
$(\Omega,\mathcal B_{M+N})$ that satisfies (i), (ii), and 
obeys 
$$
\int_{\Omega}z_1\bar{z_2}\,dQ=\langle L^{-1}z_1,L^{-1}z_2\rangle_{M\oplus N}.  
\qquad z_1, z_2\in M+N,   
$$
and the assertion (iii) follows from this formula. 
The last sentence follows from basic properties of Gaussian 
spaces and  
the uniqueness assertion of Theorem \ref{tsirGauThm2}.  
\end{proof}

\begin{remark}
Theorem \ref{tsirGauThm3} gives a precise sense in which 
certain pairs of Gaussian subspaces 
of $L^2(\Omega,\mathcal B,P)$ 
can be ``straightened" by replacing $P$ with another Gaussian 
measure that is equivalent to it.  The same thing 
can be done for any 
finite set of $n$ quasiorthogonal subspaces of a Gaussian space, 
the proof being a straightforward variation of the one above.   
We require only the case $n=2$ for this discussion.  
\end{remark}

Because of these remarks 
we can make an identification 
$$
L^2(\Omega,\mathcal B_{M+N},Q_{M+N})\cong 
L^2(\Omega,\mathcal B_M,Q_M)\otimes L^2(\Omega,\mathcal B_N,Q_N).  
$$
In order to make use of this in the construction to follow it will 
be necessary to make such identifications in a more 
invariant way, in terms of generalizations of measure 
spaces called measure classes.  We 
now describe that procedure.  

\subsubsection{The $L^2$ Space of a Measure Class}
By a {\it measure class} we mean a triple $(X,\mathcal B,\mathcal M)$ 
consisting of a Borel space $(X,\mathcal B)$ together with   
a nonempty set $\mathcal M$ of finite 
positive measures on $(X,\mathcal B)$ with the 
property 
$$
\mu\in \mathcal M, \qquad \nu\sim \mu\implies \nu\in \mathcal M,   
$$
$\mu\sim\nu$ denoting mutual absolute continuity.  
Given two positive finite measures 
$\mu$, $\nu$ on $(X,\mathcal B)$ there is a notion of the geometric 
mean $\sqrt{\mu\nu}$ due to Kakutani:  $\sqrt{\mu\nu}$ is characterized 
as the largest positive measure $\sigma$ with the property 
$$
|\int_X f\bar g\,d\sigma|^2\leq \int_X|f|^2\,d\mu\int_X|g|^2\,d\nu, 
$$
for all bounded measurable functions $f,g$.  It can be defined in 
more concrete terms as the measure $\sqrt{uv}(\mu+\nu)$ where 
$u$, $v$ are the Radon-Nikodym derivatives 
$$
u=\frac{d\mu}{d(\mu+\nu)},\qquad v=\frac{d\nu}{d(\mu+\nu)}.       
$$
The map $\mu,\nu\mapsto\sqrt{\mu\nu}$ has the following  
property: for every set of finite positive measures 
$\mu_1,\dots,\mu_n$ on $(X,\mathcal B)$ and every 
set $f_1,\dots,f_n$ of bounded measurable functions on $X$ 
we have 
$$
\sum_{j,k=1}^n\int_X f_j\bar f_k\,d\sqrt{\mu_j\mu_k}\geq 0, 
$$
(see Chapter 14 of \cite{arvMono} for more detail).  

Fixing a measure 
class $(X,\mathcal B,\mathcal M)$, we form the 
complex vector space $V$ of all formal finite sums 
$$
f_1\sqrt{\mu_1}+\cdots +f_n\sqrt{\mu_n}
$$
where $\mu_1,\dots,\mu_n\in \mathcal M$ and $f_1,\dots,f_n$ 
are bounded measurable functions.  
The preceding inequalities imply 
that we can define a positive semidefinite inner product 
$\langle\cdot,\cdot\rangle$ on $V$ uniquely by setting
$$
\langle f\sqrt{\mu},g\sqrt{\nu}\rangle=\int_X f\bar g\,d\sqrt{\mu\nu},   
$$
where $\mu,\nu\in\mathcal M$ and $f,g$ are bounded functions.  
After dividing out by elements of norm zero and completing, we obtain a 
Hilbert space $L^2(X,\mathcal B,\mathcal M)$.  

There is a natural ``square root" 
map $\mu\in\mathcal M\mapsto \sqrt\mu$ of 
$\mathcal M$ into $L^2(X,\mathcal B,\mathcal M)$ and these 
square roots can be shown to span $L^2(X,\mathcal B,\mathcal M)$.  
There is also a natural $*$-representation $\pi$ of the 
$C^*$-algebra $B(X)$ of all bounded Borel functions on 
$X$ on $L^2(X,\mathcal B,\mathcal M)$, defined uniquely 
by requiring $\pi(f)g\sqrt\mu=fg\sqrt\mu$.  
If $\mathcal M=[\mu_0]$ consists of all finite measures 
$\nu$ that are mutually absolutely continuous with 
respect to a fixed 
finite positive measure $\mu_0$, then it is possible 
to identify $L^2(X, \mathcal B,\mathcal M)$ with 
$L^2(X,\mathcal B,\mu_0)$.  Considering such formal 
expressions as elements of $L^2(X,\mathcal B,\mathcal M)$, 
one finds that for any strictly positive bounded Borel function $f$ and 
a measure $\mu\in\mathcal M$, one has 
$\sqrt{f^2\mu}=f\sqrt\mu$ where $f^2\mu$ denotes the 
obvious measure in $\mathcal M$.  

More significantly, this assignment of a Hilbert space 
to a measure class has the following property: 
Given two Borel spaces $(X,\mathcal A)$ 
and $(Y,\mathcal B)$ and a measure class $\mathcal P$ on 
the cartesian product of Borel spaces 
$(X\times Y,\mathcal A\times\mathcal B)$ which is generated 
as a measure class by a 
finite product measure $\mu\times\nu$, where 
$\mathcal A=[\mu]$ and $\mathcal B=[\nu]$, then 
$L^2(X\times Y,\mathcal A\times\mathcal B, \mathcal P)$ decomposes 
into a tensor product of Hilbert spaces in a 
way analogous to (\ref{prodDecompEq}) above: 
there is a unique unitary operator   
$$
W: L^2(X,\mathcal A,[\mu])\otimes L^2(Y,\mathcal B,[\nu]) 
\to L^2(X\times Y, \mathcal A\times \mathcal B,\mathcal P)
$$
that satisfies 
\begin{equation}\label{newProdEq}
W(f\sqrt{\mu}\otimes g\sqrt{\nu})=f(x) g(y)\sqrt{\mu\times \nu}, 
\qquad (x,y)\in X\times Y
\end{equation}
for all bounded Borel functions $f: X\to\mathbb C$, 
$g:Y\to \mathbb C$.  

The key feature of this 
identification is that while the 
unitary operator $W$ of (\ref{newProdEq}) {\it appears} to 
depend on the particular choice of $\mu$ and $\nu$, 
it actually does not.  Indeed, if we choose other measures 
$\mu_1$ on $X$ and $\nu_1$ on $Y$ such that 
$\mathcal A=[\mu_1]$ and $\mathcal B=[\nu_1]$, 
then $\mathcal P=[\mu_1\times\nu_1]$ and 
one can verify that the operator $W$ of (\ref{newProdEq}) 
also satisfies 
$$
W(f\sqrt{\mu_1}\otimes g\sqrt{\nu_1})=f(x) g(y)\sqrt{\mu_1\times \nu_1}, 
\qquad (x,y)\in X\times Y.  
$$ 
This follows from the Radon-Nikodym theorem, since 
one can check that for every 
measure $\lambda\in\mathcal M$ and 
every bounded nonnegative function $u$, 
the two formal expressions 
$\sqrt{u\lambda}$ and $\sqrt{u(x)}\sqrt{\lambda}$ 
define the same element of $L^2(X,\mathcal A,\mathcal M)$.  

With these preparations, we can now define the 
Tsirelson-Vershik product systems.  Fix $\theta>1$, choose 
a correlation function $C_\theta$ satisfying conditions 
(i), (ii), (iii) above, and let 
$H_\theta$ be the completion of $L^2(\mathbb R)$ in 
the inner product (\ref{innProdEq}).  There is a standard 
construction whereby we can realize 
the Hilbert space $H_\theta$ as a Gaussian space 
$H_\theta\subseteq L^2(\Omega,\mathcal B,P)$.  For every 
interval $I=(a,b)\subseteq \mathbb R$ let 
$H_\theta(I)$ be the closure of the set $L^2(I)$ of functions 
in $L^2(\mathbb R)$ that vanish almost everywhere outside 
$I$, let $\mathcal B_I$ be the corresponding $\sigma$ 
subalgebra of $\mathcal B$ and let $\mathcal P_I$ be the 
set of all measures on $\mathcal B_I$ that are mutually 
absolutely continuous with the restriction of $P$ 
to $\mathcal B_I$.  

Each $(\Omega,\mathcal B_I,\mathcal P_I)$ 
is a measure class, and we can form its canonical 
Hilbert space $L^2(\Omega,\mathcal B_I,\mathcal P_I)$.  
For every $t>0$ $E_\theta(t)$ is defined as the Hilbert space
$$
E_\theta(t)=L^2(\Omega,\mathcal B_{(0,t)},\mathcal P_{(0,t)}),   
$$
and $p:E_\theta\to (0,\infty)$ is the assembled family of Hilbert spaces 
$$
E_\theta=\{(t,\xi): t>0, \xi\in E_\theta(t)\}.  
$$
There is a natural Borel structure on $E_\theta$, and it is 
a nontrivial result of \cite{tsirelSecond} 
that this Borel structure is standard.   

Multiplication is defined 
in $E_\theta$ as follows.  Given $s,t>0$, $\xi\in E_\theta(s)$ and
$\eta\in E_\theta(t)$ we define the product $\xi\cdot \eta\in E_\theta(s+t)$ 
as follows.  The translation operator $T_s\in\mathcal B(L^2(\mathbb R))$ 
restricts to a unitary map of $L^2(0,t)$ onto 
$L^2(s,s+t)$, and thus its closure 
defines a unitary 
operator 
$$
U_{(s,t)}:H_\theta(0,t)\to H_\theta(s,s+t).
$$
It follows from the properties of Gaussian random variables 
that there is a unique unitary operator 
$$
\tilde U_{(s,t)}: L^2(\Omega,\mathcal B_{(0,t)},\mathcal P_{(0,t)})
\to L^2(\Omega,\mathcal B_{(s,s+t)},\mathcal P_{(s,s+t)}); 
$$  
indeed, $\tilde U_{(s,t)}$ 
is implemented by an isomorphism of 
measure classes.  Now up to sets of measure zero, 
$(0,s+t)=(0,s)\cup (s,s+t)$ decomposes into a disjoint 
union of such intervals.  However, the restriction 
of $P$ to 
$\mathcal B_{(0,s+t)}\sim\mathcal B_{(0,s)}\times\mathcal B_{(s,s+t)}$ 
does not decompose 
into a product measure.  
Nevertheless, since the 
subspaces $H_\theta(0,s)$ and $H_\theta(s,s+t)$ are 
quasiorthogonal, Theorems \ref{tsirFirstThm} 
and \ref{tsirGauThm3} imply that the {\it measure class}
$\mathcal P_{(0,s+t)}$ decomposes 
into a product of measure classes 
$\mathcal P_{(0,s)}\times\mathcal P_{(s,s+t)}$.  
Therefore 
we have a natural way of making the 
identification 
$$
L^2(\Omega,\mathcal B_{(0,s+t)},\mathcal P_{(0,s+t)})\sim 
L^2(\Omega,\mathcal B_{(0,s)},\mathcal P_{(0,s)})\otimes 
L^2(\Omega,\mathcal B_{(s,s+t)},\mathcal P_{(s,s+t)})
$$
using a unitary operator $W$ of the form (\ref{newProdEq}).  
At this point the product $\xi\cdot\eta$ can be defined 
as the image of 
$$
\xi\cdot\eta = \xi\otimes \tilde U_{(s,t)}\eta  
$$
under the latter unitary operator.  This is perhaps not the 
best definition of the multiplication in $E_\theta$, but it 
is the quickest.  There are several technical details that 
must be carefully checked to establish 
that $E_\theta$ is a product system, see \cite{tsirelSecond}.

The main two results of \cite{tsirelSecond} are the  following:  

\begin{theorem}[Tsirelson]\label{TsirMainThm}
For every $\theta>1$, the product system 
$E_\theta$ is of type $III$.  $E_{\theta}$ and $E_{\theta^\prime}$ are 
isomophic only if $\theta=\theta^\prime$.  
\end{theorem}

\subsection{Dimension Function}\label{SS:dimFn}

We have defined the numerical index of 
an \esg\ $\alpha$ 
in terms of 
certain structures that are associated with 
its concrete product system $\mathcal E_\alpha$.  
We now describe briefly 
how all of those considerations carry over 
to the setting of abstract product systems.  
In that context, it appears more appropriate 
to think of this numerical invariant as 
a logarithmic dimension function.  

\begin{definition}\label{ctpIndUnitDef}
Let $E$ be a product system.  A unit of $E$ is 
a measurable cross section 
$t\in(0,\infty)\mapsto u(t)\in E(t)$ that 
is multiplicative
\begin{equation}\label{ctpIndMultipEq}
u(s+t)=u(s)u(t), \qquad s,t>0,
\end{equation}
and is not the trivial section $u\equiv 0$.  
\end{definition}

The set of units of $E$ is denoted $\mathcal U_E$. 
The following result provides a 
covariance function for product 
systems that have units; we refer the reader 
to \cite{arvMono} for the proof.  

\begin{prop}\label{ctpIndThm4.2Prop}
Let $E$ be a product system and let 
$u,v\in\mathcal U_E$.  Then there 
is a unique complex number $c_E(u,v)$ 
satisfying 
$$
\langle u(t),v(t)\rangle =e^{tc_E(u,v)},
\qquad t>0.  
$$
The function $c_E: \mathcal U_E\times \mathcal U_E\to \mathbb C$ 
is conditionally positive definite.  
\end{prop}

\begin{definition}\label{ctpIndCovDef}
Let $E$ be a product system for which 
$\mathcal U_E\neq\emptyset$.  The function 
$c_E:\mathcal U_E\times\mathcal U_E\to\mathbb C$ 
is called the covariance function of $E$.  
\end{definition} 

The covariance function is conditionally positive 
definite, thus there is a natural way to use it to 
construct a Hilbert space $H(\mathcal U_E,c_E)$.   
It is known that $H(\mathcal U_E,c_E)$ is 
separable whenever $\mathcal U_E\neq\emptyset$, 
see \cite{arvMono}.  

\begin{definition}\label{ctpIndDimDef}
Let $E$ be a product system.  The dimension 
$\dim E$ of $E$ is defined as the dimension of 
the Hilbert space $H(\mathcal U_E,c_E)$ if 
$E$ has units, and is defined as 
$\dim E=2^{\aleph_0}$ otherwise.  
\end{definition}

Suppose now that we are given an \esg\ 
$\alpha$ acting on $\mathcal B(H)$, with 
concrete product system $\mathcal E_\alpha$.  
Paraphrasing the definition 
of units given in Section \ref{SS:pertCCR}, 
one has semigroups $\{U_t:t>0\}$ that are 
strongly continuous sections 
$t\in(0,\infty)\mapsto U_t\in\mathcal E_\alpha(t)$ 
that also tend strongly to $\mathbf 1$ as $t\to 0+$.  
On the other hand,  
Definition \ref{ctpIndUnitDef} above merely 
requires measurable sections 
of $\mathcal E_\alpha$ defined on 
the open interval $(0,\infty)$ that form a 
nonzero semigroup.  It can 
be shown that the two definitions coincide.  
Moreover, once this 
identification is made, it becomes obvious 
that the covariance 
function defined by Proposition \ref{ctpIndThm4.2Prop}
is identical with the covariance function 
defined in Section \ref{SS:pertCCR}.  Thus the 
two Hilbert spaces are the same, and in particular:
\begin{prop}\label{ctpProdIndEqualsDimProp}
For every \esg\ $\alpha$, we have 
$$
\ind(\alpha)=\dim \mathcal E_\alpha.  
$$
\end{prop}

\begin{remark}[Index of the CCR Flows]
In order to calculate the dimension of a product system, 
one has to calculate its set of units, its covariance function, 
and the dimension of the associated Hilbert space.  In 
order to illustrate the procedure, we show how to 
compute the dimension of the exponential product systems 
$E_N$ and the index of the CCR flows.  

Letting 
$K$ to be an $N$-dimensional Hilbert space, we have 
$$
E_N(t)=e^{L^2((0,t);K)}\subseteq e^{L^2((0,\infty);K)}, 
\qquad t>0, 
$$
and the multiplication of functions $f\in E(s)$ 
and $g\in E(t)$ is defined by 
$$
f\cdot g =f\otimes U_sg\in E(s+t)
$$
where $U_\lambda=\Gamma(S_\lambda)$ is the second-quantized 
shift semigroup of multiplicity $N$, and where 
the tensor product is interpreted in the sense 
of (\ref{basicsCompProdOpDef}).  

Writing $\chi_{(0,t)}\otimes\zeta$ for the function 
in $L^2((0,t);K)$ that has the constant value 
$\zeta\in K$ on the interval $0<x\leq t$ and is 
zero elsewhere, we find that 
$$
\exp(\chi_{(0,s)}\otimes\zeta)\cdot 
\exp(\chi_{(0,t)}\otimes\zeta)
=\exp(\chi_{(0,s+t)}\otimes\zeta),  
$$
hence $u(t)=\exp(\chi_{(0,t)}\otimes\zeta)$, $t>0$, 
defines  a unit of $E_N$.  
Moreover, using 
Proposition \ref{ctpProdCCRidentProp} it is 
possible to show that 
the most general unit of $E_N$ is given by  
\begin{equation}\label{ctpProdUnitEq1}
u^{(a,\zeta)}(t)=e^{ta}\exp(\chi_{(0,t)}\otimes\zeta), 
\qquad t>0, 
\end{equation}
where $a$ is a complex number and $\zeta$ is 
a vector in $K$.  From the 
formula 
$$
\langle\exp(\chi_{(0,t)}\otimes\zeta),
\exp(\chi_{(0,t)}\otimes\omega) \rangle_{E(t)} =
e^{t\langle\zeta,\omega\rangle_K}, \qquad t>0,
$$ 
we find that the covariance function 
of $E_N$ is 
\begin{equation}\label{ctpProdUnitEq2}
c_{E_N}(u^{(a,\zeta)},u^{(b,\omega)}) = 
a+\bar b+\langle \zeta,\omega\rangle_K.  
\end{equation}

A straightforward calculation based on 
(\ref{ctpProdUnitEq2}) 
shows that the Hilbert 
space $H(\mathcal U_{E_N},c_{E_N})$ is 
naturally identified 
with $K$, and therefore $\dim(E_N)=N$.  From 
Propositions \ref{ctpProdIndEqualsDimProp}  and 
\ref{ctpProdCCRidentProp}, we deduce the 
following result, which implies, for example, 
that the CCR flow of rank $2$ cannot be realized 
as a cocycle perturbation of the CCR flow of 
rank $1$.  
\begin{cor}
Let $\alpha$ be the CCR flow of rank $N=1,2,\dots,\infty$.  
Then 
$$
\ind(\alpha)=N.  
$$
\end{cor}
\end{remark}

\subsection{The Classifying Structure $\Sigma$.}
The fundamental problem in this subject is the 
classification of \esg s up to cocycle conjugacy.  The results 
of Section \ref{SS:prodSysCC} imply 
that the problem reduces to the problem of classifying 
product systems up to isomorhism, and therefore 
{\it one should approach the classification problem for 
\esg s by examining the structure of product systems 
on their own terms}.  

We now make this more precise by introducing 
a classifying structure $\Sigma$ for \esg s.  
The elements of $\Sigma$ are isomorphism classes of 
product systems.  The formation of tensor products of 
product systems 
gives rise to a commutative ``addition" in $\Sigma$.   
There is also a natural involution in $\Sigma$  
which makes $\Sigma$ into an involutive abelian semigroup with a 
zero element, and we discuss the significance  
of this involution for dynamics.  
We also describe how $\Sigma$
can be  naturally identified 
with the set of all cocycle conjugacy classes of 
\esg s.  One may conclude that the 
problem of classifying \esg s up to cocycle conjugacy 
reduces to that of determining the 
structure of the involutive semigroup $\Sigma$ and 
discovering computable invariants for its elements.

In order to fully carry out 
the discussion of this section, we must depart from the logical 
development by making use of a key 
result that has not yet 
been discussed, namely that for every 
product system $E$ there is an \esg\ whose concrete 
product system is isomorphic to $E$.  We will discuss 
that result in Lecture 4.  

The {\it trivial product system} is the 
trivial family of one-dimensional Hilbert spaces 
$$
Z=(0,\infty)\times \mathbb C, 
$$
where $\mathbb C$ has its usual inner product 
$\langle z,w\rangle=z\bar w$, where multiplication 
is defined by 
$$
(s,z)(t,w)=(s+t,zw), 
\qquad s,t>0, \quad z,w\in\mathbb C,  
$$
and where the Borel structure on 
$Z=(0,\infty)\times \mathbb C$ is the 
obvious one.  It is significant that 
$Z$ is the only ``line bundle" in the 
category of product systems, as asserted 
by the following result.

\begin{theorem}\label{ctpSigLineBunProp}
Let $E$ be a product system such that $E(t)$ is 
one dimensional for every $t>0$.  Then $E$ is 
isomorphic to the trivial product system 
$Z=(0,\infty)\times\mathbb C$.  
\end{theorem}

{\it Opposite of a Product System.}  We will see momentarily 
that anti-isomorphisms of product systems play a 
significant role in noncommutative dynamics.  By an 
{\it anti-isomorphism} of product systems $\theta: E\to F$ we 
mean a Borel isomorphism that restricts to a unitary 
operator on each fiber $\theta_t: E(t)\to F(t)$, $t>0$, 
such that $\theta(xy)=\theta(y)\theta(x)$ for $x,y\in E$.  

There is a natural involution $E\to E^{\text{op}}$ in the category 
of product systems, defined as follows.  For every 
product system $E$, $E^{\text{op}}$ is defined as the same 
measurable family of Hilbert spaces 
$p: E\to (0,\infty)$, but the multiplication 
in $E^{\text{op}}$ is reversed; for
$x,y\in E$, the product $x\cdot y$ in 
$E^{\text{op}}$ is defined as $yx\in E$.  
$E^{\text{op}}$ is called the {\it opposite} 
product system of $E$.  
If we consider 
the identity map of $E$ as a map of $E$ to $E^{\text{op}}$, 
then it becomes an anti-isomorphism of product 
systems.  Thus, a {\it product system is anti-isomorphic 
to $E$ iff it is isomorphic to $E^{\text{op}}$}.

{\it The tensor product.}  There is also a natural notion of tensor product 
in this category.  
Given product 
systems $E$, $F$, and $t>0$, we can form the 
Hilbert space $E(t)\otimes F(t)$, and the 
associated family of Hilbert spaces 
\begin{equation}\label{ctpSigEq1}
E\otimes F =\{(t,x): t>0,\  x\in E(t)\otimes F(t)\}.  
\end{equation}
Multiplication is defined in $E\otimes F$ in the 
natural way.  In more detail, given elementary 
tensors $x\otimes y\in E(s)\otimes F(s)$ and 
$x^\prime\otimes y^\prime\in E(t)\otimes F(t)$, 
the map
$$
(x\otimes y,x^\prime\otimes y^\prime)
\mapsto xx^\prime\otimes yy^\prime\in E(s+t)\otimes F(s+t)
$$
extends uniquely to a bounded bilinear map of 
$(E(s)\otimes F(s))\times (E(t)\otimes F(t))$ 
into $E(s+t)\otimes F(s+t)$, which in turn can 
be associated with a unitary operator 
$$
(E(s)\otimes F(s))\otimes(E(t)\otimes F(t))\to 
E(s+t)\otimes F(s+t)
$$ 
as required for the multiplication of a product system.  
The Borel structure of $E\otimes F$ has a natural 
definition that we omit.  

It is a nontrivial property of the dimension function 
that it obeys the following logarithmic addition 
formula
$$
\dim(E\otimes F)=\dim E + \dim F
$$
in all cases.  The proof of this formula amounts 
to establishing the fact 
that a unit $w$ of a tensor product $E\otimes F$ 
of product systems {\it must} 
decompose into a tensor product of units
$$
w_t=u_t\otimes v_t,\qquad t>0,
$$
where $u\in\mathcal U_E$ and $v\in\mathcal U_f$ 
(see \cite{arvMono} for more detail).  
In view of the relation between index and dimension
(Proposition \ref{ctpProdIndEqualsDimProp} ), and 
because of Proposition \ref{ctpSigTensProdProp} to follow, 
the preceding 
formula leads to the addition formula for the index 
of \esg s 
$$
\ind (\alpha\otimes \beta)=\ind(\alpha)+\ind(\beta).  
$$

Now let $\alpha$ and $\beta$ be two 
\esg s acting, respectively, 
on $\mathcal B(H)$ and $\mathcal B(K)$, and let 
$\mathcal E_\alpha$, $\mathcal E_\beta$ be their  
concrete product systems.  
For every $t>0$ consider the operator space 
$$
\mathcal E_\alpha(t)\otimes\mathcal E_\beta(t) = 
\overline{\text{span}}^{\|\cdot\|}\{A\otimes B: 
A\in\mathcal E_\alpha(t), B\in \mathcal E_\beta(t)\}, 
$$
the closure being relative to the {\it operator} 
norm.  The total family of spaces 
$$
\mathcal E_\alpha\otimes\mathcal E_\beta=
\{(t,C): t>0,\quad 
C\in\mathcal E_\alpha(t)\otimes\mathcal E_\beta(t)\}
$$
is called the {\it spatial tensor product} of the 
concrete product systems $\mathcal E_\alpha$ and 
$\mathcal E_\beta$.  The tensor product of product  
systems corresponds to the tensor product of 
\esg s as follows: 

\begin{prop}\label{ctpSigTensProdProp}
For any two \esg s $\alpha$, $\beta$, the 
product system of $\alpha\otimes\beta$ is 
the spatial tensor product 
$\mathcal E_\alpha\otimes\mathcal E_\beta$.  
Moreover, the spatial tensor product 
$\mathcal E_\alpha\otimes\mathcal E_\beta$ is 
naturally isomorphic to the tensor product of 
(\ref{ctpSigEq1}).  

\end{prop}

For every product system 
$E$, let $[E]$ denote the class of all product 
systems that are isomorphic to $E$.  The set 
$\Sigma$ of all such equivalence classes can 
be made into an abelian  semgroup by defining 
the sum of classes as follows 
$$
[E] + [F] = [E\otimes F],   
$$
and $\Sigma$ admits 
a natural involution 
$$
[E]^* = [E^{\text{op}}],   
$$
satisfying $(\xi+\eta)^*=\xi^*+\eta^*$, 
$\xi,\eta\in \Sigma$.  It is a straightforward 
exercise to verify that $E\cong E\otimes Z\cong Z\otimes E$
for every product system $E$, 
hence the class of the trivial 
product system $[Z]$ functions as a zero element for 
$\Sigma$.  

The role of $\Sigma$ in the classification problem 
for \esg s is spelled out as follows (where we 
assume the result of Theorem \ref{specExThm1} below).  

\begin{theorem}\label{ctpSigClassifProp}
For every \esg\ $\alpha$, let $[\mathcal E_\alpha]$ 
be the representative of its product system in $\Sigma$.  This 
association 
defines a bijection of the 
set of cocycle conjugacy classes  
of \esg s onto $\Sigma$, 
and one has 
$$
[\mathcal E_{\alpha\otimes\beta}]=
[\mathcal E_\alpha]+[\mathcal E_\beta].
$$   
\end{theorem}

\subsection{Role of the Involution in Dynamics}

The involution of $\Sigma$ is of fundamental importance 
for dynamics, as we now describe.  

Suppose that we are given two
\esg s $\alpha$, $\beta$ acting respectively 
on $\mathcal B(H)$ and $\mathcal B(K)$.  We seek 
conditions on the pair
$\alpha, \beta$ which imply that there 
is a one-parameter group of unitary
operators 
$W = \{W_t: t\in\mathbb R\}$ acting on the tensor product 
$H\otimes K$ whose associated automorphism group 
$\gamma_t(C) = W_tCW_t^*$ satisfies 
\begin{align}
\gamma_t(A\otimes\mathbf 1_K) &= \alpha_t(A)\otimes\mathbf 1_K, 
\quad\text{ for } t\geq 0 \label{ctpSigEq2}\\
\gamma_t(\mathbf 1_H\otimes B) &= \mathbf 1_H\otimes\beta_{-t}(B), 
\quad\text{ for } t\leq 0.\label{ctpSigEq3}
\end{align}
When such a group exists, $\alpha$ and $\beta$ are said
to be {\it paired}.  This relation was introduced by Powers and
Robinson in \cite{powRob} as an intermediate step in their definition 
of another index.  
We will not pursue the Powers-Robinson index here, 
but we do want to emphasize the importance of the pairing 
concept for dynamics.  

Let us first recall the context of 
Lecture 1.  
Considering the von Neumann algebra 
$\mathcal M=\mathcal B(H)\otimes\mathbf 1_K$ 
as a type $I$ subfactor of $\mathcal B(H\otimes K)$, 
with commutant 
$\mathcal M^\prime=\mathbf 1_H\otimes\mathcal B(K)$,   
we are given a pair of \esg s $\alpha$, $\beta$ acting, 
respectively, on $\mathcal M$ and $\mathcal M^\prime$, 
and we are asking if there is 
a one-parameter group of automorphisms 
$\gamma$ of $\mathcal B(H\otimes K)$ that satisfies 
the two conditions of (\ref{intIndEq1.1}).  

The following result implies that a necessry 
and sufficient condition for the existence of 
an automorphism $\gamma$ satisfying 
(\ref{ctpSigEq2}) and (\ref{ctpSigEq3}) 
is that the product 
systems of $\alpha$ and $\beta$ should 
satisfy $[\mathcal E_\beta]=[\mathcal E_\alpha]^*$.  

\begin{theorem}\label{ctpSigPairingProp}
Let $\mathcal M\subseteq \mathcal B(H)$ be a 
type $I$ factor and let 
$\alpha$ and $\beta$ be \esg s acting, 
respectively, on $\mathcal M$ and 
$\mathcal M^\prime$.  The following 
are equivalent.
\begin{enumerate}
\item[(i)]
There is a one-parameter automorphism group 
$\gamma=\{\gamma_t: t\in\mathbb R\}$ acting on 
$\mathcal B(H)$ that satisfies (\ref{ctpSigEq2}) 
and (\ref{ctpSigEq3}).  
\item[(ii)]
The product systems $\mathcal E_\alpha$ and 
$\mathcal E_\beta$ are anti-isomorphic.  
\end{enumerate}

More explicitly, if $U=\{U_t: t\in\mathbb R\}$ 
is a strongly continuous one parameter unitary group 
on $H$ whose automorphism group $\gamma_t(A)=U_tAU_t^*$ 
implements $\alpha$ and $\beta$ as in (i), then 
for every $t>0$ we have 
$U_t^*\mathcal E_\alpha(t) = \mathcal E_\beta(t)$, 
and the map 
$\theta: \mathcal E_\alpha\to \mathcal E_\beta$ defined 
by 
\begin{equation}\label{ctpSigExtraEq1}
\theta((t,T)) = (t,U_t^*T),
\qquad t>0, \quad T\in \mathcal E_\alpha(t)
\end{equation}
is an anti-isomorphism of product systems.  

Conversely, 
every anti-isomorphism  
$\theta: \mathcal E_\alpha\to \mathcal E_\beta$ 
has the form
(\ref{ctpSigExtraEq1}) for a  
unique family $\{U_t: t>0\}$ of unitary
operators.   This family is a strongly 
continuous semigroup tending strongly
to $\mathbf 1$ as $t\to 0+$, 
and its extension to a one-parameter unitary 
group in $\mathcal B(H)$ 
gives rise to an automorphism group $\gamma$ 
satisfying (i) as above.  
\end{theorem}

\begin{proof}[Sketch of Proof]
In order to communicate the flavor of the argument, we prove 
the implication (i)$\implies$ (ii), referring the reader to 
\cite{arvMono} for more detail.

Recalling that every one-parameter group $\gamma$ of 
automorphisms of $\mathcal B(H)$ is implemented 
by a strongly continuous one-parameter unitary 
group $U=\{U_t: t\in\mathbb R\}$ by way of 
$\gamma_t(A)=U_tAU_t^*$ for $t\in\mathbb R$, 
$A\in\mathcal B(H)$, 
it suffices to prove that any such group $\gamma$ 
that satisfies (i) must give rise to a map 
$\theta$ as in (\ref{ctpSigExtraEq1}) that defines 
an anti-isomorphism of product systems.    

Note first that 
$U_t^*\mathcal E_\alpha(t)\subseteq \mathcal M^\prime$
for every $t>0$.  
Indeed, if $A\in \mathcal M$ then for every 
$T\in\mathcal E_\alpha(t)$
we have 
$$
AU_t^*T = U_t^*\gamma_t(A)T = 
U_t^*\alpha_t(A)T = U_t^*TA.  
$$

We claim that $U_t^*\mathcal E_\alpha(t)=\mathcal E_\beta(t)$.  
For the inclusion $\subseteq$, 
choose $T\in \mathcal E_\alpha(t)$.  
The preceding paragraph implies that 
$U_t^*T\in \mathcal M^\prime$, so it remains to show that 
$\beta_t(B)U_t^*T = U_t^*TB$ for every $B\in M^\prime$.  For 
that, write
$$
\beta_t(B)U_t^*T = \gamma_{-t}(B)U_t^*T =
U_t^*BU_tU_t^*T = U_t^*BT = U_t^*TB,
$$
the last equality because $T\in \mathcal M$ 
commutes with $B\in \mathcal M^\prime$.

For the inclusion 
$\mathcal E_\beta(t)\subseteq U_t^*\mathcal E_\alpha(t)$,   
choose $S\in \mathcal E_\beta(t)$ and set $T=U_tS$.  Note that 
$T\in \mathcal M^{\prime\prime}=\mathcal M$ 
because for every $B\in \mathcal M^\prime$ we have 
$$
BT = BU_tS= U_t\gamma_{-t}(B)S=U_t\beta_t(B)S =
U_tSB = TB.  
$$
Moreover, 
$T=U_t S\in \mathcal M$ actually belongs
to $\mathcal E_\alpha(t)$, since for 
$A\in \mathcal M$ 
$$
\alpha_t(A)T=
\gamma_t(A)U_tS=U_tAS=U_tSA=TA,   
$$
hence $S=U_t^*T\in U^*\mathcal E_\alpha(t)$.  

Thus for every $t>0$ we can define a map 
$\theta_t: \mathcal E_\alpha(t)\to \mathcal E_\beta(t)$ by 
$\theta_t(T)=U_t^*T$.  By assembling these maps we get a 
bijective 
Borel-measurable function 
$$
\theta: (t,T)\in \mathcal E_\alpha\to 
(t,U_t^*T)\in \mathcal E_\beta
$$ 
that is a linear isomorphism on 
each fiber.  Each $\theta_t$ 
is actually unitary, since for $T_1,T_2\in \mathcal E_\alpha(t)$
we have 
$$
\<\theta_t(T_1),\theta_t(T_2)\>\mathbf 1=
\theta_t(T_2)^*\theta_t(T_1)=
(U_t^*T_2)^*(U_t^*T_1)=
T_2^*T_1=\<T_1,T_2\>\mathbf 1.  
$$
Finally, since 
for $S\in \mathcal E_\alpha(s)$ and  
$T\in \mathcal E_\alpha(t)$ we have 
$$
\theta_{s+t}(ST) = U_{s+t}^*ST=U_t^*(U_s^*S)T=
U_t^*\theta_s(S)T =
U_t^*T\theta_s(S)=\theta(T)\theta(S),   
$$
it follows that $\theta$ is an anti-isomorphism 
of product systems.  
\end{proof}

We now indicate how a key result from Lecture 2 is 
deduced from Theorem \ref{ctpSigPairingProp}.

\begin{proof}[Proof of Theorem \ref{intIndThm1}]
For every $n=1,2,\dots,\infty$, let $E_n$ be the 
exponential product system of dimension $n$.  
We point out that each $E_n$ is 
anti-isomorphic to itself.  
That follows, for example, from
the classification results for 
type $I$ product systems, since 
the product system opposite to $E_n$ is 
a decomposable product system of the same 
dimension $n$, and therefore isomorphic 
to $E_n$.  Alternately, 
one can simply write down an 
explicit anti-automorphism of $E_n$ using the concrete 
description of it given in 
Proposition \ref{ctpProdCCRidentProp}.

One concludes from these 
remarks that $E_m$ is anti-isomorphic to $E_n$ iff 
$m=n$.  
It follows that the product systems of two cocycle perturbations 
$\alpha$, $\beta$ 
of $CAR/CCR$ flows are anti-isomorphic iff 
$\alpha$ and $\beta$ have 
the same numerical index.  
Thus, Theorem \ref{intIndThm1} is now seen as 
the special case 
of Theorem \ref{ctpSigPairingProp} for cocycle 
perturbations of CAR/CCR flows, in the setting 
in which $\mathcal M=\mathcal B(H)\otimes \mathbf 1_K$ 
and $\mathcal M^\prime = \mathbf 1_H\otimes \mathcal B(K)$.  
\end{proof}

\subsection{Gauge Group.}

Let $\alpha=\{\alpha_t: t\geq0\}$ be an 
\esg\ acting on $\mathcal B(H)$.  A gauge cocycle 
is a cocycle $U=\{U_t:t\geq0\}$ for $\alpha$ with 
the property that the corresponding perturbation 
of $\alpha$ is is the trivial one: 
$$
U_t\alpha_t(A)U_t^*=\alpha_t(A),\qquad t\geq0,\quad A\in\mathcal B(H).  
$$
One sees that the pointwise product $UV=\{U_tV_t: t\geq0\}$ 
of two 
gauge cocycles $U,V$ is a gauge cocycle, as is 
$\{U_t^*:t\geq0\}$.  Thus the gauge cocycles form a group 
$G(\alpha)$ under pointwise multiplication, called the 
{\it gauge group} of $\alpha$.  There is also a natural 
topology on $G(\alpha)$ with respect to which it is a 
Polish topological group, but it will not be necessary to deal 
with topological issues here.  The gauge group reflects the ``internal 
symmetries" of $\alpha$ as we will see presently.

It is quite easy to see that 
if $\beta=\{\beta_t: t\geq0\}$ is another \esg\ that is 
cocycle conjugate to $\alpha$, then $G(\alpha)$ and 
$G(\beta)$ are isomorphic.  Thus the gauge group 
provides a rather subtle cocycle conjugacy invariant 
for \esg s.  The purpose of this section is to point 
out the role of the gauge group in dynamical issues, 
to clarify its status as the group of internal symmetries, 
and to exhibit the structure 
of the gauge groups of type $I$ \esg s in very 
explicit terms.  

Let us first examine the role of the gauge group in 
dynamics.  Given two \esg s $\alpha$, $\beta$ acting 
respectively on $\mathcal B(H)$, $\mathcal B(K)$, we have 
seen that a necessary and sufficient for there to exist 
a one-paramter group of 
automorphisms $\gamma=\{\gamma_t:t\in\mathbb R\}$ satisfying 
(\ref{ctpSigEq2}) and (\ref{ctpSigEq3}) is that the 
product systems of $\alpha$ and $\beta$ should 
by anti-isomorphic.  Assuming that this is the case, 
one would then like to know how to parameterize the 
set of all such groups $\gamma$ in concrete terms.  
The following result, a consequence of Theorem 
\ref{ctpSigPairingProp}, shows that the set of all 
such automorphism groups is parameterized naturally 
by the elements of the gauge group $G(\alpha)$.

\begin{theorem}\label{gaugeParamThm}
Let $\alpha$ and $\beta$ be two \esg s acting, 
respectively, on $\mathcal B(H)$ and $\mathcal B(K)$, 
and assume there is a one-parameter 
group $\gamma^0=\{\gamma^0_t: t\in\mathbb R\}$ of
automorphisms of $\mathcal B(H\otimes K)$ that satisfies 
(\ref{ctpSigEq2}) and (\ref{ctpSigEq3}).  
For every other 
one-parameter group of automorphisms $\gamma$ of 
$\mathcal B(H)$, the following are equivalent:
\begin{enumerate}
\item[(i)] $\gamma$ satisfies (\ref{ctpSigEq2}) 
and (\ref{ctpSigEq3}).  
\item[(ii)] There is a gauge cocycle $U=\{U_t: t\geq 0\}$ 
in $G(\alpha)$ such that the actions of $\gamma$ 
and $\gamma^0$ on 
$\mathcal B(H)$ are related as follows
\begin{equation}\label{basicsGaugeEq3}
\gamma_t(X)=U_t\gamma^0_t(X) U_t^*, \qquad X\in\mathcal B(H), 
\quad t\geq 0.  
\end{equation}
\end{enumerate}
Moreover, for every gauge cocycle $U\in G(\alpha)$, 
formula (\ref{basicsGaugeEq3}) defines a semigroup 
$\{\gamma_t: t\geq 0\}$ of automophisms of $\mathcal B(H)$ 
whose unique extension to a one-parameter automorphism group 
of $\mathcal B(H)$ 
satisfies (\ref{ctpSigEq2}) and (\ref{ctpSigEq3}).  
\end{theorem}

While Theorem \ref{gaugeParamThm} pins down 
the lack of uniqueness that accompanies the 
automorphism groups $\gamma$ that 
solve (\ref{ctpSigEq2}) and (\ref{ctpSigEq3}), it provides 
no insight into the structure 
or even the cardinality of  
gauge group, and does not lead 
to an explicit parameterization of 
the these automorphism groups.  
We now show how to identify the gauge group 
in a more concrete way that reveals the role 
of elements of $G(\alpha)$ 
as internal symmetries.  

For every gauge cocycle $U=\{U_t:t\geq 0\} \in G(\alpha)$ 
let 
$\theta^U:\mathcal E_\alpha\to(0,\infty)\times\mathcal B(H)$ 
be the map defined by  
$$
\theta^U(t,T)=(t,U_tT), \qquad t>0, 
\quad T\in\mathcal E_\alpha(t).  
$$
Note first that 
$\theta^U(\mathcal E_\alpha)\subseteq\mathcal E_\alpha$.  
Indeed, for every $t>0$ and $A\in\mathcal B(H)$ we have 
$$
\alpha_t(A)U_tT=U_t\alpha_t(A)T=U_tTA, 
$$
since $U_t$ commutes with $\alpha_t(A)$.  
Indeed, $\theta^U$ 
is a bijection of $\mathcal E_\alpha$ 
onto itself that is unitary 
on fibers.  It is multiplicative because 
\begin{align*}
\theta^U(s,S)\theta^U(t,T)&=(s,U_sS)(t,U_tT)= 
(s+t,U_sSU_tT)\\
&= (s+t,U_s\alpha_s(U_t)ST)=
(s+t,U_{s+t}ST)
\end{align*}
since $U$ is an $\alpha$-cocycle.  
Obviously, $\theta^U\circ\theta^V=\theta^{UV}$, 
and in fact we have an isomorphism 
of groups:

\begin{theorem}\label{gaugeAutThm}
The mapping 
$U\in G(\alpha)\mapsto \theta^U$ 
defines an isomorphism of the gauge group
$G(\alpha)$ onto the group ${\rm{aut}}\,\mathcal E_\alpha$ of 
all automorphisms of the product system of $\alpha$.  
\end{theorem}

Thus, in order to calculate gauge groups one should 
attempt to compute the automorphism groups of product 
systems.  However, little 
is known about type $III$  product  
systems, and even less is known about their automorphism groups.  
The type $II$ case is somewhat less mysterious, but 
still poorly understood.  
A type $II$ product system $E$ of dimension $n=1,2,\dots,\infty$ 
must contain a type $I$ part $E_I$ of dimension $n$ as a 
subsystem, but 
as yet we lack a clear understanding of how 
to assemble $E$ out of its type $I$ part $E_I$ 
and perhaps other data.

On the other hand, for type $I$ product systems these 
calculations can be carried out in very explicit 
terms.  We now describe this result from \cite{arvI}, 
passing lightly over topological issues.  
For every separable Hilbert space $H$, we define a 
group $G_H$ as follows.  As a set, $G_H$ is the 
cartesian product
\begin{equation}\label{ctpAut8.1}
G_H=\mathbb R\times H\times\mathcal U(H),
\end{equation}
$\mathcal U(H)$ being the unitary group 
of $H$.  The multiplication in $G_H$ is 
defined by 
$$
(\lambda,\xi,U)(\mu,\eta,V)= 
(\lambda+\mu+\omega(\xi,U\eta), \xi+U\eta,UV), 
$$
$\omega$ denoting the symplectic 
form on $H\times H$ given by the 
imaginary part of the inner product 
$\omega(\xi,\eta)=\Im\langle\xi,\eta\rangle$.  
The identity of $G_H$ is $(0,0,\mathbf 1)$, 
and inverses are given by the formula 
$$
(\lambda,\xi,U)^{-1}= 
(-\lambda+\omega(\xi,U\xi),U^{-1}\xi,U^{-1}).  
$$
There is a natural topology on $G_H$ which 
makes it into a Polish topological group; it
becomes a Lie group when $H$ is finite dimensional.

\begin{remark}[$G_H$ and the Canonical Commutation Relations]
\index{canonical commutation relations!and automorphisms}%
The map $\lambda\in\mathbb R\mapsto (\lambda,0,\mathbf 1)$ 
is an isomorphism of $\mathbb R$ onto 
the center of $G_H$, and the map 
$\pi(\lambda,\xi,U)=U$ defines a surjective homomorphism 
of $G_H$ onto $\mathcal U(H)$ whose kernel is 
$K=\mathbb R\times H\times\{\mathbf 1\}$.  Noting  
the multiplication rule in $K$,
\begin{equation}\label{ctpAutEq8.2}
(\lambda,\xi,\mathbf 1)(\mu,\eta,\mathbf 1) =
(\lambda+\mu+\omega(\xi,\eta),\xi+\eta,\mathbf 1), 
\end{equation}
we see that when $H$ is finite dimensional, $K$ is 
isomorphic to the universal covering group of 
the Heisenberg group of 
an appropriate dimension \cite{howeBams}.  
The split exact sequence of groups 
$$
 1\longrightarrow K\longrightarrow 
G_H\longrightarrow \mathcal U(H)\longrightarrow 1
$$
exhibits $G_H$ as a semi direct product 
of $K$ with the unitary group $\mathcal U(H)$.  
\index{Heisenberg group}%
\index{covariance function!automorphism of}%

Perhaps one can see the relation to the 
canonical commutation relations 
in more concrete terms 
by looking at the unitary representations of 
$G_H$.  The most general strongly continuous 
unitary representaation of $G_H$ on a 
Hilbert space $L$ is obtained 
as follows.  Let $(U,W,\Gamma)$ be a triple 
consisting of (i) a strongly continuous 
one-parameter unitary group 
$U=\{U_t: t\in\mathbb R\}$ acting on $L$, 
(ii) a strongly continuous mapping 
$$
\xi\in H\mapsto W(\xi)\in\mathcal U(L)
$$
satisfying the commutation relations 
\begin{equation}\label{ctpAutEq8.3}
W(\xi)W(\eta)=U(\omega(\xi,\eta))W(\xi+\eta),
\end{equation}
and (iii) a strongly continuous representation 
$\Gamma$ of $\mathcal U(H)$ on $L$, all of 
which satisfy the compatibility relations 
\begin{enumerate}
\item[(iv)]\label{ctpAutEq8.4iv}
$U(\lambda)$ commutes with 
$U(\mathbb R)\cup W(H)\cup \Gamma(\mathcal U(H))$ 
for every $\lambda\in\mathbb R$,
\item[(v)]\label{ctpAutEq8.4v}
$\Gamma(U)W(\xi)\Gamma(U)^{-1}=W(U\xi)$, $\quad \xi\in H$, 
$U\in\mathcal U(H)$.  
\end{enumerate}
Given such a triple $(U,W,\Gamma)$, 
one readily verifies that 
\begin{equation}\label{ctpAutEq8.5}
R(\lambda,\xi,V) = U(\lambda)W(\xi)\Gamma(V)
\end{equation}
defines a strongly continuous unitary representation 
$R:G_H\to\mathcal B(L)$.  Conversely, every unitary 
representation $R$ of $G_H$ decomposes uniquely 
into a product 
as in (\ref{ctpAutEq8.5}).  

Notice that if $R$ is an {\it irreducible} representation, 
then necessarily $U(\lambda)$ must be a scalar for every 
$\lambda\in \mathbb R$, hence there is a real 
constant $c$ such that $U(\lambda)=e^{ic\lambda}\mathbf 1$, 
$\lambda\in\mathbb R$.  
In this case, (\ref{ctpAutEq8.3}) reduces to 
\begin{equation*}
W(\xi)W(\eta)=e^{ic\omega(\xi,\eta)}W(\xi+\eta),
\qquad \xi,\eta\in H.  
\end{equation*}
Thus, $W$ is simply a Weyl system  
for the symplectic form on $H\times H$ given 
by $c\cdot\omega(\xi,\eta)$, with the additional 
property that it is equivariant under the action 
of the unitary group of $H$ in the sense 
of (v) above.  We conclude that $G_H$ 
{\it is the appropriate group for the analysis 
of representations 
of the canonical commutation relations 
that include rotational symmetry.}
\end{remark}

\begin{theorem}\label{gaugeGpTypIThm}
For every $N=1,2,\dots,\infty$, the gauge group 
of the CAR/CCR flow $\alpha$ of index $N$ is isomorphic 
to the group $G_H$ of (\ref{ctpAut8.1}), where $H$ is a 
Hilbert space of dimension $N$.  
\end{theorem}

We have not indicated a specific isomorphism of groups 
$\theta: G(\alpha)\to G_H$ in the discussion 
above, but it is not 
hard to do so (see \cite{arvMono}).  
Once this is done, Theorem \ref{gaugeGpTypIThm} 
provides an explicit 
parameterization of the set of anti-isomorphisms 
of Theorem \ref{ctpSigPairingProp}, and therefore 
a parameterization of the set of all one-parameter 
automorphism groups $\gamma$ that satisfy 
(\ref{ctpSigEq2}) and (\ref{ctpSigEq3}) 
in cases where 
both $\alpha$ and $\beta$ are cocycle perturbations 
of CAR/CCR flows of the same index $N$.

\specialsection{Spectrum of an $E_0$-Semigroup.}

Every product system $E$ is associated with an intrinsic 
Hilbert space $L^2(E)$.  There is also a natural way of 
representing $E$ as a concrete product system
$\mathcal E$ acting on $L^2(E)$, and  $\mathcal E$ is 
isomorphic to $E$.  In turn, 
$\mathcal E$ is associated with a semigroup of endomorphisms 
$\alpha=\{\alpha_t: t\geq 0\}$ 
of $\mathcal B(L^2(E))$ by Proposition
\ref{ctpFocRecoverProp}.   
However, this $E$-semigroup is  not 
an \esg\ because the projections $\alpha_t(\mathbf 1)$ 
decrease to $0$ as $t\to\infty$.  

This lecture focuses on the problem of constructing 
an \esg\ whose product system is isomorphic to $E$ by 
exploiting properties of a \cstar\ naturally 
associated with $E$, called the spectral \cstar\  
of $E$.  
We describe the basic theory 
of spectral \cstar s, emphasizing their 
role in the representation theory of product 
systems, and we discuss the issue of simplicity.  
The results were originally obtained in the 
series \cite{arvII}, \cite{arvIII}, \cite{arvIV}; 
and have been substantially rewritten in \cite{arvMono}.

\subsection{The $C^*$-algebra of a Product System.}

We now describe the fundamental 
properties of the regular representation and 
antirepresentation of a product 
system $E$, we define the spectral \cstar\ 
$C^*(E)$, and we discuss the most 
basic properties of these structures.

With every product system $E$ there is a naturally 
associated Banach algebra $L^1(E)$ of integrable sections 
$t\in(0,\infty)\mapsto f(t)\in E(t)$.  The norm on $L^1(E)$ 
is 
$$
\|f\|_1=\int_0^\infty \|f(t)\|\,dt, 
$$
and multiplication is defined by convolution
$$
f*g(t)=\int_0^t f(s)g(t-s)\,ds, \qquad f,g\in L^1(E), 
\quad t>0.  
$$
Notice that $f(s)g(t-s)\in E(t)$ for every 
$s$ satisfying $0< s<t$, so that $f*g$ is 
a well-defined section; one verifies the inequality 
$\|f*g\|\leq \|f\| \|g\|$ by a familiar application 
of the Fubini theorem and the fact that for 
$v\in E(s)$, $w\in E(t-s)$ one has 
$\|vw\|=\|v\|\|w\|$.

A {\it representation} of a product system $E$ 
is a Borel measurable operator valued map $\phi: E\to\mathcal B(H)$ 
satisfying 
\begin{enumerate}
\item[(i)]$\phi(y)^*\phi(x)=\langle x,y\rangle\mathbf 1$, for 
every $x,y\in E(t)$, $t>0$, and 
\item[(ii)]
$\phi(z)\phi(w)=\phi(zw)$, for arbitrary $z,w\in E$.  
\end{enumerate}
Condition (i) implies that the restriction of $\phi$ to 
each fiber $E(t)$ is a linear map.  
Anti-representations of $E$ are defined similarly, 
with (ii) replaced by its opposite $\phi(z)\phi(w)=\phi(wz)$.  
Given a representation 
$\phi$ of $E$ on a Hilbert space $H$, 
then for every 
integrable section $f\in L^1(E)$ the operator 
function $t\in(0,\infty)\mapsto \phi(f(t))\in\mathcal B(H)$
is measurable with respect to the weak opeator 
topology, and a standard application of the 
Riesz Lemma allows one to form the operator integral
\begin{equation}\label{specRegOpIntEq0}
\int_0^\infty \phi(f(t))\,dt, 
\end{equation}
which is an operator of norm at most $\|f\|_1$.  
We abuse notation slightly by writing this operator 
as $\phi(f)$; thus
$$
\phi(f)=\int_0^\infty \phi(f(t))\,dt, \qquad 
f\in L^1(E).  
$$
One verifies that $\phi(f)\phi(g)=\phi(f*g)$, and 
in fact $\phi$ is a contractive representation of 
the Banach algebra $L^1(E)$.  In a similar way, 
an anti-representation of $E$ on $H$ can be integrated 
to give a contractive {\it anti\/}representation of 
$L^1(E)$ on $H$.  

We now exhibit a representation and an anti-representation 
of $E$ on the intrinsic $L^2$ space of $E$.  
Consider the Hilbert space 
$$
L^2(E)=\int_{(0,\infty)}^\oplus E(t)\,dt
$$ 
of square-integrable sections 
$\xi: t\in(0,\infty)\mapsto \xi(t)\in E(t)$, $t>0$, 
with inner product 
$$
\langle \xi,\eta\rangle=\int_0^\infty\langle \xi(t),\eta(t)\rangle\,dt.  
$$
For  each $v\in E(t)$, $t>0$ and $\xi\in L^2(E)$ let 
$v\xi$ be the following function in $L^2(E)$ 
$$
v\xi(x) = 
\begin{cases}
v\xi(x-t), \quad &x>t,\\
0,  &0< x\leq t.  
\end{cases}
$$
One verifies easily that the map 
$\ell: E\to \mathcal B(L^2(E))$ defined by  
$\ell_v\xi=v\xi$ is a representation 
of $E$.  As with any representation  
of $E$, this can be integrated to a 
representation $\ell: L^1(E)\to\mathcal B(L^2(E))$, 
$$
\ell_f\xi=\int_0^\infty f(t)\xi\,dt,  
$$
and of course when $\xi\in L^1(E)\cap L^2(E)$, 
$\ell_f \xi$ is seen to be convolution
$$
\ell_f\xi\,(t)=f*\xi\,(t)=\int_0^t f(s)\xi(t-s)\,ds.  
$$
Similarly, for every $\xi\in L^2(E)$ 
we can define $\xi v\in L^2(E)$ by 
$$
\xi v(x) = 
\begin{cases}
\xi(x-t)v, \quad &x>t,\\
0,  &0< x\leq t.  
\end{cases}
$$
This defines an antirepresentation of $E$ on $L^2(E)$.  
The corresponding antirepresentation of $L^1(E)$  on 
$L^2(E)$ is written 
$f\in L^1(E)\mapsto r_f\in\mathcal B(L^2(E))$, 
$$
r_f\xi=\int_0^\infty \xi\,f(t)\,dt, 
\qquad f\in L^1(E),\quad \xi\in L^2(E); 
$$
and when $\xi\in L^1(E)\cap L^2(E)$ 
we have $r_f\xi=\xi*f$.  

\begin{remark}[Left and Right Semigroups]\label{specReg2semigpRem}
There are {\it two} concrete product systems that 
act naturally on $\mathcal L^2(E)$, one associated 
with the regular representation and the other 
associated with the regular antirepresentation
$$
\mathcal E_\ell(t)=\{\ell_v: v\in E(t)\},\quad 
\mathcal E_r(t)=\{r_v: v\in E(t)\},\qquad t>0.  
$$
Correspondingly, Proposition 
\ref{ctpFocRecoverProp} implies that 
there are two semigroups of endomorphisms $\alpha$, 
$\beta$ associated with these concrete 
product systems.  One obtains a more 
explicit expression for $\alpha_t$ and $\beta_t$
by choosing an arbitrary orthonormal basis 
$\{e_1(t),e_2(t),\dots\}$ for $E(t)$, letting 
$U_1(t), U_2(t),\dots$, 
$V_1(t), V_2(t),\dots$ be the sequences 
of isometries $U_n(t)=\ell_{e_n(t)}$, $V_n(t)=r_{e_n(t)}$, 
and writing 
$$
\alpha_t(A)=\sum_{n=1}^\infty U_n(t)AU_n(t)^*, \quad 
\beta_t(A)=\sum_{n=1}^\infty V_n(t)AV_n(t)^*,\quad 
A\in\mathcal B(L^2(E)).  
$$
Since $U_m(s)$ commutes with $V_n(t)$ for every 
$m,n$ and every $s,t>0$, we have 
$\alpha_s\circ\beta_t=\beta_t\circ\alpha_s$.  

For every $t>0$, let $P_t$ be the projection 
onto the subspace 
$L^2((t,\infty);E)$ of $L^2(E)$ consisting of 
all square summable sections that vanish almost 
everywhere on $(0,t]$.  Then we have 
\begin{equation}\label{specRegRangEq}
\alpha_t(\mathbf 1)=\beta_t(\mathbf 1)=P_t, \qquad t>0.  
\end{equation}
Since the intersection $\cap_t L^2((t,\infty);E)$ is 
the trivial subspace $\{0\}$, it follows that both 
$\alpha_t(\mathbf 1)$ and $\beta_t(\mathbf 1)$ 
decrease to $0$ as $t\to\infty$.  In particular, 
{\it neither of the semigroups $\alpha$, $\beta$ is an \esg}.  
Nevertheless these semigroups, and especially 
$\beta$, play a central role in the analysis 
of the spectral \cstar\ $C^*(E)$ introduced below.  
\end{remark}

Before introducing the spectral \cstar\ of $E$, 
we point out that for every pair of functions 
$f,g\in L^1(E)\cap L^2(E)$, the product $\ell_f^*\ell_g$ decomposes 
into a sum 
\begin{equation}\label{leftStarEq}
\ell_f^*\ell_g = \ell_h + \ell_k^*
\end{equation}
where $h$ and $k$ belong to $L^1(E)\cap L^2(E)$.  One verifies 
this by a direct computation, and while it is possible 
give explicit formulas for $h$ and $k$, we shall not 
do so here.

\begin{remark}[On Morita Equivalence]\label{specRegRem0}
Consider the algebra 
$$
\mathcal A=\{\ell_f: f\in L^1(E)\}\subseteq \mathcal B(L^2(E))
$$
of left convolution operators.  
There are two \cstar s of operators on 
$L^2(E)$ that one might associate with 
$\mathcal A$.  Perhaps the most natural one 
is $C^*(\mathcal A)$, the \cstar\ spanned by 
finite products of elements of $\mathcal A$ 
and their adjoints.  However, $C^*(\mathcal A)$ 
is somewhat more cumbersome than a certain 
subalgebra of it that we now describe, and 
it is the latter subalgebra 
that we choose to work with.  The purpose 
of these remarks is to point out that, since 
the two \cstar s are strongly Morita 
equivalent, 
they become isomorphic after tensoring with 
the compact operators.  In particular,  
they have the same representation theory, 
the same lattice of closed two-sided ideals, 
the same $K$-theory, etc.  
 
We write $[\mathcal S]$  for the norm-closed 
linear span of a set 
$\mathcal S\subseteq \mathcal B(L^2(E))$
of operators.  
Formula (\ref{leftStarEq}) implies that 
$\mathcal A^*\mathcal A$ is contained in 
$[\mathcal A+\mathcal A^*]$, and it follows 
that $[\mathcal A\mathcal A^*]$ is 
a \cstar.  
Actually, $[\mathcal A\mathcal A^*]$ 
is a full hereditary subalgebra 
of $C^*(\mathcal A)$.  To see that, note that 
$$
\mathcal R=
[\mathcal A+\mathcal A\mathcal A^*]
$$
is a right ideal in $C^*(\mathcal A)$ because 
$\mathcal A^*\mathcal A\subseteq [\mathcal A+\mathcal A^*]$,
and it has the properties 
\begin{equation*}
[\mathcal R^*\mathcal R]=C^*(\mathcal A), 
\qquad [\mathcal R\mathcal R^*]=[\mathcal A\mathcal A^*].  
\end{equation*}
Thus 
$[\mathcal A\mathcal A^*]$ is a full hereditary 
subalgebra of $C^*(\mathcal A)$; in particular, 
it is strongly Morita equivalent to $C^*(\mathcal A)$.  
\end{remark}

\begin{definition}\label{specRegCstarDef}
The spectral \cstar\ of $E$ is defined as the 
norm-closed linear span of the 
set of operators 
$\mathcal A\mathcal A^*=\{\ell_f\ell_g^*:f,g\in L^1(E)\}$, 
and is denoted  $C^*(E)$.  
\end{definition}

It is significant that representations 
of $C^*(E)$ can always be obtained 
by integrating representations 
of the simpler structure $E$.  

\begin{theorem}\label{specRegRepThm}
For every nondegenerate representation $\pi$ of $C^*(E)$ 
on a Hilbert space $H$ there is a representation 
$\phi: E\to\mathcal B(H)$ such that 
\begin{equation}\label{specRegRepThmEq1}
\pi(\ell_f\ell_g^*)=\phi(f)\phi(g)^*, \qquad f,g\in L^1(E),  
\end{equation}
where for $f\in L^1(G)$, $\phi(f)$ is defined as 
the integral (\ref{specRegOpIntEq0}).  
\end{theorem}

\begin{cor}\label{specRegRepCor}
For every 
nondegenerate representation $\pi: C^*(E)\to \mathcal B(H)$ 
there is a  concrete product 
system $\mathcal E\subseteq (0,\infty)\times \mathcal B(H)$ 
that is
naturally  associated with $\pi$, and which is isomorphic to $E$.  
\end{cor}

We emphasize that it is not at all apparent at this point 
that every concrete product system $\mathcal E$, that 
is isomorphic to $E$, can be obtained from a 
representation of $C^*(E)$ as in Corollary \ref{specRegRepCor}.  
Fortunately, this is
true;  it is a consequence of 
results on amenability in \cite{arvMono}.  
Because of these amenability results, it is possible 
to realize spectral $C^*$-algebras in the concrete 
form given in Definition \ref{specRegCstarDef}, 
rather than in the more 
intangible way they were originally 
defined in \cite{arvII} using universal properties.

We now discuss a useful 
formula that expresses operators of the 
form $\ell_f\ell_g^*$ in terms of the 
rank-one operator $f\otimes\bar g$, as 
an absolutely convergent weak integral  
\begin{equation}\label{specRegOpIntEq}
\ell_f\ell_g^*=\int_0^\infty \beta_t(f\otimes \bar g)\,dt 
\end{equation}
when $f,g\in L^1(E)\cap L^2(E)$.  Notice that 
the left side involves {\it left} convolution operators, 
while the integral on the right involves the semigroup of 
endomorphisms $\beta$ associated with {\it right} 
multiplication operators.  
Formula (\ref{specRegOpIntEq}) is a consequence 
of the following more precise assertion:

\begin{prop}\label{specRegTraceIntProp}
Consider the natural 
action $\beta_*$ of the $E$-semigroup $\beta$ on 
the Banach space $\mathcal L^1$ of all trace class 
operators on 
$L^2(E)$, defined for $t\geq 0$ by 
$$
\tr (\beta_{*t}(A)B)=\tr (A\beta_t(B)),  \qquad 
A\in\mathcal L^1, \quad B\in\mathcal B(L^2(E)).
$$   
Then for every $A\in\mathcal L^1$ and every pair of 
functions  $f,g\in L^1(E)\cap L^2(E)$,  
$$
\int_0^\infty|\langle\beta_{t*}(A)f,g\rangle|\, dt\leq 
{\rm{trace}}|A|\cdot \|f\|_1\|g\|_1, 
$$ 
and moreover 
\begin{equation*}
\tr (A\ell_f\ell_g^*)=
\int_0^\infty\langle \beta_{*t}(A)f,g\rangle\,dt=
 \int_0^\infty\tr (A\beta_t(f\otimes\bar g))\,dt.  
\end{equation*}
\end{prop}

\begin{remark}[$C^*(E)$ versus $C^*(E^{\text{op}})$]
We have encountered the opposite product system 
$E^{\text{op}}$ in Lecture 3, and we now 
relate it to the current discussion.  The regular 
antirepresentation of $E$ on $L^1(E)$ gives rise to 
a representation of $E^{\text{op}}$ on $L^2(E)$; and 
in fact we can identify 
$L^2(E^{\text{op}})$ with $L^2(E)$ in such a way 
that the regular representation of $E^{\text{op}}$ 
on $L^2(E^{\text{op}})$ is unitarily equivalent to the 
representation of $E^{\text{op}}$ on $L^2(E)$ 
associated with right multiplications.  

Thus, we can identify the spectral 
\cstar\ of $E^{\text{op}}$ with the following 
\cstar\ associated with right convolution operators on $L^2(E)$
$$
C^*(E^{\text{op}})=\overline{\text{span}}\{r_f r_g^*: 
f,g\in L^1(E)\}.  
$$
Since this realizes both 
$C^*(E)$ and $C^*(E^{\text{op}})$ on 
the same Hilbert space, we can look for concrete relations 
between them, and the most basic ones follow.  
\end{remark}

\begin{theorem}\label{specRegCompactThm}
Let $E$ be a nontrivial product system and let 
$\mathcal K$ be the algebra of all compact operators 
on $L^2(E)$.  Then $C^*(E)$ and $C^*(E^{\rm op})$ are 
irreducible \cstar s with the following properties:
\begin{enumerate}
\item[(i)]
$C^*(E)\cap\mathcal K=C^*(E^{\text{op}})\cap\mathcal K=\{0\}$.  
\item[(ii)]
For every $A\in C^*(E)$, 
$B\in \mathcal C^*(E^{\text{op}})$ and $v\in E$, 
$$
[A,r_v]\in\mathcal K,
\quad [B,\ell_v]\in\mathcal K,
$$ 
$[X,Y]$ denoting the commutator bracket $XY-YX$.   
\item[(iii)]
$[C^*(E),C^*(E^{\text{op}})]\subseteq\mathcal K$.  
\end{enumerate}
\end{theorem}

\subsection{Infinitesimal Description of $C^*(E)$.}

Let $E$ be a nontrivial product system, 
fixed throughout this section.  
Corollary \ref{specRegRepCor} implies that, 
in order to study 
\esg s $\alpha$ such that $\mathcal E_\alpha\cong E$, one 
should look closely at the representatation theory 
of the spectral \cstar\ $C^*(E)$.  Since every representation 
of a \cstar\ is a direct sum of cyclic representations and 
since cyclic representations 
are associated with positive linear functionals, 
we are led to examine the state space 
of $C^*(E)$ where by definition, a {\it state} of 
a \cstar\ is a positive linear functional 
of norm $1$.  
In the following sections we describe an 
analysis of the state space of $C^*(E)$, culminating 
in Theorems \ref{specStMainThm} and 
\ref{specStMainThm2}.  The  
results are applied to establish the existence of 
\esg s whose product systems have arbitrary structure 
in Theorem \ref{specExThm1} below.

In order to obtain an effective description of 
the state space of $C^*(E)$, we require an alternate 
description of this \cstar\ in terms of the infinitesimal 
generator of the semigroup 
$\beta=\{\beta_t: t\geq 0\}$ 
of endomorphisms of $\mathcal B(L^2(E))$ 
associated with the {\it right} antirepresentation 
of $E$ on $L^2(E)$, see Remark \ref{specReg2semigpRem}.  
In this section we examine the properties of this 
generator, its inverse, and their relation to the 
spectral $C^*$-algebra.  

Let $\delta$ be the generator of $\beta$.  Thus, the 
domain of $\delta$ is the set $\mathcal D$ of all 
operators $A\in\mathcal B(L^2(E))$ with the property 
that the limit 
\begin{equation}\label{specDaEq0}
\delta(A)=\lim_{t\to0+}t^{-1}(A-\beta_t(A))
\end{equation}
exists in the strong operator topology of $\mathcal B(L^2(E))$.  
This definition of $\delta(A)$ misses the usual definition 
up to a sign, but the definition above will be more 
convenient for our purposes here.  

Since $\beta$ is a semigroup of endomorhisms, $\mathcal D$ 
is a $*$-algebra of operators and $\delta$ is 
a derivation from $\mathcal D$ to $\mathcal B(H)$ 
\begin{equation}\label{specDaEq1}
\delta(AB)=\delta(A)B + A\delta(B), \qquad 
A,B\in \mathcal D.  
\end{equation}
Moreover, since the individual maps $\beta_t$ are 
normal, $\mathcal D$ is strongly dense in $\mathcal B(H)$.

For every interval $I\subseteq (0,\infty)$ there 
is a corresponding subspace $L^2(I;E)\subseteq L^2(E)$,
consisting of all square summable sections that vanish almost 
everywhere off $I$, and we write $P_I$ for 
the projection onto $L^2(I;E)$.  
An operator $A\in \mathcal B(L^2(E))$ is said to 
be {\it supported in } $I$ if 
$A=P_IA=AP_I$.  An operator is said to have 
{\it bounded support} if there is a $t>0$ such 
that it is supported in $(0,t]$.  The set 
$$
\mathcal B_0=\bigcup_{t>0} P_{(0,t]}\mathcal B(L^2(E))P_{(0,t]}
$$ 
of all operators of bounded support is obviously a 
$*$-subalgebra of $\mathcal B(L^2(E))$ which contains 
the $*$-algebra $\mathcal K_0=\mathcal K\cap\mathcal B_0$ 
of all compact operators of bounded support 
as a self adjoint two-sided ideal, and in particular 
it is strongly dense.   

The following result asserts that
$C^*(E)$ can be defined purely 
in terms of the generator of the semigroup
of endomorphisms of $\mathcal B(L^2(E))$ associated 
with {\it right} multiplication operators.  This 
fact is central to the description of the state 
space of $C^*(E)$ that we will describe below.

\begin{theorem}\label{cpecDaAlgThm}
 Let 
$\mathcal A$ be the space of all operators $A$ in 
the domain of $\delta$ such that $\delta(A)$ is a compact 
operator with bounded support.  If $A,B\in\mathcal A$ 
are such that $\delta(A)$ is supported in 
$(0,a]$ and $\delta(B)$ 
is supported in $(0,b]$, then $\delta(AB)$  is
supported in $(0,a+b]$.  

$\mathcal A$ is a $*$-algebra whose norm
closure is $C^*(E)$.  
\end{theorem}

\subsection{Decreasing Weights.}In this section 
we introduce a family of linear functionals 
defined on the $*$-algebra $\mathcal A\subseteq C^*(E)$  of 
Theorem \ref{cpecDaAlgThm}, and we characterize 
those among them that are formally positive.

\begin{definition}\label{specWtDef1}
A locally normal weight is a linear functional 
$$
\omega: \mathcal B_0\to\mathbb C
$$
defined 
on the local algebra 
$\mathcal B_0\subseteq\mathcal B(L^2(E))$ with the 
property that for every $t>0$ the restriction 
of $\omega$ to $P_{(0,t]}\mathcal B(L^2(E))P_{(0,t]}$ 
is a positive normal linear functional.  
\end{definition}
\index{locally normal weight}%
\index{weight!locally normal}%

\begin{remark}\label{specWtRem1}
For example, given a 
normal weight $\tilde\omega: \mathcal B(L^2(E))^+\to[0,+\infty]$
with the property that $\tilde\omega(P_{(0,t]})<\infty$ 
for every $t>0$, the restriction of $\tilde \omega$ 
to the cone of positive operators 
in $\mathcal B_0$ extends uniquely to a linear
functional on $\mathcal B_0$ 
that is a locally normal weight.  On the 
other hand, there are locally normal weights on $\mathcal B_0$ 
that cannot be associated with 
normal weights of $\mathcal B(L^2(E))$ 
(a class of examples is constructed in Appendix A of 
\cite{arvIV}).  Thus one should consider locally normal 
weights as somewhat more general than normal weights. 
\end{remark} 

Notice that the semigroup $\beta=\{\beta_t: t\geq 0\}$ 
associated with the left antirepresentation of $E$ 
on $L^2(E)$ acts naturally on $\mathcal B_0$; 
indeed, if $B$ is a bounded operator supported in the 
interval $(0,b]$ then, 
for every $s\geq 0$, $\beta_s(B)$ is a bounded  
operator supported in the interval 
$(s,s+b]\subseteq (0,s+b]$.  

\begin{definition}\label{specWtdef2}
A locally normal weight $\omega$ is called decreasing 
if for every $B\in\mathcal B_0$ and $t\geq 0$ we have 
$\omega(\beta_t(B^*B))\leq\omega(B^*B)$.  
\end{definition}
To maintain euphony, we refer to such 
objects simply as {\it decreasing weights}.  
\index{decreasing weight}%
\index{weight!decreasing}%
The set of all decreasing weights is a cone of linear 
functionals on $\mathcal B_0$, and it is partially 
ordered by the relation 
$\omega_1\leq \omega_2$ iff $\omega_2-\omega_1$ is 
a decreasing weight.  

Now let $\mathcal A\subseteq C^*(E)$ be the 
$*$-algebra of Theorem \ref{cpecDaAlgThm},  let 
$\delta$ be the generator of $\beta$ as defined 
in (\ref{specDaEq0}), and let $\omega$ be a locally 
normal weight.  Since $\delta(\mathcal A)$ is 
contained in the domain of $\omega$ we can define 
a linear functional $d\omega$ on $\mathcal A$ as 
follows
\begin{equation}\label{specWtEq1}
d\omega(A)=\omega(\delta(A)),\qquad A\in \mathcal A.  
\end{equation}
One may interpret $d\omega$ as the 
derivative of $\omega$ in the direction 
opposite to the 
flow of the semigroup of endomorphisms 
$\beta$.  Typically, both $\omega$ and 
$d\omega$ are unbounded linear functionals 
on their respective domains $\mathcal B_0$ 
and $\mathcal A$.  The following result 
identifies the linear functionals $d\omega$ 
that are formally positive on $\mathcal A$.  

\begin{theorem}\label{specWtThm1}
Let $\omega$ be a locally normal weight and let 
$d\omega:\mathcal A\to\mathbb C$ be the linear functional
of (\ref{specWtEq1}).  The following 
are equivalent.  
\begin{enumerate}
\item[(i)]
$d\omega(A^*A)\geq 0$ for every $A\in\mathcal A$.  
\item[(ii)]
$\omega$ is decreasing.  
\end{enumerate}
\end{theorem}

\subsection{State Space of $C^*(E)$.}

Given that $C^*(E)$ is exhibited as the norm closure of 
the $*$-algebra $\mathcal A$ as in Theorem \ref{cpecDaAlgThm}, 
we now give a concrete description of the cone of 
positive linear functionals on $C^*(E)$ in terms 
of the infinitesimal structure of $\mathcal A$ 
(Theorem \ref{specStMainThm}), 
and we determine which positive linear 
functionals on $C^*(E)$ give rise to \esg s 
(Theorem \ref{specStMainThm2}).  

\begin{remark}[Growth of a Decreasing Weight]\label{specStRem1}
Let $\omega: \mathcal B_0\to \mathbb C$ be a decreasing 
weight.  
For every nondegenerate bounded 
interval $I\subseteq (0,\infty)$ let $L^2(I;E)$ 
be the corresponding subspace of $L^2(E)$, 
consisting of square integrable 
sections that vanish almost everywhere 
off $I$, with projection $P_I: L^2(E)\to L^2(I;E)$.  
We will be concerned with decreasing weights 
that satisfy the growth condition
\begin{equation*}
\sup_I \frac{\omega(P_I)}{|I|}<+\infty, 
\end{equation*}
$|I|$ denoting the length of $I$, the supremum 
being taken over all bounded 
intervals $I\subseteq(0,\infty)$.  
Using the fact that $\omega$ is 
decreasing, it is not hard to show that the supremum 
can be restricted to intervals of the 
form $(0,\epsilon]$ where $\epsilon$ is 
arbitrarily small; consequently 
\begin{equation}\label{specStEq00}
\sup_{I\subseteq(0,\infty)}\frac{\omega(P_I)}{|I|} = 
\limsup_{t\to0+}\frac{\omega(\mathbf 1-\beta_t(\mathbf 1))}{t}.  
\end{equation}
The common value of (\ref{specStEq00}) is a number 
in $[0,+\infty]$, called 
the {\it growth} of $\omega$.  
\end{remark}

After these preparations, one has the 
following description of the state 
space of $C^*(E)$ in terms of decreasing weights.  

\begin{theorem}\label{specStMainThm}
Let $\Omega$ be the partially ordered 
cone of all locally normal decreasing 
weights $\omega$ on $\mathcal B_0$ of finite growth.  
For every 
$\omega\in\Omega$, 
let $d\omega$ be the linear functional defined on 
$\mathcal A$ by 
$$
d\omega(A) = \omega(\delta(A)), \qquad A\in\mathcal A.  
$$
Then $d\omega$ is bounded and extends to a positive 
linear functional on  the norm closure $C^*(E)$ 
of $\mathcal A$.  
The map $\omega\mapsto d\omega$ defines an affine 
order isomorphism of $\Omega$ onto the cone of all 
positive linear functionals on $C^*(E)$, and one has 
\begin{equation}\label{specStThmEq1}
\limsup_{t\to0+}\frac{\omega(\mathbf 1-\beta_t(\mathbf 1)}{t}
\leq 
\|d\omega\|\leq 4\cdot 
\limsup_{t\to0+}\frac{\omega(\mathbf 1-\beta_t(\mathbf 1)}{t}.
\end{equation}
\end{theorem}

\begin{remark}[Essential States]\label{specStRem1.5}
States of $C^*(E)$ 
give rise to semigroups of endomorphisms of $\mathcal B(H)$, and 
now we need to make that correspondence 
quite explicit.   Let $E$ be a product system, 
and let 
$\rho$ be a positive linear functional 
on $C^*(E)$.  We assert that there is a triple 
$(\phi,\xi, H)$ consisting of a representation 
$\phi$ of $E$ on a Hilbert space $H$ and 
a vector $\xi\in H$ with the following 
properties
\begin{align}
\rho(\ell_f\ell_g^*)&=\langle\phi(f)\phi(g)^*\xi,\xi\rangle,
\label{specStEq3.a}\\
H &=\overline{span}\{\phi(f)\phi(g)^*\xi: f,g\in L^1(E)\}, 
\label{specStEq3.b}
\end{align}
where $\phi$ is the associated representation 
of the Banach algebra $L^1(E)$ 
$$
\phi(f)=\int_0^\infty \phi(f(t))\,dt, \qquad f\in L^1(E).  
$$
Indeed, the GNS construction gives rise to a triple
$(\pi,\xi,H)$ consisting of a representation $\pi$ of 
$C^*(E)$ on a Hilbert spaces $H$ and a cyclic vector 
$\xi\in H$, such that $\rho(A)=\langle\pi(A)\xi,\xi\rangle$
for $A$ in $C^*(E)$; 
Theorem \ref{specRegRepThm} provides 
a representation $\phi: E\to\mathcal B(H)$ 
such that $\pi(\ell_f\ell_g^*)=\phi(f)\phi(g)^*$ 
for $f,g\in L^1(E)$, and 
(\ref{specStEq3.a}) and (\ref{specStEq3.b}) follow.  
There is a uniqueness assertion that goes with 
such ``GNS" triples $(\phi,\xi,H)$ for $\rho$, but 
we will not require that.  

In particular, given such a triple $(\phi,\xi, H)$ 
for $\rho$, there is an associated 
$E$-semigroup $\alpha=\{\alpha_t: t\geq 0\}$ 
acting on 
$\mathcal B(H)$ by way of 
$$
\alpha_t(A)=\sum_{n=1}^\infty \phi(e_n(t))A\phi(e_n(t))^*,
\qquad t\geq 0, \quad A\in\mathcal B(H)  
$$
where $\{e_1(t),e_2(t),\dots\}$ is an orthonormal 
basis for $E(t)$.  $\alpha$ will be an 
\esg\ iff $\alpha_t(\mathbf 1)=\mathbf 1$ for 
every $t\geq 0$.  Since $\alpha_t(\mathbf 1)$ 
is the projection onto the subspace of 
$H$ spanned by the ranges 
of the operators in $\phi(E(t))$, this 
will be the case iff 
\begin{equation}\label{specStEq4}
[\phi(E(t))H]=H,\qquad t\geq 0.  
\end{equation}
\end{remark}

\begin{definition}\label{specStEssStatDef}
A positive linear functional $\rho$ on $C^*(E)$ 
is called essential if the 
representation $\phi: E\to\mathcal B(H)$ 
associated with $\rho$ satisfies (\ref{specStEq4}), 
and therefore gives rise to an \esg.  
\end{definition}
\index{state!essential}%
\index{essential state}%

\begin{remark}[Invariant Weights]\label{specStRem2}
Let $\omega:\mathcal B_0\to\mathbb C$ be a 
locally normal weight that is invariant 
under $\beta$ in the sense that 
$\omega(\beta_t(B))=\omega(B)$, 
$t\geq 0$, $B\in\mathcal B_0$.  It is 
obvious that $\omega$ is decreasing, and 
we assert that the growth (\ref{specStEq00}) 
is finite.  
Indeed, we claim 
that when $\omega\neq 0$ 
there is a positive constant $c$ 
such that 
\begin{equation}\label{specStEq5}
\omega(\mathbf 1-\beta_t(\mathbf 1))=ct,\qquad 
t\geq 0.  
\end{equation}
To see that, let $P$ be the spectral measure defined 
on $[0,\infty)$ by the property 
$$
P([a,b])=\beta_a(\mathbf 1)-\beta_b(\mathbf 1), 
\qquad 0\leq a<b<+\infty.  
$$
Than we can define a positive measure 
$\mu$ on the Borel subsets of $[0,\infty)$ 
by $\mu(S)=\omega(P(S))$, $S\subseteq[0,\infty)$.  
$\mu$ is finite on compact sets and positive 
on some nondegenerate intervals 
because $\omega\neq 0$.  Moreover, 
for every interval $I=[a,b]\subseteq [0,\infty)$ 
and every $t\geq 0$ we have 
$$
\mu(I+t)=\omega(P_{I+t})=\omega(\beta_t(P_I))
=\omega(P_I)=\mu(I),   
$$
and it follows that $\mu(S+t)=\mu(S)$ for every 
Borel set $S$ and $t\geq 0$.  Such a measure 
must be a nonzero multiple of Lebesgue measure 
and (\ref{specStEq5}) follows.  
\end{remark}

The key fact is that the states of 
$C^*(E)$ that give rise to \esg s are 
precisely the derivatives of 
$\beta$-invariant weights.  

\begin{theorem}\label{specStMainThm2}
Let $\omega$ be a locally normal weight satisfying 
$$
\omega(\beta_t(B))=\omega(B), \qquad B\in\mathcal B_0, 
\quad t\geq 0. 
$$  
Then $\omega$ belongs to $\Omega$, and 
$d\omega$ is an essential positive linear 
functional on $C^*(E)$.  
Conversely, if $\omega\in\Omega$ is such that 
$d\omega$ is essential, then $\omega\circ\beta_t=\omega$
for every $t\geq 0$.  
\end{theorem}

\subsection{Existence of $E_0$-Semigroups.}

We now indicate how the results of 
the preceding discussion are applied to construct 
essential states of $C^*(E)$.

\begin{prop}\label{specExProp}
Let 
$\beta=\{\beta_t: t\geq 0\}$ be the semigroup of endomorphisms 
of $\mathcal B(L^2(E))$ associated with 
the antirepresentation of $E$ on $L^2(E)$.   
There is a 
normal weight  $\omega$ of $\mathcal B(L^2(E))$ 
with the property that $\omega(\beta_t(B))=\omega(B)$ for every 
$t\geq 0$, $B\in\mathcal B(L^2(E))^+$, 
and which satisfies 
\begin{equation}\label{specExEq1}
\omega(\mathbf 1-\beta_t(\mathbf 1))=t,\qquad t\geq 0.  
\end{equation}
\end{prop}

\begin{proof}[Sketch of Proof]
Consider the single 
endomorphism $\beta_1$.  Since 
$\mathbf 1-\beta_1(\mathbf 1)$ is a nonzero 
projection, we may choose a normal state 
$\nu_0$ on $\mathcal B(L^2(E))$ such 
that $\nu_0(\mathbf 1-\beta_1(\mathbf 1))=1$.  
Let $V$ be any isometry satisfying 
\begin{equation}\label{specExEq2}
\beta_1(A)V=VA,\qquad A\in\mathcal B(L^2(E)), 
\end{equation}
and define a sequence of normal states 
$\nu_1, \nu_2,\dots$ on $\mathcal B(L^2(E))$ by 
$$
\nu_n(A)=\nu_0(V^{*n}AV^n),\qquad A\in\mathcal B(L^2(E)),
\quad n=1,2,\dots.  
$$
Since $\nu_0$ annihilates the projection 
$\beta_1(\mathbf 1)$ we have 
$\nu_0\circ\beta_1=0$; and for 
$n\geq 1$ the commutation relation 
(\ref{specExEq2}) implies 
$\nu_n\circ\beta_1=\nu_{n-1}$.  
Hence 
\begin{equation*}\label{specExEq3}
\nu=\sum_{n=0}^\infty \nu_n
\end{equation*}
defines a normal weight of $\mathcal B(L^2(E))$ 
satisfying $\nu\circ\beta_1=\nu$.  
We can now define a normal weight $\omega$ on 
$\mathcal B(L^2(E))^+$ as follows:
\begin{equation*}\label{specExEq4}
\omega(A)=
\int_0^1\nu(\beta_s(A))\,ds,\qquad 
A\in \mathcal B(L^2(E))^+.     
\end{equation*}

One finds that $\omega$ is 
invariant under the full semigrouop 
$\{\beta_t: t\geq 0\}$, and that 
$\omega(P_I)<\infty$ for every bounded 
interval $I\subseteq (0,\infty)$.  
Remark \ref{specStRem2} implies that there 
is a positive constant $c$ such that 
$\omega(\mathbf 1-\beta_t(\mathbf 1))=ct$ for $t>0$, 
and that is obviously sufficient.  
\end{proof}

The central result on the existence of \esg s now follows:

\begin{theorem}\label{specExThm1}
For every product system $E$, there is an 
\esg\ $\alpha$ such that $\mathcal E_\alpha\cong E$.  
\end{theorem}

\begin{proof}
Proposition \ref{specExProp} implies
that there is a locally normal 
weight $\omega$ on $\mathcal B_0$ satisfying 
$$
\omega\circ\beta_t=\omega, \quad\text{and}\quad
\omega(\mathbf 1-\beta_t(\mathbf 1))=t,\qquad 
t\geq 0.  
$$
In particular, $\omega$ is a decreasing 
weight satisfying the growth requirement for 
membership in $\Omega$, and Theorem \ref{specStMainThm} 
implies that $d\omega$ is a positive linear 
functional on $C^*(E)$.  It must be essential 
by Theorem \ref{specStMainThm2}, so that 
the representation $\phi: E\to\mathcal B(H)$  
associated with $d\omega$ as in Remark
\ref{specStRem1.5} gives 
rise to a concrete product system $\mathcal E$ 
that is isomorphic to $E$ and which is 
the concrete product system of an \esg. 
\end{proof}

\subsection{Simplicity}\label{S:specSimp}

Let $E$ be a product system.  The 
one-parameter unitary 
group $\Gamma=\{\Gamma(\lambda):\lambda\in\mathbb R\}$ 
defined on $L^2(E)$ by 
$$
\Gamma(\lambda)\xi\,(t)=e^{it\lambda}\xi(t), \qquad 
t>0,\quad \xi\in L^2(E)   
$$
implements a one parameter group 
of $*$-automorphisms of the spectral \cstar\ by  way of 
$$
\gamma_\lambda(A)=\Gamma(\lambda) A\Gamma(\lambda)^*, 
\qquad A\in C^*(E), \quad\lambda\in\mathbb R.  
$$

Experience has led us to believe that $C^*(E)$ 
has no nontrivial closed two-sided ideals; however, 
the issue remains unresolved in
general.  What we do know is summarized 
in the following two results from \cite{arvII}.  

\begin{theorem}\label{specSimpThm1}
Let $E$ be an arbitrary product system and 
consider the one parameter group 
$\gamma=\{\gamma_\lambda: \lambda\in\mathbb R\}$ 
of gauge automorphisms of $C^*(E)$.  
Then $C^*(E)$ is $\gamma$-simple in the sense that 
the only closed $\gamma$-invariant ideals in 
$C^*(E)$ are the trivial ones $\{0\}$ and $C^*(E)$.  
\end{theorem}

\begin{theorem}\label{specSimpThm2}
Let $E$ be a product system that is not of 
type $III$.  Then $C^*(E)$ is a simple \cstar.  
\end{theorem}

The following problem remains:
\vskip0.1in
{\bf Problem:} Is the spectral \cstar\ of a type 
$III$ product system simple?

\bibliographystyle{amsalpha}

\begin{thebibliography}{Arv97b}

\bibitem[AK92]{arvKish}
W.~Arveson and A.~Kishimoto, \emph{A note on extensions of {$E_0$}-semigroups},
  Proc. Amer. Math. Soc. \textbf{116} ((1992)), no.~3, 769--774.

\bibitem[AL82]{AFL}
A.~Accardi, L.~Frigerio and J.~T. Lewis, \emph{Quantum stochastic processes},
  Publ. RIMS, Kyoto Univ. \textbf{18} ((1982)), 97--133.

\bibitem[Arv89a]{arvI}
W.~Arveson, \emph{Continuous analogues of {F}ock space}, Memoirs Amer. Math.
  Soc. \textbf{80} ((1989)), no.~3.

\bibitem[Arv89b]{arvIII}
W.~Arveson, \emph{Continuous analogues of {F}ock space {$III$}: singular
  states}, J. Oper. Th. \textbf{22} ((1989)), 165--205.

\bibitem[Arv90a]{arvII}
W.~Arveson, \emph{Continuous analogues of {F}ock space {$II$}: the spectral
  {$C^*$}-algebra}, J. Funct. Anal. \textbf{90} ((1990)), 138--205.

\bibitem[Arv90b]{arvIV}
W.~Arveson, \emph{Continuous analogues of {F}ock space {$IV$}: essential
  states}, Acta Math. \textbf{164} ((1990)), 265--300.

\bibitem[Arv97a]{arvPaths}
W.~Arveson, \emph{Path spaces, continuous tensor products, and
  {$E_0$}-semigroups}, Operator Algebras and Applications (A.~Katavolos, ed.),
  Series C: Math. and Phys. Sci., vol. 495, Kluwer Academic Publishers, (1997),
  pp.~1--112.

\bibitem[Arv97b]{arvPureAbsorbing}
W.~Arveson, \emph{Pure {$E_0$}-semigroups and absorbing states}, Comm. Math.
  Phys. \textbf{187} ((1997)), 19--43.

\bibitem[Arv99]{arvIndDil}
W.~Arveson, \emph{On the index and dilations of completely positive
  semigroups}, Int. J. Math. \textbf{10} ((1999)), no.~7, 791--823.

\bibitem[Arv00]{arvInteractions}
W.~Arveson, \emph{Interactions in noncommutative dynamics}, Comm. Math. Phys.
  \textbf{211} ((2000)), 63--83.

\bibitem[Arv03]{arvMono}
W.~Arveson, \emph{Noncommutative dynamics and ${E}$-semigroups}, Monographs in
  Mathematics, Springer-Verlag, New York, 2003, in press.

\bibitem[AW69]{arakiWoods}
H.~Araki and E.~J. Woods, \emph{Complete boolean algebras of type {$I$}
  factors}, Publ. RIMS (Kyoto University) \textbf{2, Series A} ((1969)), no.~2,
  157--242.

\bibitem[Bha99]{bhatMin}
B.~V.~R. Bhat, \emph{Minimal dilations of quantum dynamical semigroups to
  semigroups of endomorphisms of {$C^*$}-algebras}, J. Ramanujan Math. Soc.
  \textbf{14} ((1999)), no.~2, 109--124.

\bibitem[BP94]{BhatPar}
B.~V.~R. Bhat and K.~R. Parthasarathy, \emph{Kolmogorov's existence theorem for
  {M}arkov processes in {$C^*$}-algebras}, Proc. Indian Acad. Sci. (Math. Sci.)
  \textbf{104} ((1994)), 253--262.

\bibitem[BS00]{BhSkeide}
B.~V.~R. Bhat and M.~Skeide, \emph{Tensor product systems of {H}ilbert modules
  and dilations of completely positive semigroups}, Infinite Dimensional
  Analysis, Quantum Probability and Related Topics \textbf{3} ((2000)),
  519--575.

\bibitem[EL77]{EvLewis}
D.~E. Evans and J.~T. Lewis, \emph{Dilations of irreversible evolutions in
  algebraic quantum theory}, Comm. Dublin Inst. Adv. Studies Series A ((1977)),
  no.~24.

\bibitem[Fel58]{feldmGau}
J.~Feldman, \emph{Equivalence and perpendicularity of {G}aussian processes},
  Pac. J. Math. \textbf{8} ((1958)), 699--708.

\bibitem[GV64]{GelVil}
I.~M. Gelfand and N.~Ya. Vilenkin, \emph{Generalized functions}, vol. 4:
  Applications of Harmonic Analysis, Academic Press, New York, (1964).

\bibitem[H{\'a}j58]{haj}
J.~H{\'a}jek, \emph{On a property of normal distributions of any stochastic
  process}, Czech. Math. J. \textbf{83} ((1958)), no.~8, 610--618.

\bibitem[How80]{howeBams}
R.~Howe, \emph{The role of the {H}eisenberg group in harmonic analysis}, Bull.
  A.M.S. \textbf{3} (1980)), no.~2, 821--843.

\bibitem[Kak48]{kakProd}
S.~Kakutani, \emph{On equivalence of infinite product measures}, Ann. Math.
  \textbf{79} ((1948)), no.~2, 214--224.

\bibitem[K{\"u}m85]{kumMarkovDil}
B.~K{\"u}mmerer, \emph{Markov dilations on {$W^*$}-algebras}, J. Funct. Anal.
  \textbf{63} ((1985)), 139--177.

\bibitem[Pow87]{powTypeIII}
R.~T. Powers, \emph{a non-spatial continuous semigroup of {$*$}-endomorphisms
  of {$\mathcal B(H)$}}, Publ. RIMS (Kyoto University) \textbf{23} ((1987)),
  no.~6, 1054--1069.

\bibitem[Pow99]{powNewEx}
R.~T. Powers, \emph{New examples of continuous spatial semigroups of
  endomorphisms of {$\mathcal B(H)$}}, Int. J. Math. \textbf{10} ((1999)),
  no.~2, 215--288.

\bibitem[PR89]{powRob}
R.~Powers and D.~Robinson, \emph{An index for continuous semigroups of
  {$*$}-endomorphisms of {$\mathcal B(H)$}}, J. Funct. Anal. \textbf{84}
  ((1989)), 85--96.

\bibitem[Sau86]{sauv86}
J.-L. Sauvageot, \emph{Markov quantum semigroups admit covariant
  {$C^*$}-dilations}, Comm. Math. Phys. \textbf{106} ((1986)), 91--103.

\bibitem[Seg58]{segCan}
I.~E. Segal, \emph{Distributions in {H}ilbert space and canonical systems of
  operators}, Trans. Amer. Math. Soc. \textbf{88} ((1958)), 12--41.

\bibitem[SeL97]{selThesis}
D.~SeLegue, \emph{Minimal dilations of cp maps and a {$C^*$}-extension of the
  szeg\"o limit}, Ph.D. thesis, University of California, Berkeley, June
  (1997).

\bibitem[Tsi00a]{tsirelFirst}
B.~Tsirelson, \emph{From random sets to continuous tensor products: answers to
  three questions of {W}. {A}rveson}, preprint arXiv:math.FA/0001070, 12 Jan
  (2000).

\bibitem[Tsi00b]{tsirelSecond}
B.~Tsirelson, \emph{From slightly coloured noises to unitless product systems},
  preprint arXiv:math.FA/0006165 v1, 22 June (2000).

\bibitem[VT98]{verTs}
A.~Vershik and B.~Tsirelson, \emph{Examples of nonlinear continuous tensor
  products of measure spaces and non-{F}ock factorizations}, Rev. Math. Phys.
  \textbf{10} ((1998)), no.~1, 81--145.

\end{thebibliography}

\newcommand{\noopsort}[1]{} \newcommand{\printfirst}[2]{#1}
  \newcommand{\singleletter}[1]{#1} \newcommand{\switchargs}[2]{#2#1}
\providecommand{\bysame}{\leavevmode\hbox to3em{\hrulefill}\thinspace}

\end{document}